\documentclass[12pt, reqno, dvipsnames]{amsart}
\usepackage{graphicx} 
\usepackage{extarrows} 
\usepackage{enumitem} 
\usepackage[utf8]{inputenc}
\usepackage{amssymb}
\usepackage{graphicx}
\usepackage{color}
\usepackage{amsmath}
\usepackage{amsfonts}
\usepackage{geometry}
\usepackage{graphicx}
\usepackage{mathrsfs}
\usepackage{amscd}
\usepackage{epic}
\usepackage{shuffle}[2008/10/27 shuffle product symbol]
\usepackage{xcolor}
\definecolor{ForestGreen}{RGB}{34,139,34}
\usepackage[colorinlistoftodos,textsize=tiny]{todonotes}
\usepackage{xargs} 
\usepackage[all]{xy}
\usepackage{tikz-cd}
\usepackage{comment}
\usetikzlibrary{arrows.meta}

\usepackage[mathcal]{euscript}

\theoremstyle{plain}

\newtheorem{theorem}{Theorem}[section]

\newtheorem{prop}[theorem]{Proposition}

\newtheorem{definition}[theorem]{Definition}

\newtheorem{lemma}[theorem]{Lemma}
\newtheorem{remark}[theorem]{Remark}

\newtheorem{cor}[theorem]{Corollary}

\newtheorem*{thm}{Theorem}

\newcounter{introtheorem}
\renewcommand{\theintrotheorem}{\Alph{introtheorem}}
\newenvironment{introtheorem}[2]{%
  \refstepcounter{introtheorem}%
  \par\addvspace{0.65\baselineskip}%
  \noindent\textbf{Theorem~\theintrotheorem\ (#1) [#2].}%
  \ \itshape\ignorespaces
}{\par\addvspace{0.55\baselineskip}}

\theoremstyle{definition}

\newtheorem{example}[theorem]{Example}

\numberwithin{equation}{section}

\long\def\ignore#1{}

\def\N{\mathbb N}

\def \srestr#1{\mathstrut_{\scriptstyle |}\hbox to
-1.5pt{}\raise-4pt\hbox{$\hskip 1pt\scriptscriptstyle #1$}}

\def\c{\mathfrak{c}}

\def\v{\mathfrak v}

\def\N{\mathbb{N}}

\def\N{{\mathbb N}}

\def\KK{{\mathbb K}}

\def\Id{{\mathbb I}}

\def\B{{\mathcal B}}

\def\hb#1{\hbox{#1}}

\def\hb#1{\hbox{#1}}

\def\dim#1{\hb{dim}(#1)}

\def\L1#1{L^1(#1)}

\def\L#1#2{L^{#1}(#2)}

\newcommand{\Fset}{\mathsf{Set^{\times}}}
 
\newcommand{\Vect}{\mathsf{Vect}}

\newcommand{\Ss}{\mathsf{Sp}} 
\newcommand{\Ssk}{\mathsf{Sp}_{\mathbb{K}}} 
\newcommand{\Spgr}{\mathsf{Sp}^{\rm grft}} 
\newcommand{\Spgrdec}{\mathsf{Sp}^{\rm grft}_{\rm dec}} 
\newcommand{\Spcogr}{\mathsf{Sp}_{\rm cogrft}} 
\newcommand{\Spogr}{\mathsf{Sp}^{\rm ogrft}} 
\newcommand{\Spocogr}{\mathsf{Sp}_{\rm ocogrft}} 

\newcommand{\Moup}{\mathsf{Mon}^{\rm up}} 
\newcommand{\Mogr}{\mathsf{Mon}^{\rm grft}} 

\newcommand{\fO}{\mathcal{O}}
\newcommand{\fU}{\mathcal{U}}
\newcommand{\fH}{\mathcal{H}}
\newcommand{\fF}{\mathcal{F}}
\newcommand{\fG}{\mathcal{G}}

\newcommand{\thh}{\mathbf{h}}

\newcommand{\tp}{\mathbf{p}} 
\newcommand{\tq}{\mathbf{q}}

\newcommand{\taa}{\mathbf{a}}

\newcommand{\wX}{\mathbf{X}}
\newcommand{\wE}{\mathbf{E}} 
\newcommand{\wL}{\mathbf{L}} 
\newcommand{\tPi}{\mathbf{\Pi}} 
\newcommand{\tSig}{\mathbf{\Sigma}} 
\newcommand{\tcG}{\mathbf{G_c}} 

\usepackage[colorinlistoftodos,textsize=tiny]{todonotes}
\newcommandx{\unsure}[2][1=]{\todo[linecolor=red,backgroundcolor=red!25,bordercolor=red,#1]{#2}}
\usepackage[colorlinks=true,final,backref=page,hyperindex, pdfstartview=FitV,,, linkcolor=violet, citecolor=red, urlcolor=cyan]{hyperref}

\newcommand{\Bb}{\mathsf{B}}
\newcommand{\Cc}{\mathsf{C}}
\newcommand{\bB}{\reflectbox{$\mathsf{B}$}}
\newcommand{\cC}{\reflectbox{$\mathsf{C}$}}

\newcommand{\sBb}{\overrightarrow{\mathsf{B}}}
\newcommand{\sCc}{\overrightarrow{\mathsf{C}}}

\newcommand{\bBs}{\overrightarrow{\reflectbox{$\mathsf{B}$}}}

\definecolor{ColA}{HTML}{0061b7} 
\definecolor{ColB}{HTML}{b75600} 
\definecolor{FG}{RGB}{34,139,34} 

\newcommand{\Def}[1]{{\bf \color{NavyBlue}{\em #1}}}

\makeatletter
\def\paragraph{\@startsection{paragraph}{4}%
  \z@\z@{-\fontdimen2\font}%
  {\normalfont\bfseries}}
\makeatother

\begin{document}
\title{Graftings and cograftings for species}
\author[Clavier, Vargas]{Pierre Clavier, Yannic Vargas}
\address[P. Clavier]{IRMA, Université de Strasbourg, 
7 rue René-Descartes, 67084 Strasbourg Cedex, France. \newline
IRIMAS - Département de Mathématiques, Université de Haute Alsace, 18 rue des Frères Lumière
68 093, Mulhouse, France}
\email{pierre.clavier@uha.fr}
\address[Y. Vargas]{CUNEF Universidad, Department of Mathematics, C. de Almansa, 101, 28040 Madrid, Spain.}
\email{yannic.vargas@cunef.edu}
\begin{abstract}
We introduce graftings and cograftings for vector species, modelled
on the elementary operations of adjoining and deleting a root in a
rooted forest.  We relate these maps to operads, up and down
operators, and
obtain nontrivial balanced pairs, admitting pair of creation/annihilation operators with commuting relations.  We establish universal
properties for the species of rooted forests.
Using Fock functors, we recover the classical Connes--Kreimer Hopf
algebra, its noncommutative planar analogue, and Foissy's
$q$-infinitesimal Hopf algebra. 
\end{abstract}

\maketitle{}


\tableofcontents

\section{Introduction}\label{sec:Introduction}

\subsection{Grafting of trees and the Connes-Kreimer Hopf algebra}

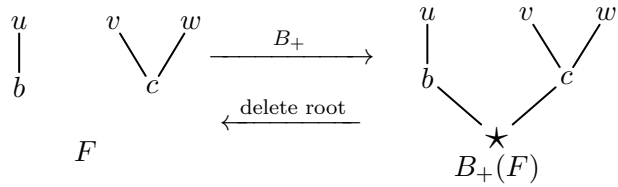
\begin{figure}[htbp]
\centering
\begin{tikzpicture}[
  x=0.80cm,y=0.83cm,
  every node/.style={font=\small,inner sep=1.5pt},
  treeedge/.style={line width=0.75pt,line cap=round}
]
\begin{scope}[shift={(-3.4,0)}]
  \node (b) at (-1.1,0.80) {$b$};
  \node (c) at ( 1.1,0.80) {$c$};
  \node (u) at (-1.1,1.82) {$u$};
  \node (v) at ( 0.46,1.82) {$v$};
  \node (w) at ( 1.74,1.82) {$w$};
  \draw[treeedge] (b)--(u) (c)--(v) (c)--(w);
  \node at (0,-0.22) {$F$};
\end{scope}
\node[inner sep=0pt] at (0,1.45)
  {$\xrightarrow{\qquad B_+\qquad}$};
\node[inner sep=0pt] at (0,0.38)
  {$\xleftarrow{\ \text{delete root}\ }$};
\begin{scope}[shift={(3.4,0)}]
  \node (star) at (0,-0.02) {\Large$\star$};
  \node (bb) at (-1.15,0.94) {$b$};
  \node (cc) at ( 1.15,0.94) {$c$};
  \node (uu) at (-1.15,1.96) {$u$};
  \node (vv) at ( 0.49,1.96) {$v$};
  \node (ww) at ( 1.81,1.96) {$w$};
  \draw[treeedge] (star)--(bb) (star)--(cc)
       (bb)--(uu) (cc)--(vv) (cc)--(ww);
  \node at (0,-0.50) {$B_+(F)$};
\end{scope}
\end{tikzpicture}
\caption{Root grafting and root deletion.  Roots are drawn at the
bottom.  For a planar tree, the children of the new root inherit the
left-to-right order of the forest components.}
\label{fig:intro_grafting}
\end{figure}

A rooted tree is a connected loopless graph with a distinguished vertex, called
its \Def{root}.  A rooted forest is a collection of rooted trees.
One of the simplest ways to construct a tree is to start with a forest,
adjoin a new vertex, and connect this vertex to the roots of all its
components. The new vertex becomes the root of the resulting tree.
We denote this operation, called the \Def{grafting map for trees}, by $B_+$. Deleting the root of a tree reverses
the construction and leaves a forest. Figure~\ref{fig:intro_grafting}
shows this elementary correspondence.

\medbreak

Every rooted forest is a product of rooted trees, and every rooted tree is obtained
uniquely by grafting the roots of a forest below a new vertex.  The Connes--Kreimer Hopf algebra exploits precisely this decomposition,
with a coproduct that records admissible cuts of trees
\cite{connes1999hopf,FoissyInf2009}. More precisely, let $\mathcal{H}_{CK}$ the free algebra generated by non-planar rooted trees. Monomials are naturally associated with (commutative) forests, that is, a set of non-planar rooted trees. This space is endowed with the linear extension of the grafting map $B_+:\mathcal{H}_{CK} \to \mathcal{H}_{CK}$. For every monomial $f=\{t_1, t_2, \hdots, t_r\}$, $B_+(f)$ is the tree obtained by grafting the roots of $t_1, t_2, \hdots, t_r$ on a common new root.  

\

The pair $(\mathcal{H}_{CK}, B_+)$ satisfy the following universal property (see Section 3.1 in \cite{Moerdijk})

\begin{thm}[Kreimer-Moerdijk]
\label{UnivPropHCK}
Let $\mathcal{A}$ be a commutative algebra and $L: \mathcal{A} \to \mathcal{A}$ a linear map. Then there exists a unique algebra map $\phi: \mathcal{H}_{CK} \to \mathcal{A}$ such that $\phi \circ B_+ = L \circ \phi$. Additionally, if $\mathcal{A}$ is a Hopf algebra, then $\phi$ is a map between Hopf algebras.
\end{thm}

The universality of $(\mathcal{H}_{CK}, B_+)$ allows to construct the coproduct for the Connes-Kreimer Hopf algebra as follows: consider the linear map $L: \mathcal{H}_{CK} \otimes \mathcal{H}_{CK} \to \mathcal{H}_{CK} \otimes \mathcal{H}_{CK}$ given by
\begin{equation}
        L(x \otimes y):= B_+(x) \otimes (\varepsilon(y)\cdot 1) + x \otimes B_+(y).
\end{equation}
Hence, there exists a unique algebra map $\Delta: \mathcal{H}_{CK} \to \mathcal{H}_{CK} \otimes \mathcal{H}_{CK}$ such that
\begin{equation}\label{deltaB}
        \Delta \circ B_+ = L \circ \Delta.
\end{equation}

The resulting algebra map $\Delta: \mathcal{H}_{CK} \to \mathcal{H}_{CK} \otimes \mathcal{H}_{CK}$ is the classical coproduct of the Connes-Kreimer Hopf algebra, defined in term of \emph{admissible cuts} on trees and forests.





\medbreak

Our aim is to express this recursive structure at the level of
\Def{combinatorial species} \cite{Joyal,BLL}.  A species $\tp$
assigns a vector space $\tp[I]$ to every finite set $I$, naturally
with respect to relabeling.  Its derivative $\tp'$ records an additional
\Def{ghost} label; hence a map into $\tp'$ can represent the
attachment of a new, distinguished vertex.  Write $\taa$ for the
species of non-planar rooted trees, and $\wE\circ\taa$ for the
species of their forests.  The root-grafting map is then
\[
  \Bb:\wE\circ\taa\to\taa'.
\]

\medbreak

\subsection{The universal property of rooted forests}

Our first result explain why the grafting operation is so well suited to universal
constructions.  Suppose that a commutative monoid $\tp$ in species
has a natural map $u:\tp\to\tp'$ (an \Def{up operator}).  To assign
to each rooted forest an element of $\tp$, there are only two rules
to follow: products of trees must be sent to products in $\tp$, and
adjoining a new root must be sent to the operation $u$.  The
recursive decomposition above makes these rules sufficient to
define the assignment, and also forces its uniqueness.  This is the
starting point of the paper.

\begin{introtheorem}{Universal rooted forests}{Theorem~\ref{Universality}, Corollaries~\ref{coro:univ_grafting} and~\ref{coro:CK_Hopf}}
The free commutative monoid $\wE\circ\taa$, equipped with the up
operator induced by root grafting, is initial among commutative
monoids with an up operator.  Dually, it has the corresponding
universal property for finite-dimensional comonoids with a down
operator.  With the compatible grafting cocycle, it also carries the
canonical Connes--Kreimer Hopf-monoid structure.
\end{introtheorem}

\newpage

\subsection{From admissible cuts to a $q$-deformation of the Connes--Kreimer Hopf algebra}

The universal property in species leads naturally to Hopf structures, in the category of graded vector spaces.  On a
rooted forest $F$, an admissible cut selects edges in such a way that
no path from a root to a leaf meets more than one selected edge.
The forest above the cuts is denoted $P^c(F)$, while the forest
containing the original roots is $R^c(F)$. The ordinary
Connes--Kreimer coproduct is built by summing
$P^c(F)\otimes R^c(F)$ over all cuts:
\[\Delta_{\mathrm{CK}}(F)=
  \sum_{c\in\operatorname{Adm}(F)}P^{c}(F)\otimes R^{c}(F).\]

\begin{figure}[htbp]
\centering
\begin{tikzpicture}[
  x=0.95cm,y=0.95cm,
  vertex/.style={circle,fill=black,inner sep=1.55pt},
  prunedvertex/.style={circle,fill=blue!65!black,inner sep=1.55pt},
  branch/.style={draw=black!85,line width=0.85pt,line cap=round},
  prunedbranch/.style={draw=blue!65!black,line width=0.95pt,line cap=round},
  smoothcut/.style={
    draw=red!75!black,
    line width=1.15pt,
    dash pattern=on 1.3pt off 2.8pt,
    line cap=round,
    line join=round
  }
]

\begin{scope}[shift={(0,3.20)}]

\draw[branch]
  (-4.45,0)--(-5.08,1.08)
  (-4.45,0)--(-3.82,1.08);

\node[vertex] at (-4.45,0) {};
\node[prunedvertex] at (-5.08,1.08) {};
\node[prunedvertex] at (-3.82,1.08) {};
\node at (-4.45,-0.42) {$T_1$};

\draw[branch]
  (-0.75,0)--(-1.70,1.12)
  (-0.75,0)--( 0.20,1.12)
  ( 0.20,1.12)--( 0.75,2.30);

\draw[prunedbranch]
  (-1.70,1.12)--(-2.20,2.30)
  (-1.70,1.12)--(-1.20,2.30)
  ( 0.20,1.12)--(-0.35,2.30);

\node[vertex] at (-0.75,0) {};
\node[prunedvertex] at (-1.70,1.12) {};
\node[vertex] at ( 0.20,1.12) {};

\node[prunedvertex] at (-2.20,2.30) {};
\node[prunedvertex] at (-1.20,2.30) {};
\node[prunedvertex] at (-0.35,2.30) {};
\node[vertex]       at ( 0.75,2.30) {};

\node at (-0.75,-0.42) {$T_2$};

\draw[branch]
  (3.75,0)--(3.75,1.12);

\draw[prunedbranch]
  (3.75,1.12)--(3.08,2.30)
  (3.75,1.12)--(4.42,2.30);

\node[vertex] at (3.75,0) {};
\node[vertex] at (3.75,1.12) {};
\node[prunedvertex] at (3.08,2.30) {};
\node[prunedvertex] at (4.42,2.30) {};

\node at (3.75,-0.42) {$T_3$};

\draw[smoothcut]
  (-6.35,-0.28)
    .. controls (-5.92,0.12) and (-5.55,0.50) .. (-5.12,0.60)
    .. controls (-4.72,0.69) and (-4.22,0.69) .. (-3.85,0.60)
    .. controls (-3.00,0.40) and (-2.05,0.30) .. (-1.25,0.58)
    .. controls (-0.92,0.70) and (-0.58,1.12) .. (-0.10,1.78)
    .. controls (0.18,2.16) and (0.55,2.56) .. (1.15,2.58)
    .. controls (2.00,2.60) and (2.95,1.62) .. (3.75,1.76)
    .. controls (4.42,1.88) and (5.18,2.24) .. (6.28,2.72);

\node[red!75!black,fill=white,inner sep=1pt] at (1.95,2.45) {$c$};

\end{scope}

\begin{scope}[shift={(0,-0.78)}]

\node[prunedvertex] at (-5.18,0.55) {};
\node[prunedvertex] at (-4.55,0.55) {};

\draw[prunedbranch]
  (-2.55,0.55)--(-3.05,1.55)
  (-2.55,0.55)--(-2.05,1.55);

\node[prunedvertex] at (-2.55,0.55) {};
\node[prunedvertex] at (-3.05,1.55) {};
\node[prunedvertex] at (-2.05,1.55) {};

\node[prunedvertex] at (-0.78,0.55) {};

\node[prunedvertex] at (0.55,0.55) {};
\node[prunedvertex] at (1.25,0.55) {};

\node[blue!65!black] at (-1.90,-0.32) {$P^c(F)$};

\end{scope}

\node at (2.10,-0.18) {$\otimes$};

\begin{scope}[shift={(0,-0.78)}]

\node[vertex] at (3.05,0.55) {};

\draw[branch]
  (4.70,0.55)--(4.70,1.55)
  (4.70,1.55)--(4.70,2.55);

\node[vertex] at (4.70,0.55) {};
\node[vertex] at (4.70,1.55) {};
\node[vertex] at (4.70,2.55) {};

\draw[branch]
  (7.05,0.55)--(7.05,1.55);

\node[vertex] at (7.05,0.55) {};
\node[vertex] at (7.05,1.55) {};

\node at (5.15,-0.32) {$R^c(F)$};

\end{scope}

\end{tikzpicture}

\caption{
A nontrivial admissible cut of a rooted forest with three trees.
The cut separates a forest $P^c(F)$ consisting of the two leaves of $T_1$,
the subtree rooted at the left child of the root of $T_2$, one further leaf of $T_2$,
and the two leaves of $T_3$.
The remaining forest $R^c(F)$ contains the roots of $T_1$ and $T_3$,
and the branch of $T_2$ containing its right-most leaf.
}
\label{fig:CK-admissible-cut-corrected}
\end{figure}
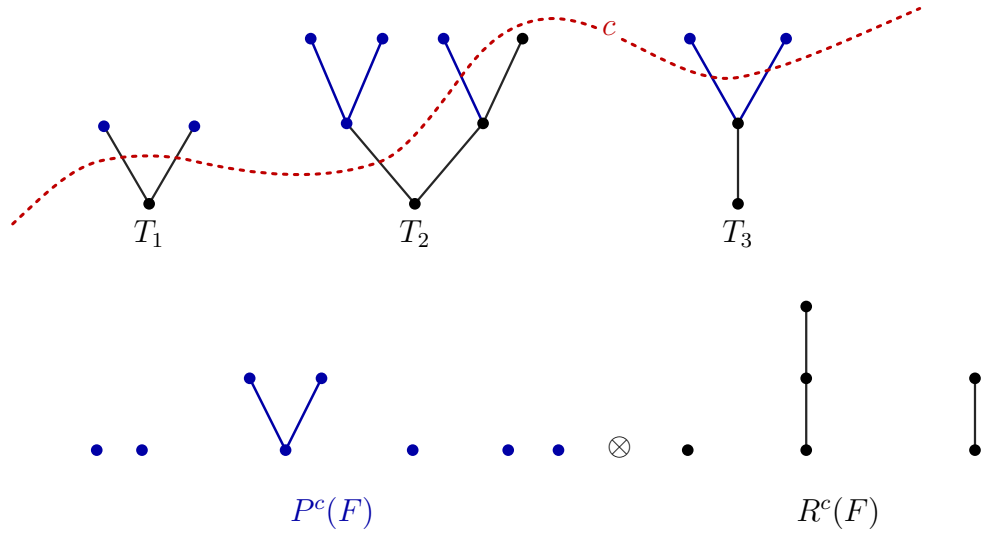

\medbreak

There are three related, but distinct, realizations.  Non-planar
forests yield the commutative Connes--Kreimer Hopf algebra.  If the
children of each vertex are ordered and the forest product is
concatenation, one obtains the noncommutative Connes--Kreimer algebra.
In the other hand, Foissy's infinitesimal version has the same underlying planar forests
but satisfies a different product--coproduct compatibility, the so-called \Def{infinitesimal compatibility} (in the sense of Loday-Ronco, \cite{LRCofree}):
\[
\Delta_0(FG)
=(F\otimes1)\Delta_0(G)+\Delta_0(F)(1\otimes G)-F\otimes G.
\]

\medbreak

A $q$-braided compatibility interpolates between the ordinary and
infinitesimal laws.  The deformation itself is known from Foissy's
work. Our contribution here is to explain it through the species-level
universal grafting construction and the deformed Fock functor
\cite{FoissyQ2003,AM2010}.  This recovers both endpoints without
changing the basic operation $B_+$.

\begin{introtheorem}{The $q$-Connes--Kreimer interpolation}{Theorems~\ref{thm:qCK_Fock_interpolation}}
The deformed Fock construction on planar rooted forests yields the
weighted-cut coproduct
\[
 \Delta_q(F)=\sum_{c\in\operatorname{Adm}(F)}
 q^{\operatorname{inv}_F(c)}P^c(F)\otimes R^c(F),
\]
compatible with the $q$-braided tensor product and the root-grafting
cocycle.  Its specializations at $q=1$ and $q=0$ give, respectively,
the noncommutative Connes--Kreimer and Foissy infinitesimal
coproducts.
\end{introtheorem}

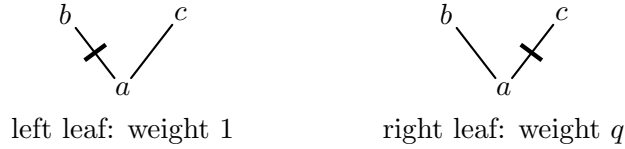
\begin{figure}[htbp]
\centering
\begin{tikzpicture}[
  x=0.9cm,y=0.9cm,
  every node/.style={font=\small,inner sep=1.5pt},
  edge/.style={line width=0.8pt,line cap=round}
]
\begin{scope}[shift={(-2.8,0)}]
 \node (aa) at (0,0) {$a$};
 \node (bb) at (-0.85,1.10) {$b$};
 \node (cc) at ( 0.85,1.10) {$c$};
 \draw[edge] (aa)--(bb) (aa)--(cc);
 \draw[line width=1.65pt] (-0.56,0.42)--(-0.25,0.65);
 \node at (0,-0.62) {left leaf: weight $1$};
\end{scope}
\begin{scope}[shift={(2.8,0)}]
 \node (aaa) at (0,0) {$a$};
 \node (bbb) at (-0.85,1.10) {$b$};
 \node (ccc) at ( 0.85,1.10) {$c$};
 \draw[edge] (aaa)--(bbb) (aaa)--(ccc);
 \draw[line width=1.65pt] (0.25,0.65)--(0.56,0.42);
 \node at (0,-0.62) {right leaf: weight $q$};
\end{scope}
\end{tikzpicture}
\caption{Two single-edge cuts of the planar corolla
$B_+(\bullet\bullet)$. The short stroke marks the cut edge; its
left-to-right position controls the inversion exponent. Together,
these cuts contribute $(1+q)\,\bullet\otimes B_+(\bullet)$ to the
$q$-deformed coproduct.}
\label{fig:intro_qcuts}
\end{figure}

For example, the two single-edge cuts of a planar corolla contribute
with weights $1$ and $q$ as shown in Figure~\ref{fig:intro_qcuts}.
The coefficient $1+q$ is therefore already visible in the smallest
nontrivial planar example.

\subsection{The categorical version of graftings}

The preceding examples begin with a particular species of trees.
We now ask which aspects survive if rooted forests are replaced by
other combinatorial structures.  For a species $\tp$, a \Def{grafting} is a natural map
\[\Bb:\wE\circ\tp\to\tp'.\] 
It combines a family of $\tp$-objects,
on pairwise disjoint label sets, into a new object with an additional
ghost point.  An up operator only acts on one $\tp$-structure,
so the two notions are not literally the same.

\medbreak

When $\tp$ is already a monoid with product $\mu$, however, an up
operator $u:\tp\to\tp'$ produces a grafting by first multiplying the
inputs and then applying $u$:
\[
\Bb_u(x_1\cdots x_r)
=u\bigl(\mu(x_1,\ldots,x_r)\bigr).
\]
This factorization is exactly what identifies the appropriate
subcategory of graftings and leads to an adjunction for the free
commutative monoid construction.

\begin{introtheorem}{Categories of grafting maps and up operators}{Theorem~\ref{thm:comparison_up_grafting_NEW}}
The category of commutative monoids with up operators is isomorphic
to the subcategory of monoids with graftings that factor through the
monoid product, with the prescribed decomposable morphisms.  The free
commutative monoid construction moreover gives the adjunction between
the corresponding categories of up operators and graftings.
\end{introtheorem}

This result, proved in Section~\ref{subsec:up_down_cat}, explains
why one can sometimes replace a grafting map by an up operator
without losing any information, and also identifies when the
larger, partition-sensitive notion of a grafting is necessary.

\medbreak

\subsection{From operads to grafting operations}

There is another, conceptually different, source of graftings.
An operad organizes the substitution of operations into one another.
For rooted trees, one can think of its composition as replacing a
vertex by another rooted tree while respecting the attachment of
subtrees.  The NAP operad is the basic example that produces ordinary
root grafting.  More generally, operadic composition endows the
derivative of an operad with a monoid structure, so an up operator can
be extended from one input to an entire forest of inputs.

\begin{introtheorem}{Grafting from an operad}{Theorem~\ref{thm:operad_grafting}, Proposition~\ref{prop:grafting_from_NAP}}
Let $\tq$ be a positive symmetric operad and let
$u:\tq\longrightarrow\tq'$ be an up operator.  The monoid
structure induced by the operad on $\tq'$ extends $u$ uniquely to a
monoidal grafting $\Bb_u:\wE\circ\tq\longrightarrow\tq'$.  For
the NAP operad, this construction recovers the usual grafting of
rooted forests.
\end{introtheorem}

\medbreak

\subsection{Up/down operators and balanced pairs}

Grafting creates a new vertex; a cografting removes a distinguished
vertex and records the components left behind.  On species these
become \Def{up and down operators},
\[
  u:\tp\to\tp',\qquad
  d:\tp'\to\tp.
\]
A particularly useful situation occurs when their successive
applications satisfy the balanced identities of Gu\c{t}\v{a}--Maassen~\cite{gucta2002symmetric}, later developed by
Aguiar--Mahajan~\cite[Section~19.4]{AM2010}.  These identities can
be interpreted as the combinatorial counterparts of the relations
between creation and annihilation operators.  Rather than checking
them only after passing to a Fock space, we express them directly in
terms of grafting and cografting.

\begin{introtheorem}{Balanced up/down pairs}{Proposition~\ref{prop:balanced_grafting_cografting}, Examples~\ref{ex:balanced_graphs} and~\ref{ex:balanced_rooted_trees}}
The up and down operators induced by a grafting and a cografting
form a balanced pair precisely when four explicit species-level
compatibility relations are satisfied.  Rooted trees and graphs
provide nontrivial examples, with balance constants $\lambda_n=1$
for the constructions considered here.
\end{introtheorem}

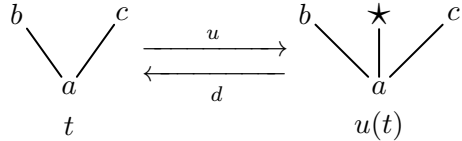
\begin{figure}[htbp]
\centering
\begin{tikzpicture}[
  x=0.82cm,y=0.87cm,
  every node/.style={font=\small,inner sep=1.5pt},
  treeedge/.style={line width=0.75pt,line cap=round}
]
\begin{scope}[shift={(-3.0,0)}]
  \node (a1) at (0,0) {$a$};
  \node (b1) at (-0.85,1.12) {$b$};
  \node (c1) at ( 0.85,1.12) {$c$};
  \draw[treeedge] (a1)--(b1) (a1)--(c1);
  \node at (0,-0.58) {$t$};
\end{scope}

\node[inner sep=0pt] at (-0.65,0.70)
  {$\xrightarrow[]{\qquad u\qquad}$};
\node[inner sep=0pt] at (-0.65,0.01)
  {$\xleftarrow[\qquad d \qquad]{}$};
\begin{scope}[shift={(2.0,0)}]
  \node (a2) at (0,0) {$a$};
  \node (b2) at (-1.20,1.12) {$b$};
  \node (s2) at (0,1.12) {\Large$\star$};
  \node (c2) at ( 1.20,1.12) {$c$};
  \draw[treeedge] (a2)--(b2) (a2)--(s2) (a2)--(c2);
  \node at (0,-0.58) {$u(t)$};
\end{scope}
\end{tikzpicture}
\caption{A rooted-tree up/down pair: $u$ adds the ghost $\star$
as a child of the root, and $d$ removes it and glue all possible children of the ghost to the root of $t$. }
\label{fig:intro_balanced}
\end{figure}

The simplest tree illustration is to attach a ghost leaf to the root
and then remove it.  This makes the relation $d\,u=\mathrm{id}$
visible on each tree, while the remaining balanced identities govern
the possible orders of creating and removing ghost vertices.
Figure~\ref{fig:intro_balanced} displays this operation. The pair $(u,d)$ endows the species of (non-planar) trees with a balanced pair, in the sense of Aguiar-Mahajan. This construction on trees seems to be new.

\medbreak

\subsection{Ordered branching}

The noncommutative counterpart of a forest is an \emph{ordered}
forest.  For a positive species $\tq$, its species is the free
associative monoid $\wL\circ\tq$, and an ordered grafting has the
form $\overrightarrow{\Bb}:\wL\circ\tq\to\tq'$.
There are two different ways to split an ordered forest: unshuffle
places arbitrary components into different ordered subwords, while
deconcatenation splits only at consecutive positions.  The former
leads to the operator studied in
Section~\ref{subsec:ordered_grafting_unshuffle}.

For the latter, the twice-iterated deconcatenation gives the natural
composite
\[
\begin{gathered}
\wL\circ\tq \xrightarrow{\ \Delta_{\mathrm{dec}}^{(2)}\ }
(\wL\circ\tq)^{\cdot3}
\xrightarrow{\ \mathrm{id}\cdot\overrightarrow{\Bb}\cdot\mathrm{id}\ }
(\wL\circ\tq)\cdot\tq'\cdot(\wL\circ\tq)
\cong\ (\wL\circ\tq)'.
\end{gathered}
\]

\medbreak

\subsection{Organization of the paper}

Section~\ref{sec:species} recalls the language of species, their
monoidal structures, and their derivatives.  Graftings and
cograftings, operadic constructions, and examples are introduced in
Section~\ref{sec:graft_cograft}.  Section~\ref{sec:up_down} treats
up/down operators, their categorical comparison, and balanced pairs.
Section~\ref{sec:trees} establishes the universal properties of
rooted forests.  The Connes--Kreimer and Foissy Hopf algebras and
their $q$-interpolation are recovered in
Section~\ref{sec:recovering_CK_Foissy}. 

\subsection{Acknowledgements}
This project started thanks to 
Prof.~Sylvie Paycha, who organized a research stay with the authors
at the Perimeter Institute. Both of us had the opportunity to work
with her as postdoctoral researchers, an experience that significantly
shaped our mathematical trajectories. We warmly thank her for her
unconditional support. Merci and Gracias, Sylvie!

\section{Elements of the theory of species} \label{sec:species}

In this section we will give the preliminaries of combinatorial species and some of its algebraic siblings. We will follow closely \cite{AM2010}.

\subsection{Preliminaries}\label{sec:Preliminaries}

Let $I$ be a finite set. A \Def{set composition} (or \emph{ordered partition}) of $I$ is a finite sequence $F$ denoted by
\[F=(B_1, \hdots, B_k) \qquad \text{ or } \qquad F=B_1|\cdots | B_k,\]
such that $B_1, \hdots, B_k$ are non-empty, pairwise disjoint and $\bigcup_{i=1}^kB_i=I$. The subsets $B_i$ are called the \Def{blocks} of the set composition $F$. If $F$ is a set composition of I , we write $F \vDash I$. Also, $\Sigma[I]$ denotes the set of all set compositions of $I$. The \Def{lenght} of a set composition $F=B_1|\cdots | B_k$ is given by $\ell(F):=k$, the number of blocks of $F$. A \Def{set decomposition} (or simply \Def{decomposition}) is a set composition where empty blocks are allowed. If $F=B_1|\cdots|B_k$ is a decomposition of $I$, we write $I=B_1 \sqcup \cdots \sqcup B_k$.

\medbreak

If $S \subseteq I$, the \Def{restriction} of a set composition $F \vDash I$ to $S$ is the new set composition
\[F|_S:=(B_1 \cap S, \hdots, B_k \cap S)^{\times} \vDash S,\]
where the symbol $\times$ indicates that empty blocks are removed from the list. 

\medbreak

Given two disjoint finite sets $S$ and $T$, the \Def{concatenation} of two set compositions $F=B_1|\cdots|B_r \vDash S$ and $F'=B_1'|\cdots| B_s' \vDash T$ is the set composition of $S\sqcup T$ given by
\[F\cdot G:= B_1|\cdots | B_r|B_1'|\cdots | B_s' \vDash S \sqcup T.\]

A \Def{set partition} $\pi$ of a finite set $I$ is a set of disjoint nonempty subsets of $I$, denoted by
\[\pi =\{B_1, \hdots, B_k\} \qquad \text{ or } \qquad \pi = B_1.\hdots . B_k,\]
such that $\bigcup_{i=1}^rB_i=I$. The \Def{length} of $\pi$ is the number of blocks $\ell(\pi):=k$ of $\pi$. The set of all set partition of $I$ is denoted by $\Pi[I]$. Also, if $\pi \in \Pi[I]$, we write $\pi \vdash I$. Given two disjoint finite sets $S$ and $T$, the \Def{union} $\pi \cup \pi'$ of two set partitions $\pi \vdash S$ and $\pi' \vdash T$ is the set partition of $S\sqcup T$ whose blocks belong either to $\pi$ or $\pi'$. 

\medbreak

Given $\pi \in \Pi[I]$, the \Def{restriction} of $\pi$ to a subset $S \subseteq I$ is the new set partition
\[\pi|_S:=\{B_1\cap S. \hdots, B_k\cap S\}^{\times} \in \Pi[S],\]
where we use the same notation than for set compositions. We say that that the subset $S \subseteq I$ is \Def{$\pi$-admissible} if $\pi=\pi|_S \cup \pi|_T$, where $T$ is the complement of $S$ in $I$.

\medbreak

\medbreak

Given a finite set $I$, the \Def{support map}
\[\text{supp}_I: \Sigma[I]\to \Pi[I]\]
is the function that forget the order of the blocks of a set composition of $I$. It is an order preserving map for the refinement order, which do not preserve meets: $13|24 \wedge 24|13 = 1|2|3|4$ and $\text{supp}_{[4]}(1|2|3|4)=1.2.3.4$, but $\text{supp}_{[4]}(13|24)= \text{supp}_{[4]}(24|13)= 13.24$, so $\text{supp}_{[4]}(13|24)\wedge \text{supp}_{[4]}(24|13)=13.24\neq 1234$.



\subsection{Categories of species: definitions, structures, examples}

In this section we recall the basic definitions of the general theory of \emph{combinatorial species}. Following \cite{AM2010}, we will focus in \emph{vector species}.

\subsubsection{Basic definitions} \label{subsubsec:basics}

Let $\mathbb{K}$ be a field of arbitrary characteristic. Let $\Fset$ be the category of finite sets and bijections between finite sets, and $\Vect_{\mathbb{K}}$ be the category of $\mathbb{K}$-vector spaces and linear maps between vector spaces. A \Def{vector species}, or simply a \Def{species}, is a functor $\tp: \Fset \to \Vect_{\mathbb{K}}$.

\medbreak

More precisely, a species is formed by a family of vector spaces $\tp[I]$, one for each finite set $I$, together with linear maps $\tp[\sigma]: \tp[I] \to \tp[J]$, one for each bijective map $\sigma: I \to J$. By the very definition of functor, such maps must satisfy
\[\tp[\tau \circ \sigma] = \tp[\tau]\circ \tp[\sigma] \qquad \text{ and } \qquad \tp[\text{id}_I]=\text{id}_{\tp[I]},\]
whenever $I \xrightarrow{\sigma} J \xrightarrow{\tau}K$ are composable bijections. In particular, every map $\tp[\sigma]: \tp[I] \to \tp[J]$ is an isomorphism, with inverse $\tp[\sigma^{-1}]$.

\medbreak

Two species $\tp$ and $\tq$ are said to be \Def{combinatorially equivalent} if there is a natural isomorphism $\alpha: \tp \Rightarrow \tq$ between functors. Explicitly, $\alpha$ consists of a family of isomorphisms $\alpha_I: \tp[I]\to \tq[I]$, one for each finite set $I$, such that $\tq[\sigma]\circ\alpha_I = \alpha_J \circ \tp[\sigma]$, for every bijection $\sigma: I \to J$. In this case, we writes $\tp = \tq$.

\medbreak

A species $\tp$ is said to be \Def{finite dimensional} if $\dim{\tp[I]}<\infty$ for each finite set $I$, \Def{connected} if $\tp[\emptyset] \cong \mathbb{K}$, and \Def{positive} if $\tp[\emptyset]=(0)$ where $(0)$ is the trivial vector space containing only the neutral element. The \Def{positive part} of a species $\tp$ is the positive species $\tp_+$ given by
\begin{equation} \label{eq:pos_part}
\tp_+[I]=\begin{cases}
\tp[I],& \text{ if } I\neq \emptyset\\
(0), & \text{ otherwise}
\end{cases}.
\end{equation}

Given a vector space $V$, let ${\bf 1}_V$ be the vector species defined by
\[\mathbf{1}_V[I]=\begin{cases}
V,& \text{ if } I= \emptyset\\
(0), & \text{ otherwise}
\end{cases}.\]

\

A \Def{morphism} between species $\tp$ and $\tq$ is a natural transformation between the functors. That is, a morphism $f: \tp \to 
\tq$ of species is a collection of linear maps $f=(f_I)_I$, one for every finite set $I$, with $f_I: \tp[I]\to \tq[I]$ satisfying $\tq[\sigma]\circ f_I = f_J \circ \tp[\sigma]$, for all bijection between finite sets $\sigma: I \to J$. This defines a category on species (over the field $\mathbb{K}$), denoted by $\Ss$. Unless specified, we do not makes distinction between the category of species and the full subcategory of finite-dimensional species.

There are several possible monoidal structures on this category. 
We will consider two of them: the \Def{Cauchy} and \Def{substitution} products $\cdot$ and $\circ$, respectively: for any finite set $I$,
\begin{equation}\label{cauchy}
(\tp \cdot \tq)[I]:=\bigoplus_{I = S \sqcup T}\tp[S]\otimes \tq[T];
\end{equation}
\begin{equation}\label{comp}
(\tp \circ \tq)[I]:=\bigoplus_{X \vdash I}\tp[X]\otimes \left(\bigotimes_{B \in X} \tq[S] \right).
\end{equation}

We denote by $(\Ss, \cdot)$ and $(\Ss, \circ)$ the resulting monoidal categories obtained from the Cauchy and substitution operations, respectively. The monoidal category $(\Ss, \cdot)$ is braided and symmetric, while $(\Ss, \circ)$ is not braided. 

\subsubsection{Algebraic structures on vector species} \label{subsubsec:alg}
Classical algebraic structures on vector spaces, such as associative algebras, coassociative coalgebras, bialgebras, Lie algebras and others, have their counterpart in the symmetric monoidal category of vector species, the tensor product of vector spaces being replaced by the Cauchy product of species.

\medskip
\paragraph{Monoids}
A \Def{monoid} in $(\Ss, \cdot)$  consist of a species $\tp$ equipped with morphisms of species
\begin{equation*}
    \mu: \tp \cdot \tp \to \tp \qquad \text{ and } \qquad \iota: \mathbf{1}_{\mathbb{K}} \to \tp
\end{equation*}
making both diagrams below commute
\begin{equation}\label{diag-monoid}    \xymatrix{\tp\cdot\tp\cdot\tp\ar[rr]^{\mu\cdot\text{id}_\tp}\ar[d]_{\text{id}_\tp\cdot\mu}
    &&\tp\cdot\tp\ar[d]^\mu\\
    \tp\cdot\tp\ar[rr]_\mu &&\tp
    }
    \hskip 15mm
    \xymatrix{\mathbf{1}_{\mathbb{K}}\cdot\tp\ar[r]^{\iota\cdot\text{id}_\tp}\ar[dr]_\simeq
    & \tp\cdot\tp\ar[d]_\mu&\tp\cdot\mathbf{1}_{\mathbb{K}}\ar[l]_{\text{id}_\tp\cdot\iota}\ar[dl]^\simeq\\
    &\tp&
    }
\end{equation}

The morphism $\mu=(\mu_I)_I$ consists on a collection of linear maps $\mu_I: (\tp \cdot \tp)[I]\to \tp[I]$, one for each finite set $I$, which in turns are defined as the sum of all possible maps
\[\mu_{S,T}: \tp[S] \otimes \tp[T]\to \tp[I],\]
for every $I=S\sqcup T$.  If $x \in \tp[S]$, $y \in \tp[T]$,  we writes $x \cdot y \in \tp[I]$ to denotes the image of $x\otimes y$ under $\mu_{S,T}$. For the unit map $\iota=(\iota_I)_I$, all components are defined as the zero map, except the $\emptyset$-component $\iota_\emptyset: \mathbb{K} \to \tp[\emptyset]$. The monoid $(\tp, \mu, \iota)$ is said to be \Def{commutative} if $x\cdot y = y \cdot x$, for all $I=S \sqcup T$, $x \in \tp[S]$ and $y \in \tp[T]$. Monoids in species are also called \emph{twisted algebras}.

\medbreak

Given two monoids $(\tp, \mu, \iota)$ and $(\tq, \mu', \iota')$, a \Def{morphism of monoids} $f: (\tp, \mu, \iota) \to (\tq, \mu', \iota')$ consists of a family $f=(f_I)_I$ of linear maps $f_I: \tp[I]\to \tq[I]$, one for every finite set $I$, such that the following diagrams commute:
\begin{equation}\label{monoidmorphism}    
\xymatrix{
\tp[S]\otimes \tp[T]\ar[rr]^{f_S \otimes f_T}\ar[d]_{\mu_{S,T}}
    &&\tq[S]\otimes \tq[T]\ar[d]^{\mu'_{S,T}}\\
    \tp[I]\ar[rr]_{f_I} &&\tq[I]
    }
    \hskip 15mm
    \xymatrix{   \tp[\emptyset]\ar[rr]^{f_\emptyset}&&\tq[\emptyset]\\   &\mathbb{K}\ar[lu]^{\iota_\emptyset}\ar[ru]_{\iota'_\emptyset}&
    }
\end{equation}

\medskip

\paragraph{Comonoids}
A \Def{comonoid} in $(\Ss, \cdot)$ corresponds to the dual notion of monoid in $(\Ss, \cdot)$. Formally, a comonoid consists of a species $\tp$ equipped with morphisms of species
\begin{equation*}
    \Delta: \tp \to \tp \cdot \tp \qquad \text{ and } \qquad \varepsilon: \tp \to \mathbf{1} .
\end{equation*}
such that the duals of both diagrams \eqref{diag-monoid} commute. That is, for each finite set $I$ and for each decomposition $I=S \sqcup T$, we have a linear map 
\begin{equation*}
    \Delta_{S,T}: \tp[I] \to \tp[S] \otimes \tp[T] \quad \text{ and } \quad\varepsilon_\emptyset: \tp[\emptyset] \to \mathbb{K}.
\end{equation*}

\medbreak

A comonoid $(\tp, \mu, \iota)$ in $(\Ssk, \cdot)$ is \Def{cocommutative} if for any finite disjoint sets $S, T$, and element $x\in \tp[S\sqcup T]$, we have $\Delta_{S,T}(x)=\Delta_{T,S}(x)$.

\medbreak

Given two comonoids $(\tp, \Delta, \varepsilon)$ and $(\tq, \Delta', \varepsilon')$, a \Def{morphism of comonoids} $f: (\tp, \Delta, \varepsilon) \to (\tq, \Delta', \varepsilon')$ is a family of linear maps $f=(f_I)_I$ satisfying the dual notions of \eqref{monoidmorphism}.

\medskip

\paragraph{Bimonoids  and Hopf monoids}
A \Def{bimonoid} $(\thh, \mu, \Delta, \iota, \varepsilon)$ in $(\Ss, \cdot)$ is a monoid and comonoid such that the diagram
\[\xymatrix{
\thh[A] \otimes \thh[B] \otimes \thh[C] \otimes \thh[D] \ar[rr]^-{\cong} && \thh[A] \otimes \thh[C] \otimes \thh[B] \otimes \thh[D] \ar[dd]^-{\mu_{A,C}\otimes\mu_{B,D}}\\
&&\\
\thh[S_1]\otimes \thh[S_2]\ar[uu]^-{\Delta_{A,B}\otimes \Delta_{C,D}} \ar[r]_-{\mu_{S_1, S_2}} & \thh[I] \ar[r]_-{\Delta_{T_1, T_2}} & \thh[T_1]\otimes \thh[T_2]
}\]

commutes, for every two decompositions $I=S_1\sqcup S_2=T_1 \sqcup T_2$ of a finite set $I$ with the following resulting pairwise intersections:
\[A:=S_1\cap T_1 \quad , \quad B:=S_1 \cap T_2 \quad , \quad C:=S_2 \cap T_1 \quad , \quad D:=S_2 \cap T_2.\]

\medbreak

We can define the convolution algebra $\text{End}_{\Ss}(\thh)$ as the set of natural transformations $l: \thh \to \thh$ with the convolution product $\star$ defined by
\begin{equation}
    (l\star m)_I=\sum_{S\sqcup T=I}\mu_{S,T}(l_S\otimes m_T)\Delta_{S,T}
\end{equation}
for any finite set $I$.

\medbreak

A morphism of species $s: \thh \to \thh$, is called and \Def{antipode} of $\thh$, if $\thh[\emptyset]$ is a Hopf algebra with antipode $s_\emptyset: \thh[\emptyset] \to \thh[\emptyset]$, and for each nonempty set $I$, we have
\begin{equation}
    \sum_{S \sqcup T = I} \mu_{S,T}(\text{id}_S \otimes s_T)\Delta_{S,T} = 0 = \sum_{S \sqcup T = I} \mu_{S,T}(s_S \otimes \text{id}_T)\Delta_{S,T}.
\end{equation}

A \Def{Hopf monoid} in $(\Ssk, \cdot)$ is a bimonoid along with an antipode $s:\thh \to \thh $.
Equivalently, the identity map $\text{id}_{\thh}$ is invertible in $\text{End}_{\Ssk}(\thh)$ (c.f. \cite{AM2010}, section 2.7).
Recall that this is also the case in the classical Hopf algebras.

\

\subsubsection{Examples of species} 

We give some examples of species which will play important roles in this paper.

\paragraph{Species of singletons}

Let $\mathbf{X}$ be the \Def{species of singletons} defined on any finite set $I$ by
\begin{equation} \label{eq:species_sing}
 \mathbf{X}[I] = \begin{cases} 
 & \mathbb{K}\quad\text{if }|I|=1 \\
 & (0)\quad\text{otherwise}.
 \end{cases}
\end{equation}
The action of $\mathbf{X}$ on bijections is trivial, and $\mathbf{X}$ is a species.

\paragraph{Species of decorated exponential}

Let $V$ be a vector space, referred as the \emph{space of decorations}. Given a finite set $I$, the set of all functions $f:I \to V$, denoted by $V^{\otimes I}$, is a vector space. Moreover, every bijection $\sigma: I \to J$ between finite sets induces an isomorphism $V^I \cong V^J$. This defines a species $\wE_V$, called the \Def{decorated exponential species}. When $V \cong \mathbb{K}$, we obtains the \Def{exponential species} $\wE_V = \wE$.

\medskip

\paragraph{Species of lists}

Let $I$ be a finite set and let $n:=|I|$. A \Def{linear order} on $I$ is a bijection $\ell: [n] \to I$. When $n=0$, there is a unique linear order $\varnothing: \emptyset \to \emptyset$.  If $\sigma: I \to J$ is a bijection between finite sets, the map $\ell \mapsto \sigma \circ \ell$ sends a linear order $\ell$ of $I$ to a new linear order on $J$. This defines the connected species $\wL$ of linear orders. 

\medbreak

If $\ell$ is a linear order on $I$, we write $i_1 \leq_\ell i_2$ if $\ell^{-1}(i_1)\leq \ell^{-1}(i_2)$, for all $i,j \in I$. This endows the set $I$ with a total order, for each linear order $\ell$. In general, we write $\ell = i_1|i_2| \cdots |i_n$ to denote the linear order $\ell$ of $I$, for which $I=\{i_1, i_2, \hdots, i_n\}$ and $i_1 <_\ell i_2 <_\ell \cdots <_\ell i_n$.

\medbreak

If $S \subseteq I$ and $\ell$ is a linear order on $I$, there exists a unique linear order on $S$, denoted by $\ell \cap S$, such that $s_1 \leq_{\ell \cap S} s_2$ if and only if $s_1 \leq_\ell s_2$. Writing $\ell$ as a word, $\ell \cap S$ is the new word obtained by removing elements in $I \setminus S$ and keeping the elements from $S$  in the same relative order.

\medskip

\paragraph{Species of rooted trees and forests} \label{subsubsec:trees}

A \Def{rooted tree} is a connected graph which is acyclic and that has a distinguished vertex, called the \emph{root}. A set of rooted trees is called a \Def{rooted forest}. 

\medbreak

Let $\taa$ be the species of rooted trees. A linear basis of $\taa[I]$ is given by the set of rooted trees on the vertex-set $I$. We put $\taa[\emptyset]= (0)$, so $\taa$ is positive. This allows to define the species $\mathbf{F}:=\wE\circ \taa$ of \Def{rooted forests}. By definition, $\mathbf{F}$ is connected.

\

\subsection{Universal constructions on species}

\subsubsection{Species applied to compositions and partitions} \label{subsubsec:comp_part}

Let $\tq$ be a species. Given a finite set $I$ and set composition $F=B_1|B_2|\cdots | B_k$ of $I$, define
\begin{equation}\label{qf}
    \tq(F):=\tq[B_1]\otimes \tq[B_2] \otimes \cdots \otimes \tq[B_k].
\end{equation}

Let $J$ a finite set, and let $G=C_1|C_2|\cdots | C_\ell$ be a set composition of $J$. There is a canonical identification
\begin{equation}\label{identificationFG}
    \tq(F)\otimes \tq(G) \cong \tq(F \cdot G).
\end{equation}
If $x_{B_1} \in \tq[B_1], \hdots, x_{B_k} \in \tq[B_k]$, the element $x=x_{B_1}\otimes \cdots \otimes x_{B_k} \in \tq(F)$ can be seeing as a \emph{non-commutative monomial} on the elements $x_{B_1}, \hdots, x_{B_k}$. We sometimes writes 
\begin{equation}\label{nc-monomial}
x= \bigotimes_{B \in F}x_B  
\end{equation}
to represent the non-commutative monomial $x$.

\medbreak

We now focus on a similar construction for set partitions. The category of vector spaces over $\mathbb{K}$ is endowed with a symmetric monoidal structure $(\mathsf{Vec}, \otimes, \beta)$. In there, we can consider unbracketed tensor products due to the symmetry $\beta$. Also, for every bijection $\sigma$ on $k$ elements $i_1, \hdots, i_k$, there is a well-defined isomorphism
\[\theta_\sigma: V_{i_1}\otimes \cdots \otimes V_{i_k} \to V_{\sigma(i_1)}\otimes \cdots \otimes V_{\sigma(i_k)}\]
defined via the symmetry $\beta$. If $F,G \vDash I$ are such that $\text{supp}_I(F) = \text{supp}_I(G)$, let $\theta_{F,G}: \tq(F)\to \tq(G)$ be the induced isomorphism from the (unique) permutation $\sigma\in \mathfrak{S}_{\ell(F)}$
sending $F$ to $G$. 

\medbreak

Given a partition $\pi=\{B_1, \hdots, B_k\} \vdash I$, let $\tq(\pi)$ be the vector space such that, for every set composition $F=B_1|\cdots|B_k \vDash I$, there is an isomorphism
\[\Theta_F: \tq(F)=\tq[B_1]\otimes \cdots \otimes \tq[B_k] \to \tq(\pi),\]
such that $\Theta_G \circ \theta_{F,G}=\Theta_F$, for any $G \vDash I$ such that $\text{supp}_I(G)=\text{supp}_I(F)$. If $x_B \in \tq[B]$ for every $B \in \pi$, we denotes by 
\begin{equation}\label{c-monomial}
x= \bigotimes_{B \in \pi}x_B
\end{equation}
an element in $\tq(\pi)$, which can be seeing as a \emph{commutative monomial} on the elements $\{x_B: B \in \pi\}$ (compare this with its non-commutative analogue \eqref{nc-monomial}). The space $\tq(\pi)$ can be seen as a quotient of the direct sum of the spaces $\tq(F)$ with $F$ having a fixed support. With the approach presented above, we are able to treat together species applied to set compositions and set partitions.  For an explicit description of the unordered tensor product on $(\mathsf{Vec}, \otimes, \beta)$, see \cite{AM2010}, Example 1.30.


\subsubsection{Free monoid} \label{FreeMonoid}
Given a positive species $\tq \in \Ss_+$, an $(\wL\circ \tq)$-structure on $I$ is formed by a partition $X \vdash I$, with a linear order on the set blocks of $X$, and finally a $\tq$-structure on each block of $X$. This is equivalent to a set composition such that each block of $F$ is endowed with a $\tq$-structure. Therefore, we have
\begin{equation}
    (\wL \circ \tq)[I]=\bigoplus_{F \vDash I}\tq(F).
\end{equation}
We consider now a product in $\wL\circ\tq$ as follows. Let $I$ be a finite set and let $I=S \sqcup T$ be a decomposition of $I$. If $F \vDash S$ and $G \vDash T$, recall from \eqref{identificationFG} that we have a canonical map $\tq(F)\otimes \tq(G) \to \tq(F \cdot G)$. Let $\mu_{S,T}$ be the direct sum of all the maps associated to every pair of set compositions $F$ and $G$ of the set $I$ and $J$, respectively. Hence, the map $\mu_{S,T}$ is of the form
\[\mu_{S,T}: (\wL\circ \tq)[S]\otimes (\wL\circ \tq)[T] \to (\wL\circ \tq)[I].\]

Taking $\mu$ as the collection of the sums of all possible maps $\mu_{S,T}$, for every finite set $I$ and every decomposition $I =S\sqcup T$, gives a morphism of species
\[\mu: (\wL\circ \tq)\cdot (\wL\circ \tq) \to \wL\circ \tq. \]
    
The product $\mu$ is called \Def{concatenation}. Hence, $\mathcal{T}(\tq):=(\wL\circ \tq, \mu, \iota)$ is a monoid, where $\iota$ is the canonical embedding of the species $\mathbf 1_{\mathbb K}$ into $\wL \circ \tq$, which identifies $\mathbb K=\mathbf 1_{\mathbb K}[\emptyset]$ with $(\wL \circ \tq)[\emptyset]$.

\medskip

\medbreak

\subsubsection{Free commutative monoid}\label{FreeCMonoid} 
Given a positive species $\tq \in \Ss_+$, an $(\wE\circ \tq)$-structure on $I$ is formed by a partition $\pi \vdash I$, with a $\wE$-structure on the set of blocks of $\pi$, and finally a $\tq$-structure on each block of $X$. Therefore, we have
\begin{equation}
    (\wE \circ \tq)[I]=\bigoplus_{\pi \vdash I}\tq(\pi).
\end{equation}

The \Def{commutative concatenation} is the map
\[\mu: (\wE \circ \tq)\cdot (\wE \circ \tq) \to \wE \circ \tq\]
defined as follows: each component $\mu_I$ is the direct sum of maps
\[\mu_{S,T}: (\wE \circ \tq)[S]\otimes (\wE \circ \tq)[T]\to (\wE \circ \tq)[I], \]
one for each decomposition $I = S\sqcup T$, such that $\mu_{S,T}$ is the direct product of the maps
\[\tq(\pi_1)\otimes \tq(\pi_2) \xrightarrow{\cong} \tq(\pi_1 \sqcup \pi_2),\]
where $\pi_1 \vdash S$ and $\pi_2 \vdash T$.

\medskip

As $\wL\circ\tq$, the species $\wE\circ\tq$ can be endowed with a comonoidal structure. This will be used explicitly to construct up coderivation from branching maps in section~\ref{subsec:up_down_grafting}. We describe it explicitly in the following lines.

\medbreak

Given $\pi \in \Pi[I]$ and a decomposition $I = S \sqcup T$, recall from the Preliminaries Section \eqref{sec:Preliminaries} that $S$ is $\pi$-admissible if $\pi=\pi|_S \cup \pi|_T$. The \Def{cocommutative deconcatenation} is the map
\[\Delta: \wE \circ \tq \to (\wE \circ \tq)\cdot (\wE \circ \tq),\]
where each $\Delta_I$ is the direct sum of the morphisms
\[\Delta_{S,T}: (\wE \circ \tq)[I] \to (\wE \circ \tq)[S]\otimes (\wE \circ \tq)[T].\]

In turn, each $\Delta_{S,T}$ is completely determined by the following rule:
\begin{align*}
\tq(\pi)  & \quad \to \quad \tq(\pi|_S) \otimes \tq(\pi|_T)\\
\bigotimes_{B \in \pi}x_B  & \quad \mapsto \quad 
\begin{cases}
\displaystyle \bigotimes_{B \in \pi|_S}x_B \, \otimes  \, \bigotimes_{B \in \pi|_T}x_B, & \text{ if } S \text{ is } \pi\text{-admissible;}\\
0, & \text{ otherwise}.
\end{cases}
\end{align*}

For every positive species $\tq$, there is a coproduct compatible with monoid $\mathcal{T}(\tq)$, called the \emph{deshuffle} (see \cite{AM2010}, section 11.4.2). The dual structure consists of the \emph{shuffle} monoid on $\wL \circ \tq$, together with the \emph{deconcatenation}. It is well-known that both bimonoids are not isomorphic. For the case of $\wE \circ \tq$, the situation is different: the dual map of the cocommutative deconcatenation is precisely the concatenation product for $\wE \circ \tq$. For a more detailed discussion, see \cite{AM2010}, section 11.6.4.





\subsection{Further operations on species}

Given a finite set $I$, we denotes by $I^+:= I \cup \{*_I\}$ the new set formed by appending to $I$ a new element\footnote{If we want to be rigorous, we can make $*_I:=I$.} $*_I$, not included on $I$. 

\subsubsection{The derivative operation} \label{subsubsec:deriv}
Let $\tp$ be a species. The \Def{derivative} $\tp'$ of $\tp$ is a new species defined as $\tp'[I]:=\tp[I^+]$. If $\sigma: I \to J$ is a bijection between finite sets, the action of $\tp'$ on $\sigma$ is given by $\tp'[\sigma]:=\tp[\sigma^+]$, where $\sigma^+: I^+ \to J^+$ is given by $\sigma^+(i):=\sigma(i)$ for all $i \in I$, and $\sigma^+(*_I):=*_J$.

\medbreak

Any decomposition of $I^+$ into two disjoint subsets is equivalent to a decomposition of $I$ of the same kind. Indeed, given $I=S \sqcup T$, then any decomposition of $I^+$ is of the form $I^+= S^+ \sqcup T$ or $I^+ = S \sqcup T^+$. This imply that derivatives of species behaves well with respect to the Cauchy product:
\begin{equation}\label{DerivativeRule}
(\tp \cdot \tq)' = \tp' \cdot \tq + \tp \cdot \tq'.  
\end{equation}

\medbreak

The derivative operation $\tp \mapsto \tp'$ defines an endofunctor on the category of species. Given a morphism of species $f: \tp \to \tq$, its derivative $f': \tp' \to \tq'$ is given by the following composition of maps
\begin{equation}\label{derivativef}
f'_I: \tp'[I] = \tp[I^+] \xrightarrow{\, f_{I^+}\,}\tq[I^+] = \tq'[I].
\end{equation}

Similarly, we have the following combinatorial chain rule for species: given $\tp, \tq$ species with $\tq$ positive, then
\begin{equation}
(\tp \circ \tq)'=(\tp'\circ \tq)\cdot \tq'.
\end{equation}

\medbreak

We specify now the derivative maps obtained from a product and coproduct maps from a monoid and comonoid, respectively. First, if $(\tp, \mu)$ is a monoid, then the derivative of the map $\mu: \tp \cdot \tp \to \tp$ is of the form 
\begin{equation}
\mu': \tp' \cdot \tp + \tp \cdot \tp' \to \tp'.
\end{equation}

\medbreak

Every decomposition of the set $I^+ = I \sqcup \{*_I\}$ is equivalent to a decomposition $(S,T)$ of $I$, where one of the blocks is selected (thus indicating where is located the element $*_I$). Hence, $\mu'_{I}$ is the sum of all maps
\begin{equation}\label{DerivativeMonoid}
\mu'_{S,T} := \mu_{S^+,T} + \mu_{S,T^+},
\end{equation}
one for each decomposition $I=S\sqcup T$. 

\medbreak

Similarly, if $(\tp , \Delta)$ is a comonoid, the derivative of the map $\Delta: \tp \to \tp\cdot \tp$ is the map
\begin{equation}
\Delta': \tp' \to \tp'\cdot \tp + \tp \cdot \tp',
\end{equation}

where for every finite set $I$ the linear map $\Delta'_I$ is the sum of the components
\begin{equation}\label{DerivativeComonoid}
\Delta'_{S,T} := \Delta_{S^+,T} + \Delta_{S,T^+},
\end{equation}
one for each decomposition $I = S\sqcup T$.

\

\subsubsection{Hadamard product} \label{subsubsec:Hadamard} 

Given two species $\tp$ and $\tq$, the \Def{Hadamard product} $\tp \times \tq$ is the new species defined as
\begin{equation}\label{Hadamard}
    (\tp \times \tq)[I]:=\tp[I]\otimes \tq[I],
\end{equation}
for every finite set $I$, and similarly on bijections.


\subsubsection{Dual species} \label{subsubsec:dual} 

For any species $\tq$ we denote $\tq^*$ the associated \Def{dual species} defined by $\tq^*[I]:=\tq[I]^*$ and for any bijection $\sigma:I\to J$ we set $\tq^*[\sigma]:=\tq[\sigma^{-1}]^*:\tq^*[I]\to\tq^*[J]$. Then if $\tp$ and $\tq$ are finite dimensional we have $(\tp\cdot\tq)^*=\,\tp^*\cdot\tq^*$ and $(\tp\circ\tq)^*=\,\tp^*\circ\tq^*$.

This construction lifts to a duality functor $(-)^*:\Ss^{\text{op}}\to\Ss$ between the opposite category of $\Ss$ and $\Ss$. When restricted to finite dimensional species, it is a monoidal functor for the Cauchy, the composition, and the Hadamard products.

\subsection{Species with up and down operators}

The notions of up and down operators on species are developed in Aguiar and Mahajan \cite[Section 8.12]{AM2010}. These operators are closely related to the theory of creation and annihilation operators on Fock spaces. After recalling the main definitions, we introduce special up and down operators that are compatible with monoids and comonoids in species.

\begin{definition} \label{def:up_down}
A \Def{species with up operator} (resp. a \Def{species with down operator}) is a species $\tp$ with a morphism of species $u: \tp \to \tp'$ (resp. $d: \tp' \to \tp$). A \Def{morphism of species with up operators} (resp. a \Def{morphism of species with down operators}) $f:(\tp, u) \to (\tq, v)$ (resp. $f:(\tp, d) \to (\tq, e)$) is a morphism of species $f:\tp \to \tq$ that intertwines $u$ with $v$ (resp. $d$ with $e$), i.e. such that the following diagrams commute
\begin{equation*}
\SelectTips{eu}{12}
\xymatrix@R+1pc@C+2pc{
\tp \ar[r]^{f} \ar[d]_{u}& \tq \ar[d]^{v}\\
\tp' \ar[r]_{f'}& \tq'
} \qquad  \qquad
\SelectTips{eu}{12}
\xymatrix@R+1pc@C+2pc{
\tp' \ar[r]^{f'} \ar[d]_{d}& \tq' \ar[d]^{e}\\
\tp \ar[r]_{f}& \tq
}.
\end{equation*}

We denotes $\Ss^u, \Ss_d$ and $\Ss^u_d$ the categories of species with up operators, down operators, and up and down operators, respectively.
\end{definition}

\medbreak

Whenever $\tp$ and $\tq$ are (co)monoids, we may consider up and down operators that respect the (co)monoidal structures between species.
\begin{definition} \label{def:up_down_deriv} 
Given a monoid $(\tp, \mu)$ in $(\Ss, \cdot)$, a \Def{down derivation} of $\tp$ is a down operator $d: \tp' \to \tp$ such that the following diagram commute
\begin{equation}\label{diag:down_derivation}
\SelectTips{eu}{12}\xymatrix@R+1pc@C+2pc{
(\tp \cdot \tp)' \ar[r]^{d\cdot \normalfont\text{id} \,+\, \normalfont\text{id}\cdot d} \ar[d]_{\mu'}& \tp \cdot \tp \ar[d]^{\mu}\\
\tp' \ar[r]_{d}& \tp
}
\end{equation}

Dually, let $(\tp, \Delta)$ be a comonoid in $(\Ssk, \cdot)$. An \Def{up coderivation} of $\tp$ is an up operator $u: \tp \to \tp'$ such that the following diagram commutes
\begin{equation} \label{diag:up_coderivation}
\SelectTips{eu}{12}\xymatrix@R+1pc@C+2pc{
\tp \ar[r]^u \ar[d]_{\Delta} & \tp' \ar[d]^{\Delta'}\\
\tp \cdot \tp \ar[r]^{u\cdot \normalfont\text{id} \,+\, \normalfont\text{id}\cdot u}& (\tp \cdot \tp)' 
}.
\end{equation}

\end{definition}

\medbreak

Notice that the maps $d\cdot \text{id} + \text{id}\cdot d$ and $u\cdot \text{id} + \text{id}\cdot u$ appear from the calculus rule \eqref{DerivativeRule}.

\medbreak

Let $I$ be a finite set and let $I=S \sqcup T$ be a decomposition. If $(\tp, \mu, d)$ is a monoid with down derivation, we have
\begin{equation}
d_I\mu'_{S,T}(x \otimes y)= \mu_{S,T}(d_S(x)\otimes y + x \otimes d_T(y)), \quad \forall \, x \in \tp[S], \forall \,  y \in \tq[T],
\end{equation}
where $\mu'_{S,T}$ is defined as in \eqref{DerivativeMonoid}. Similarly, whenever $(\tp, \Delta, u)$ is a comonoid with up-derivation, then
\begin{equation}
(\Delta')_{S,T}(u_I(x)) = (u_S \otimes \text{id} + \text{id}\otimes u_T)\Delta_{S,T}(x), \quad \forall \, x \in \tp[I],
\end{equation}
where $(\Delta')_{S,T}$ is given as in \eqref{DerivativeComonoid}.

More precisely, given a monoid $(\tp, \mu)$, then $d: \tp' \to \tp$ is an down derivation if, and only if
\begin{equation*}
    d_I((\mu')_{S,T}(x\otimes y)) = \mu_{S,T}\left(d_S(x)\otimes y+x\otimes d_T(y)\right)
\end{equation*}
for any finite set $I$, any decomposition $S\sqcup T=I\cup\{*\}$ and any $x\otimes y\in\tp[S]\otimes\tp[T]$.

Similarly, given an comonoid $(\tp,\Delta)$, then $u:\tp\to\tp'$ is an up coderivation if, and only if
\begin{equation*}
    (\Delta')_{S,T}(u_I(x))=(u_S\otimes \text{id}+\text{id}\otimes u_T)\Delta_{S,T}(x)
\end{equation*}
for any finite set $I$, any decomposition $S\sqcup T=I$ of $I$ and any $x\in\tp[I]$.


\begin{remark}
    One also has up derivations and down coderivations. These notions can easily be transferred from these presented above, however we will not need them in this article and we therefore omit these definitions.
\end{remark}

\section{Grafting and cografting} \label{sec:graft_cograft}

In this section we introduces an analogue of the classical grafting map of forests in the context of species.

\subsection{Species with grafting and cografting}

\subsubsection{Graftings and cograftings}

\begin{definition} \label{def:graf_cograft} 
Let $\tq$ be a positive species. A \Def{grafting map} $\Bb$ for $\tq$ is a species map $\Bb: \wE \circ \tq \to \tq'$. An \Def{ordered grafting map} $\sBb$ for $\tq$ is a map $\sBb: \wL\circ \tq \to \tq'$. Dually, a \Def{cografting map} $\bB$ for $\tq$ is a map $\bB: \tq' \to \wE\circ \tq$, and an \Def{ordered cografting map} $\bBs$ is a map $\bBs: \tq' \to \wL \circ \tq$.

\medbreak

A \Def{morphism of species with grafting maps} $f:(\wE\circ\tp, \Bb) \to (\wE\circ\tq, \Cc)$ is a morphism of species $f:\wE\circ \tp \to \wE\circ\tq$ such that the following diagram commutes:
\begin{equation} \label{eq:def_morph_graft}
\SelectTips{eu}{12}\xymatrix{
\wE \circ \tp \ar[r]^{f} \ar[d]_{\Bb} & \wE \circ \tq \ar[d]^{\Cc}\\
\tp' \ar[r]_{f'} & \tq'
}
\end{equation}
where in the bottom line, $\tp'$ and $\tq'$ are seen as included in $(\wE\circ\tp)'$ and $(\wE\circ\tq)'$ respectively. To be more precise, for each finite set $I$, the map $f_I:(\wE\circ\tp)[I]\longrightarrow(\wE\circ\tq)[I]$ decomposes as
\begin{equation*}
    f_I=\sum_{\pi\vdash I}f_\pi,\qquad f_\pi:\tp(\pi)\longrightarrow\tq(\pi).
\end{equation*}
Then the bottom line of equation~\eqref{eq:def_morph_graft} should be read as the map 
\[f_{I^+}:\tp[I^+]=\tp'[I]\longrightarrow\tq'[I]=\tq[I^+].\]

Similarly, a \Def{morphism of species with cografting maps} $f:(\wE\circ\tp, \bB) \to (\wE\circ\tq, \cC)$ is a morphism of species $f: \tp \to \tq$ such that the following diagram commutes:
\begin{equation}
\SelectTips{eu}{12}\xymatrix{
\tp' \ar[r]^{f'} \ar[d]_{\bB} & \tq' \ar[d]^{\cC}\\
\wE \circ \tp \ar[r]^{f}  & \wE \circ \tq 
}
\end{equation}

    The same definitions define morphisms of species with ordered grafting and cografting maps. Thus we have the species $\Spgr$, $\Spcogr$, $\Spogr$, and $\Spocogr$ of species with grafting, species with cograftings, species with ordered graftings, and species with ordered cografting respectively. As before, we use the same notations for their finite dimensional counterparts.
\end{definition}
\begin{remark} \label{rk:dec}
    In practice, we often build a morphism between two species with grafting $f:(\wE\circ\tp,\Bb)\to(\wE\circ\tq,\Cc)$ from a morphism $g:\tp\to\tq$. Indeed the species map $f:=\normalfont \text{id}_{\mathbf{E}} \circ g$ is a morphism of species with grafting if the following diagram commutes.
    \[\xymatrix{
    \wE \circ \tp \ar[r]^{\normalfont\text{id}_{\mathbf{E}} \circ g} \ar[d]_{\Bb} & \wE \circ \tq \ar[d]^{\Cc}\\
    \tp' \ar[r]_{g'} & \tq'
    }\]
    This construction amounts to defining $f_\pi=\bigotimes_{B\in\pi}g_B$ for any $\pi\vdash I$. The same construction works for cograftings and for the ordered versions of these objects. We call \Def{decomposable} such a morphism of species with graftings, and set $\Spgr_{\rm dec}$ (resp. $\Spogr_{\rm dec}$) the category of species with graftings (resp. ordered graftings) with these morphisms.
\end{remark}
We finish this subsection with one last definition, of grafting maps that are morphisms of monoids, where the monoid structures of $\wL\circ\tp$ and $\wE\circ\tp$ we introduced in \S \ref{FreeMonoid} and \ref{FreeCMonoid} respectively. 
\begin{definition} \label{def:graft_mon}
    If $\tq'$ is endowed with a monoidal structure, $\Bb$ and $\sBb$ (resp. $\bB$ and $\bBs$) are called \Def{monoidal grafting maps} (resp. \Def{monoidal cografting maps}) if they are morphism for the corresponding monoidal structures.
\end{definition}

\subsubsection{The ordering and unordering functors}

The categories $\Spgr$, $\Spcogr$, $\Spogr$, and $\Spocogr$ are related, as we now explain.

\

Let $\tp, \tp_1, \tp_2$ be finite-dimensional species. The following combinatorial equivalences hold
\begin{equation}
(\tp^*)^*=\tp \quad , \quad (\tp_1 \cdot \tp_2)^* = \tp_1^* \cdot \tp_2^* \quad , \quad (\tp_1 \circ \tp_2)^* = \tp_1^* \circ \tp_2^*,
\end{equation}
where $\tp_2[\emptyset]=(0)$ in the latter case. 

\medbreak

As $\wE$ and $\wL$ are finite dimensional we have $\wE =\wE^*$ and $\wL =\wL^*$. Then, given a grafting $\Bb: \wE \circ \tq \to \tq'$ (resp. cografting $\bB:\tq'\to\wE\circ\tq$) on a finite dimensional species $\tq$, we obtains a cografting (resp. grafting) map on $\tq^*$, given by the composition of the following maps:
\begin{align*}
& \hspace{1.5cm}\xymatrix{
(\tq^*)' = (\tq')^* \ar[r]^-{\Bb^*} & (\wE \circ \tq)^* \ar[r]^{=} & \wE^* \circ \tq^* \ar[r]^{=} & \wE \circ \tq^*.
}\\
& \left(\text{resp. }  \xymatrix{\wE\circ\tq^* \ar[r]^{=} & \wE^*\circ\tq^* \ar[r]^{=} & (\wE\circ\tq)^* \ar[r]^{\bB^*} & (\tq')^*=(\tq^*)'} \right)
\end{align*}
Similarly, any ordered grafting (resp. ordered cografting) on $\tq$ induces an ordered cografting (resp. ordered grafting) map on $\tq^*$. Since the duality of species is functorial, so are these constructions.

We relate now grafting and ordered graftings. First, there is a commutative diagram of bimonoids:
\begin{equation}
\SelectTips{eu}{12}
\xymatrix{
\wL \ar[r] \ar[d]_{\pi} & \wL^*\\
\wE \ar[r] & \wE^* \ar[u]_{\pi^*}
}
\end{equation}
This is diagram (8.34) in \cite{AM2010}. Given a finite set $I$, the map $\pi: \wL \to \wE$ is given by $\pi_I(\ell):= *_I$, for every linear order $\ell \in \wL[I]$; the horizontal arrow $\wE \to \wE^*$ is the isomorphism given by $*_I \mapsto (*_I)^*$, and $\wL \to \wL^*$ sends any linear order $\ell$ of $I$ to the sum $\sum_{\ell'}(\ell')^*$, where the sum is taken over all linear orders $\ell'$ of $I$. 

From a species with grafting $(\tq, \Bb)$, the map 
\[
\SelectTips{eu}{12}
\xymatrix{
\wL\circ\tq  \ar[r]^-{\pi\circ \text{ id }} & \wE \circ \tq \ar[r]^-{\Bb} & \tq'
}\]

defines an ordered grafting map on $\tq$. If $f:(\wE\circ\tp,\Bb)\to(\wE\circ\tq,\Cc)$ is a morphism of species with grafting operators, then for any finite set $I$ and any $\pi\vdash I$ define
\begin{align*}
    \tilde f_\pi:\wL[\pi]\otimes\tp(\pi) & \longrightarrow \wL[\pi]\otimes\tq(\pi) \\
    \ell\otimes x_1\cdots x_k & \longrightarrow \ell\otimes f_\pi(x_1\cdots x_k).
\end{align*}
Then
\[\SelectTips{eu}{12}\xymatrix{
    \wL \circ \tp \ar[r]^{\tilde f} \ar[d]_{\pi\circ\text{id}} & \wL \circ \tq \ar[d]_{\pi\circ\text{id}} \\
    \wE\circ\tp \ar[r]^{ f} \ar[d]_{\Bb} & \wE \circ \tq \ar[d]^{\Cc}\\
    \tp' \ar[r]_{f'} & \tq'
    }\]
    commutes. The bottom square since $f$ is by hypothesis a morphism of species, and the top square by definition of $\tilde f$. Thus $\tilde f$ is a morphism of ordered grafting operators for the grafting operators built from $\Bb$ and $\Cc$ since in $f'$ and $\tilde f'$ coincide in the bottom line. Thus we have a functor $\fO:\Spgr\to\Spogr$.

\medskip   

    On the other hand, if $(\tq,\sBb)$ is a species with ordered grafting map then the composition
    \[
\SelectTips{eu}{12}
\xymatrix{
\wE\circ\tq \ar[r]^-{\simeq} & \wE^*\circ\tq \ar[r]^-{\pi^*\circ \text{ id }} & \wL^* \circ \tq \ar[r]^-{\simeq} & \wL \circ \tq 
\ar[r]^-{\sBb} & \tq'
}\]
defines a grafting map. We denote this grafting map by $\fU\sBb:\wE\circ\tq\longrightarrow\tq'$.

For the morphisms, we need to restrict the category we are working in. A morphism of species with ordered grafting $f:(\wL\circ\tp,\Bb)\longrightarrow(\wL\circ\tq,\Cc)$ is said to be \Def{order independent} if it exists a species map $\tilde f:\wE\circ\tp\longrightarrow\wE\circ\tq$ such that for any finite set $I$, any $\pi\vdash I$ and any $\ell\otimes x_1\cdots x_k\in\wL[\pi]\otimes\tp(\pi)$ we have 
    \begin{equation*}
        f_\pi(\ell\otimes x_1\cdots x_k)=\ell\otimes\tilde f_\pi(x_1\cdots x_k).
    \end{equation*}
Then for such an order independent map, let us set $\fU f:=\tilde f$. It is straightforward to check that $\fU f:(\wE\circ\tp,\fU\sBb)\longrightarrow(\wL\circ\tq,\fU\sCc)$ is a morphism of species with grafting. Thus we have a functor $\fU:\Spogr_{\rm pres}\longrightarrow\Spgr$, where $\Spogr_{\rm pres}$ is the category of species with ordered graftings as objects and order independent morphisms as morphisms.
\begin{remark}
    The above construction can be generalised to \Def{splitting} map $f:\wL\circ\tp\longrightarrow\wL\circ\tq$, that is, map that can be written $f=f_1\times f_2$, with $f_1:\wL\longrightarrow\wL$ and $f_2:\wE\circ\tp\longrightarrow\wE\circ\tq$ two morphisms of species by setting $\fU f:=f_2$. Another approach is to define for any morphism of species with ordered grafting $f:(\wL\circ\tp,\Bb)\longrightarrow(\wL\circ\tq,\Cc)$ a morphism $\fU f:\wE\circ\tp\longrightarrow\wE\circ\tq$ by setting, for any finite set $I$, any $\pi\vdash I$ and any $x_1\cdots x_k\in\tp(\pi)$
    \begin{equation*}
        \fU f(x_1\cdots x_k):=\sum\ell\otimes f(x_1\cdots x_k).
    \end{equation*}
    However, one can check than these constructions are not functorial. In other words, we only have a functor (which we denote as the previous one to not complicate notations) $\fU:{\rm Disc}(\Spogr)\longrightarrow{\rm Disc}(\Spgr)$, where ${\rm Disc}(\Spogr)$ (resp. ${\rm Disc}(\Spgr)$) is the discrete category underlying $\Spogr$ (resp. $\Spgr$).
\end{remark}

\medskip 

The same type of constructions can be made for cograftings and ordered cograftings. We omit these details since cograftings play a less important role in this paper.




\subsubsection{Properties of the functors}

We have the following:
\begin{prop} \label{prop:functors_grft_ogrft}
\begin{enumerate}[label=(\roman*)]
    \item The duality functor lifts to involutive functors 
    \begin{equation*}
        (-)^*:(\Spgr)^{\rm op}\to\Spcogr\qquad (-)^*:(\Spogr)^{\rm op}\to\Spocogr.
    \end{equation*}
    The same statement holds for ordered (co)grafting.
    \item The functor $\fU:\Spogr_{\rm pres}\longrightarrow\Spgr$ is faithful and dense.
    \item The functor $\fO:\Spgr\longrightarrow\Spogr_{\rm pres}$ is full and faithful.
%
\end{enumerate}
\end{prop}
\begin{proof}
Through this proof, we will use the notations $\ell\otimes x_1\cdots x_k$ and $x_1\cdots x_k$ indiscriminately for elements of $\wL[\pi]\otimes\tp(\pi)\subseteq \wL\circ\tp[I]$ and of $\tp(\pi)\subseteq\wE\circ\tp[I]$ respectively, for any species $\tp$, for a finite set $I$ and $\pi\vdash I$. Thus we will omit to mention these objects, and will also sometimes omit the subscripts $I$ or $\pi$ in morphisms of species. Both are done to increase the readability of the proof.

Now, let us prove each point individually. 
\begin{enumerate}[label=(\roman*)]
    \item This was proved in the discussion above.
    \item Let $(\wL\circ\tp,\sBb)$ and $(\wL\circ\tq,\sCc)$ be any two species with ordered grafting and 
    $$f:\fU(\wL\circ\tp,\sBb)=(\wE\circ\tp,\fU\sBb)\longrightarrow\fU(\wL\circ\tq,\sCc)=(\wE\circ\tq,\fU\sCc)$$
    be any morphism of species with grafting. Then, assume there are two order independent morphisms of species with ordered grafting $g_1,g_2:(\wL\circ\tp,\sBb)\longrightarrow(\wL\circ\tq,\sCc)$ such that $\fU g_1=\fU g_2=f$. Then since the $g_i$ are order independent we have (with some abuse of notations) $g_i(\ell\otimes x_1\cdots x_k)=\ell\otimes\fU g_i(x_1\cdots x_k)$ by definition of $\fU$. Then by assumptions (and with the same abuse of notations):
    \begin{equation*}
        g_1(\ell\otimes x_1\cdots x_k)=\ell\otimes f(x_1\cdots x_k) = g_2(\ell\otimes x_1\cdots x_k).
    \end{equation*}
    Thus $g_1=g_2$ and $\fU$ is faithful.

    Now, let $(\wE\circ\tp,\Bb)$ be any species with grafting. Then, let us define $\sBb:\wL\circ\tp\longrightarrow\tp'$ (again with an abuse of notation) by
    \begin{equation*}
        \sBb(\ell\otimes x_1\cdots x_k):=\frac{1}{k!}\Bb(x_1\cdots x_k).
    \end{equation*}
    Then, still with the same abuse of notation we have
    \begin{equation*}
        \fU\sBb(\ell\otimes x_1\cdots x_k)=\sum_\ell\sBb(\ell\otimes x_1\cdots x_k)=\frac{1}{k!}\sum_\ell\Bb(x_1\cdots x_k)=\Bb(x_1\cdots x_k),
    \end{equation*}
    where we have used the definition of $\fU$ and the fact that there are $k!$ total orders on a set with $k$ elements. Thus $\fU(\wL\circ\tp,\sBb)=(\wE\circ\tp,\Bb)$ and $\fU$ is dense as claimed.

    \item Let $(\wE\circ\tp,\Bb)$ and $(\wE\circ\tq,\Cc)$ be two species with grafting operators and $f:(\wL\circ\tp,\fO\Bb)\longrightarrow(\wL\circ\tq,\fO\Cc)$ be an order independent morphism of species with ordered graftings. Then there is a map $\tilde f:\wE\circ\tp\longrightarrow\wE\circ\tq$ such that $f=id_\wL\times\tilde f$. Then since $f$ is a morphism of species with ordered graftings we have, with the conventions described at the beginning of this proof, and still some abuse of notations:
    \begin{equation*}
        f_{\{I^+\}}\left((1)\otimes\fO\Bb(\ell\otimes x_1\cdots x_k)\right) = (1)\otimes \fO\Cc\left(f_\pi(\ell\otimes x_1\cdots x_k)\right).
    \end{equation*}
    Then using the fact that $f$ is order preserving and the definition of $\fO f$ we obtain
    \begin{equation*}
        \tilde f_{\{I^+\}}(\Bb(\ell\otimes x_1\cdots x_k)) = \Cc(\tilde f_\pi(\ell\otimes x_1\cdots x_k)).
    \end{equation*}
    Thus $\tilde f:(\wE\circ\tp,\Bb)\longrightarrow(\wE\circ\tq,\Cc)$ is a morphism of species with grafting operators. Since $\fO\tilde f=f$, the functor $\fO$ is full.

    Now, let again $(\wE\circ\tp,\Bb)$ and $(\wE\circ\tq,\Cc)$ be two species with grafting operators and $f:(\wL\circ\tp,\fO\Bb)\longrightarrow(\wL\circ\tq,\fO\Cc)$ be an order independent morphism of species with ordered graftings. Assume there are two morphisms of species with graftings $g_1,g_2:(\wE\circ\tp,\Bb)\longrightarrow(\wE\circ\tq,\Cc)$ such that $\fO g_1=\fO g_2=f$. Then by hypothesis and definition of $\fO$ we have
    \begin{align*}
        \fO g_1(\ell\otimes x_1\cdots x_k)=\ell\otimes g_1(x_1\otimes x_k) = f(\ell\otimes x_1\cdots x_k)& =\fO g_2(\ell\otimes x_1\cdots x_k) \\ & =\ell\otimes g_2(x_1\otimes x_k).
    \end{align*}
    Therefore $g_1=g_2$ and $\fO$ is faithfull.
\end{enumerate}
\end{proof}
As pointed out before, the same type of constructions for cograftings allows one to build functors between these categories as well. Then these functors have the same properties. We omit these results in order to shorten the paper since they would not have major impacts on the rest of the paper.

\subsection{Examples of species with graftings and cograftings}\label{subsec:examples_graftings}

We give some examples of branching and cobranching maps, to make clear that these are naturally occurring objects.

\subsubsection{Grafting of rooted trees} \label{subsubsec:grafting_trees}

Recall from \S \ref{subsubsec:trees} that $\taa$ is the species of rooted trees. Then $\wE\circ\taa$ is the species of rooted forests. Define $\Bb:\wE\circ\taa\longrightarrow\taa'$ the map that, for any finite set $I$ and any $\pi\vdash I$, sends a forest $t_1\cdots t_k\in\taa(\pi)$ to the tree $\Bb(t_1\cdots t_k)\in\taa'[I]$ whose root is $*_I$ and is grafted to each of the roots of the trees $t_i$. 

Notice that this construction is the species version of the celebrated grafting of rooted trees, which plays important roles, for example in algebraic combinatorics (see for example \cite{FoissyBerlin}) and stochastic analysis (see \cite{hairer2015geometric}). This object is our initial motivation for this work, and we will developp its properties further in the reminder of this paper.

\medskip 

For now, notice that $\taa'$ admits a natural monoidal structure. For any finite $I$, any $S\sqcup T=I$, and any $t_1\otimes t_2\in\taa'[S]\otimes\taa'[T]$ set
\begin{equation*}
    t_1\vee t_2 = \begin{cases}
        \Bb(f_1f_2)\quad\text{ if }t_1=\Bb(f_1)\text{ and }t_2=\Bb(f_2) \\
        0\quad\text{ otherwise.}
    \end{cases}
\end{equation*}
This endows $\taa'$ with a monoidal structure and $\Bb$ is a monoidal grafting map for this structure. The same constructions can be made for planar trees and forests, and the resulting ordered grafting map would also be monoidal.

\subsubsection{Graftings for partitions and compositions}

    The species of partitions is $\tPi = \wE \circ \wE_+$. As in~\cite{AM2013}, we write $\wE[I]=\KK*_I$. There is a natural grafting map $\Bb$ for $\wE_+$ obtained by sending every partition $\pi$ of a finite set $I$ to $*_{I^+}$. 

    Recall that $\wE$ is endowed with a trivial monoid structure. Then $\wE_+'$ inherits this structure by setting, for any finite set $I$ and any $S\sqcup T=I$
    \begin{align*}
        \mu_{S,T}:\wE'_+[S]\otimes\wE'_+[T] & \longrightarrow\wE_+'[I] \\
        *_{S^+}\otimes*_{T^+} & \longmapsto *_{I^+}.
    \end{align*}
    $\tPi$ also has a monoid structure, given by the union. Indeed, if $S\sqcup T=I$, $\pi_1\vdash S$, and $\pi_2\vdash T$, then $\pi_1\cup\pi_2\vdash I$. Then the grafting map for $\wE_+$ described above is monoidal.

\medskip

Similarly, we can also endow $\wE_+$ with a monoidal ordered grafting. Recall that $\tSig = \wL \circ \wE_+$ is the species of set compositions. It admits the following monoid structure: for $S\sqcup T$, $F_1=B_1|\cdots|B_k\vDash S$, and $F_2=C_1|\cdots|C_\ell\vDash T$, define:
\begin{align*}
    \nu_{S,T}:& \tSig[S]\otimes\tSig[T]  \longrightarrow\tSig[I] \\
    & F_1\otimes F_2  \longmapsto B_1|\cdots|B_k|C_1|\cdots|C_\ell.
\end{align*}
Then the ordered grafting map $\Bb:\tSig\longrightarrow\wE_+'$ on $\wE_+$ which sends any composition $F\vDash I$ to $*_{I^+}$ is monoidal.

\subsubsection{Grafting for permutations} \label{subsubsec:perm_grafting}

It is a standard result of group theory that any permutation can be written as a product of disjoint cyclic permutations. 
    Therefore, the species of permutations ${\bf \mathfrak{S}}$ can be written as ${\bf \mathfrak{S}}= \wE \circ \mathbf{C}$, where $\mathbf{C}$ is the species of cycles. Given a finite set $I$, every cyclic permutation on $I^+$ can be identified in a natural manner with a linear order on $I$, by forgetting the element $*_I$ of the cycle: the first element of the new linear order is the element following $*_I$ on the cycle permutation (see Section 1.4, Example 7 on \cite{BLL}). Therefore, $\mathbf{C}'=\wL$. 
    
    \medskip
    
    Given a cycle $(i_1 \, i_2 \, \cdots \, i_k)$ formed with the elements of a finite set $I$, we can form $k$ linear orders on $I$ by selecting one element $i_j$ of the cycle and listing the next elements from the cyclic order. In particular, the last element of the linear order is the element preceding $i_j$ on the cycle $(i_1 \, i_2 \, \cdots \, i_k)$. Hence, given a permutation in ${\bf \mathfrak{S}}[I]$ with $n$ cycles, we define a map ${\bf \mathfrak{S}} \to \wL$ by summing over all linear orders $\ell$ of $I$ such that every linear order formed from one cycle of the permutation is a consecutive subword of $\ell$. For instance, if $I = \{a,b,c,d\}$, the image of the permutation $(a \, b)(c \, d)$ is given by
    \[(a \, b)(c \, d) \mapsto abcd + abdc + bacd + badc + cdab + dcab+ cdba + dcba.\]
    This is a grafting map $\Bb:\wE\circ\mathbf{C}\to\mathbf{C}'$ on cyclic permutations is not monoidal, neither with the concatenation operation nor the shuffle of lists on $\wL$.

\subsubsection{Grafting of connected graphs} \label{subsubsec:conn_graph}

    The species of graphs is given by $\mathbf{G}=\wE \circ \mathbf{G}_{\mathsf{c}}$, where $\mathbf{G}_{\mathsf{c}}$ is the species of connected graphs. Given a graph with vertex set $I$, and connected components $g_1, \hdots, g_n$, we construct a graph with vertex set $I^+$ by attaching every vertex of every connected graphs of $g$ to $*_I$. This defines a grafting $\wE\circ\mathbf{G} \to \mathbf{G}_{\mathsf{c}}'$. 

Let us now endow $\tcG'$ with a structure of monoid. For any finite set $I$, any $S\sqcup T=I$, and any $g\otimes h=\tcG'[S]\otimes\tcG'[T]$ define $g\Join h$ to be the graph obtain by merging $g$ and $h$ on their vertices decorated by $*_S$ and $*_T$ respectively. This defines a species map $\mu:\tcG'.\tcG'\longrightarrow\tcG'$ which is clearly associative. Then the grafting defined above is monoidal for this structure of monoid on $\tcG'$.
    

\subsubsection{Cograftings for rooted trees} \label{ex:cografting_trees} 
    We can endows the species of rooted trees $\tq = \mathbf{a}$ with a cografting map as follows. Let $I$ be a finite set. Given a tree $t \in \mathbf{a}'[I]$, we set 
\begin{align*}
    \bB(t) = \begin{cases} 
              & f \quad\text{ if }t=\Bb(f) \\
              & 0 \quad \text{ otherwise.}
            \end{cases}
\end{align*}
In particular, $\bB$ vanishes on any tree $t\in\mathbf{a}'[I]$ where the $*_I$ is not on the root. $\Bb$ is a monoidal cografting for the monoidal product of $\mathbf{a}'$ defined in \S\ref{subsubsec:grafting_trees}.

\medskip
    
One can define other cobranching on rooted trees. For example, let $\{s_1, \hdots, s_k\}$ be the (possibly empty) set of trees grafted to the vertex in $t$ labelled by $*_I$. Then, $\widetilde{\bB}(t)$ is defined as the forest 
    \[\widetilde{\bB}(t):= t' s_1  \cdots  s_k,\]
    where $t'$ is the (possibly empty) tree obtained from $t$ by erasing all the elements in $t$ greater than $*_I$, when considering the vertices of $t$ as a partial ordered set with the order $\leq_t$. Then $\widetilde\bB:\taa'\longrightarrow\wE\circ\taa$ is a cobranching map, however it is not monoidal for the monoidal structure on $\taa'$ described in \S \ref{subsubsec:grafting_trees}.

\subsection{Grafting maps from operads} \label{subsec:operads}

Let $\tq$ be a positive species. Recall from \S \ref{FreeCMonoid} that the \emph{free monoid on $\tq$} is given by the species $\wE\circ\tq$ endowed with the (commutative) concatenation. It has the following freeness property: given a monoid $\tp$ and a species map $\zeta: \tq \to \tp_+$, then $\zeta$ extends uniquely to a monoid map:
\begin{equation} \label{eq:univ_prop_E}
\xymatrix{
(\wE\circ\tq)_+\ar[rr]^{\widehat{\zeta}}
    & &\tp_+\\
    &\tq \ar[lu] \ar[ru]_{\zeta}&
    },
\end{equation}
where the unlabeled arrow corresponds to the natural inclusion. Here, $\tp_+$ denotes the positive part associated to the species $\tp$, see equation \eqref{eq:pos_part}.

\medskip

From the universality of the monoid $\wE\circ\tq$, a sufficient condition for the existence of a (necessarily unique) monoidal grafting map $\Bb: \wE\circ\tq \to \tq'$ is given by:
\begin{itemize}
    \item[$(i)$] a monoidal structure on $\tq'$, and
    \item[$(ii)$] a species map $\zeta: \tq \to \tq'$.
\end{itemize}
Items $(i)$ and $(ii)$ leads to the following definition.
\begin{definition} \label{def:induced_grafting}
Let $\tq$ be a positive species. Given a monoid structure on $\tq'$ and a species map $\zeta: \tq \to \tq'$, the \Def{induced grafting map from $\zeta$}, $\Bb_\zeta: \wE\circ\tq \to \tq'$, is the unique monoid map induced by $\widehat{\zeta}$.
\end{definition}

Notice that $\Bb_\zeta: \wE\circ\tq \to \tq'$ depends on $\zeta$ and the monoid structure on $\tq'$. We don't indicate the latter one,  which will be clear from the context.

\medskip

We can always consider the trivial monoid structure in $\tq'$ in order to apply the previous definition. A non-trivial structure on $\tq'$ coming from $\tq$ is given by M\'endez in \cite{Mendez}, using the notion of \emph{symmetric operad}.

\medbreak

A \Def{positive (symmetric) operad} on $\tq$ consists of a positive species $\tq$ together with maps
\[\gamma: \tq\circ\tq \to \tq \quad \text{ and } \quad \eta: \wX\to \tq\]
with $\wX$ the species of singletons of equation~\eqref{eq:species_sing} satisfying associativity and unitality axioms (see \cite{AM2010}, Appendix B.1.2 for precise axioms). When $(\tq, \gamma, \eta)$ is a positive operad, the chain rule from $\gamma$ and $\eta$ gives the maps
\[\gamma': (\tq'\circ \tq)\cdot \tq' \to \tq' \quad \text{ and } \eta': 1 \to \tq'.\]
The product map $\mu: \tq' \cdot \tq' \to \tq'$ on $\tq'$ is obtained by composing the maps
\[
\xymatrix{
\mu: \tq' \cdot \tq' = (\tq' \circ \wX)\cdot \tq' \ar[r]^-{\tq'(\eta)\cdot \text{id} }& (\tq'\circ \tq)\cdot \tq' \ar[r]^-{\gamma'}& \tq'.
}
\]
The product $\mu$ and the unit $\eta': 1 \to \tq'$ endows $\tq'$ with a monoidal structure.
\begin{lemma}[M\'endez, \cite{Mendez} Proposition 3.13] \label{lem:Mendez}
If $\tq$ is a positive operad, then $\tq'$ is a positive monoid.
\end{lemma}
This result, together with the freeness of $\wE\circ\tq$, leads to the following sufficient condition for the existence of a grafting map on a positive operad $\tq$.
\begin{theorem} \label{thm:operad_grafting}
Let $(\tq, \gamma, \eta)$ be a positive operad. Given an up operator $u: \tq \to \tq'$ and a monoidal structure on $\tq'$ induced by $\gamma$, there exists a unique monoidal grafting map $\Bb_u: \wE\circ\tq\to \tq'$ such that the following diagram commutes
\begin{equation*}
\xymatrix{
(\wE\circ\tq)_+\ar[rr]^{\Bb_u}
    & &\tq'\\
    &\tq \ar[lu] \ar[ru]_{u}&
    }.
\end{equation*}
\end{theorem}
\begin{proof}
    This result follows directly from applying the diagram \eqref{eq:univ_prop_E} to the monoid on $\tq'$ constructed above and taking $\zeta=u$. In particular, the universal property of $\wE\circ\tq$ implies that $\Bb_u$ is a map of monoids, thus is a monoidal grafting, as well as its unicity.
\end{proof}

\subsection{Graftings of rooted trees from the NAP operad}

An universal result constructed from the natural grafting map on trees $\Bb: \wE\circ\taa\to \taa'$ introduced in \S\ref{subsubsec:grafting_trees} will be shown in theorem~\ref{Universality}. In the next lines we will show that this grafting map is induced from the monoidal structure on $\tq'$ coming from a particular operad structure on trees.

\medskip

There is an operad on the species of rooted trees $\taa$ given as follows: any element of $(\taa \circ \taa)[I]$ is given by a partition $\pi= \{B_1, \hdots, B_k\}$ of $I$, a rooted tree with $k$ vertices labeled by the elements of $\pi$, and $k$ rooted trees labelled by the elements of the blocks of $\pi$. Pictorially, we can represent such data as a rooted tree with $k$ vertices (corresponding to the tree in $\taa[\pi]$) where every vertex contains a rooted tree, each of them labeled with the elements of the blocks $B_1, \hdots, B_k$. Then the composition map $\gamma:\taa\circ\taa\to\taa$ replaces each vertex with its decoration, and branches the roots of its children on the root of the decoration.

\medbreak

This is the \Def{NAP operad}, where NAP stands for Non-Associative Permutative. It was introduced by Livernet in \cite{livernet2006rigidity}, (see also \cite{chapoton2007operads} and~\cite{mendez2014antipode} for this operad in the context of species).

\medbreak

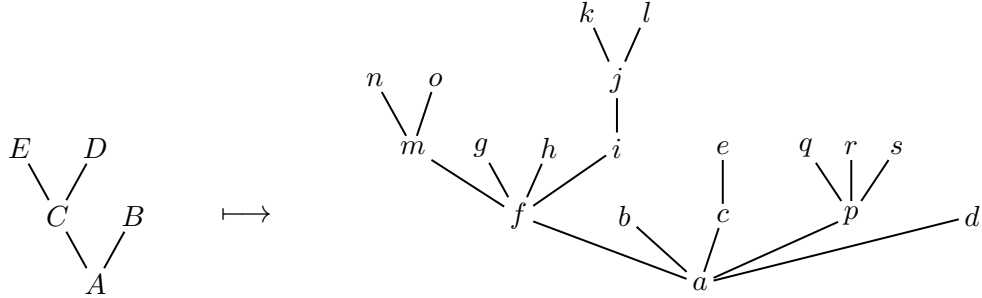
\begin{figure}[hh]

\begin{equation*}
\begin{tikzpicture}[
    x=1cm,
    y=0.9cm,
    baseline=(current bounding box.center),
    vertex/.style={
        inner sep=2pt,
        font=\small
    },
    edge/.style={
        draw=black,
        line width=0.75pt
    }
]

\node[vertex] (A) at (0,0) {$A$};

\node[vertex] (B) at (0.5,1) {$B$};
\node[vertex] (C) at (-0.5,1) {$C$};

\node[vertex] (D) at (0,2) {$D$};
\node[vertex] (E) at (-1,2) {$E$};

\draw[edge]
    (A)--(B)
    (A)--(C)
    (C)--(D)
    (C)--(E);

\node at (2,1) {$\longmapsto$};

\node[vertex] (a) at (8,0) {$a$};
\node[vertex] (f) at (5.6,1) {$f$};
\node[vertex] (b) at (7,1) {$b$};
\node[vertex] (c) at (8.3,1) {$c$};
\node[vertex] (p) at (10,1) {$p$};
\node[vertex] (d) at (11.6,1) {$d$};
\node[vertex] (m) at (4.2,2) {$m$};
\node[vertex] (g) at (5.1,2) {$g$};
\node[vertex] (h) at (6,2) {$h$};
\node[vertex] (i) at (6.9,2) {$i$};
\node[vertex] (e) at (8.3,2) {$e$};
\node[vertex] (q) at (9.4,2) {$q$};
\node[vertex] (r) at (10,2) {$r$};
\node[vertex] (s) at (10.6,2) {$s$};

\node[vertex] (n) at (3.7,3) {$n$};
\node[vertex] (o) at (4.5,3) {$o$};

\node[vertex] (j) at (6.9,3) {$j$};

\node[vertex] (k) at (6.5,4) {$k$};
\node[vertex] (l) at (7.3,4) {$l$};

\draw[edge]
    (a)--(f)
    (a)--(b)
    (a)--(c)
    (a)--(p)
    (a)--(d)

    (f)--(m)
    (f)--(g)
    (f)--(h)
    (f)--(i)

    (m)--(n)
    (m)--(o)

    (i)--(j)
    (j)--(k)
    (j)--(l)

    (c)--(e)

    (p)--(q)
    (p)--(r)
    (p)--(s);

\end{tikzpicture}
\end{equation*}

    \caption{Example of operad on rooted trees $\taa$}
    \label{f1}
\end{figure}

For example, defining the trees
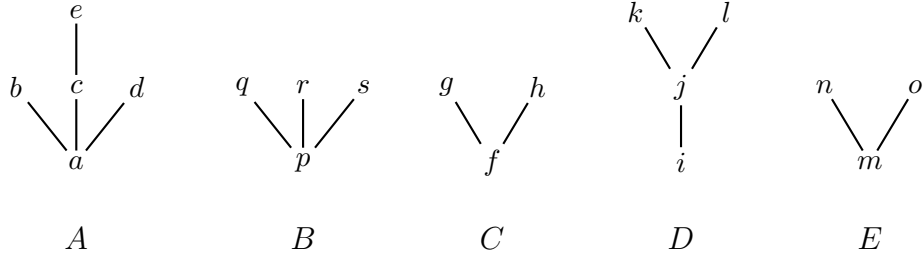
\begin{equation*}
\begin{tikzpicture}[
    baseline=(current bounding box.center),
    vertex/.style={
        inner sep=2pt,
        font=\small
    },
    edge/.style={
        draw=black,
        line width=0.8pt
    },
    treelabel/.style={
        font=\normalsize
    }
]

\begin{scope}[shift={(0,0)}]
    \node[vertex] (a) at (0,0) {$a$};
    \node[vertex] (b) at (-0.8,1) {$b$};
    \node[vertex] (c) at (0,1) {$c$};
    \node[vertex] (d) at (0.8,1) {$d$};
    \node[vertex] (e) at (0,2) {$e$};

    \draw[edge]
        (a)--(b)
        (a)--(c)
        (a)--(d)
        (c)--(e);

    \node[treelabel] at (0,-1) {$A$};
\end{scope}

\begin{scope}[shift={(3,0)}]
    \node[vertex] (p) at (0,0) {$p$};
    \node[vertex] (q) at (-0.8,1) {$q$};
    \node[vertex] (r) at (0,1) {$r$};
    \node[vertex] (s) at (0.8,1) {$s$};

    \draw[edge]
        (p)--(q)
        (p)--(r)
        (p)--(s);

    \node[treelabel] at (0,-1) {$B$};
\end{scope}

\begin{scope}[shift={(5.5,0)}]
    \node[vertex] (f) at (0,0) {$f$};
    \node[vertex] (g) at (-0.6,1) {$g$};
    \node[vertex] (h) at (0.6,1) {$h$};

    \draw[edge]
        (f)--(g)
        (f)--(h);

    \node[treelabel] at (0,-1) {$C$};
\end{scope}

\begin{scope}[shift={(8,0)}]
    \node[vertex] (i) at (0,0) {$i$};
    \node[vertex] (j) at (0,1) {$j$};
    \node[vertex] (k) at (-0.6,2) {$k$};
    \node[vertex] (l) at (0.6,2) {$l$};

    \draw[edge]
        (i)--(j)
        (j)--(k)
        (j)--(l);

    \node[treelabel] at (0,-1) {$D$};
\end{scope}

\begin{scope}[shift={(10.5,0)}]
    \node[vertex] (m) at (0,0) {$m$};
    \node[vertex] (n) at (-0.6,1) {$n$};
    \node[vertex] (o) at (0.6,1) {$o$};

    \draw[edge]
        (m)--(n)
        (m)--(o);

    \node[treelabel] at (0,-1) {$E$};
\end{scope}

\end{tikzpicture}
\end{equation*}

we get the composition of this operad for these trees represented in figure~\ref{f1}.

\medbreak

The monoidal structure on $\taa'$ is illustrated in the following figure~\ref{f2}. Given $S= \{a,b,c,d,e,f,g\}$ and $T= \{x,y,z,w\}$, the product of two elements $t_1 \in \taa'[S]$ and $t_2 \in \taa'[T]$ is obtained by applying the operad $\gamma$ of $\taa$ locally on the ghost-label of $t_1$, replacing this label by the root of $t_2$. All successors of the ghost in $t_1$ become succesors of the root of $t_2$. 
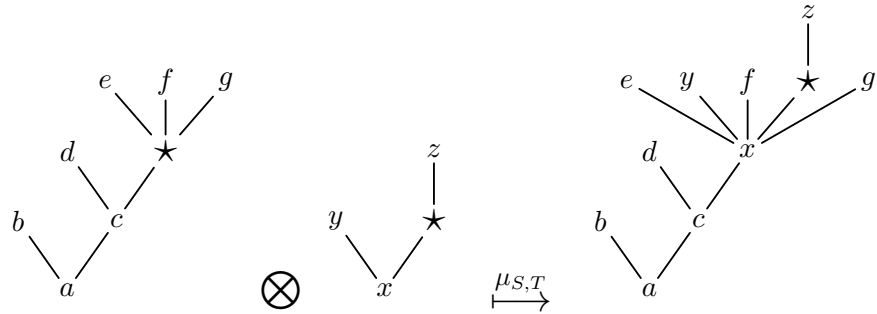
\begin{figure}[hh]

\begin{equation*}
\begin{tikzpicture}[
    x=1cm, y=0.92cm,
    baseline=(current bounding box.center),
    vertex/.style={
        font=\small,
        inner sep=2pt
    },
    edge/.style={
        draw=black,
        line width=0.7pt,
        line cap=round
    },
    symbol/.style={
        font=\large,
        inner sep=0pt
    }
]

\begin{scope}[shift={(0,0)}]
    \node[vertex] (a1) at (0,0) {$a$};
    \node[vertex] (b1) at (-0.65,1) {$b$};
    \node[vertex] (c1) at (0.65,1) {$c$};

    \node[vertex] (d1) at (0,2) {$d$};
    \node[vertex] (t1) at (1.30,2) {\Large$\star$};

    \node[vertex] (e1) at (0.50,3) {$e$};
    \node[vertex] (f1) at (1.30,3) {$f$};
    \node[vertex] (g1) at (2.10,3) {$g$};

    \draw[edge]
        (a1)--(b1)
        (a1)--(c1)
        (c1)--(d1)
        (c1)--(t1)
        (t1)--(e1)
        (t1)--(f1)
        (t1)--(g1);
\end{scope}

\node[symbol] at (2.825,0) {$\bigotimes$};

\begin{scope}[shift={(4.20,0)}]
    \node[vertex] (x2) at (0,0) {$x$};

    \node[vertex] (y2) at (-0.65,1) {$y$};
    \node[vertex] (z2) at (0.65,1) {\Large$\star$};

    \node[vertex] (t2) at (0.65,2) {$z$};

    \draw[edge]
        (x2)--(y2)
        (x2)--(z2)
        (z2)--(t2);
\end{scope}

\node[symbol] at (6.00,0) {$\overset{\mu_{S,T}}{\longmapsto}$};

\begin{scope}[shift={(7.70,0)}]
    \node[vertex] (a3) at (0,0) {$a$};

    \node[vertex] (b3) at (-0.65,1) {$b$};
    \node[vertex] (c3) at (0.65,1) {$c$};

    \node[vertex] (d3) at (0,2) {$d$};
    \node[vertex] (x3) at (1.30,2) {$x$};

    \node[vertex] (e3) at (-0.30,3) {$e$};
    \node[vertex] (y3) at (0.50,3) {$y$};
    \node[vertex] (f3) at (1.30,3) {$f$};
    \node[vertex] (z3) at (2.10,3) {\Large$\star$};
    \node[vertex] (g3) at (2.90,3) {$g$};

    \node[vertex] (t3) at (2.10,4) {$z$};

    \draw[edge]
        (a3)--(b3)
        (a3)--(c3)
        (c3)--(d3)
        (c3)--(x3)
        (x3)--(e3)
        (x3)--(y3)
        (x3)--(f3)
        (x3)--(z3)
        (x3)--(g3)
        (z3)--(t3);
\end{scope}

\end{tikzpicture}
\end{equation*}

    \caption{Example of monoid on derivative of rooted trees $\taa'$}
    \label{f2}
\end{figure}

\medskip

Now, from the previous theorem, this operad on rooted trees implies the existence of a branching map of trees provided we also have an up operator on trees. It turns out that the branching map on trees of introduced before comes from this operad and the simplest non-trivial up operator on trees. To prove this proposition, it is enough to check that the map $\widehat\zeta$ of equation~\eqref{eq:univ_prop_E} coincides with the grafting on rooted trees of \S\ref{subsubsec:grafting_trees}. We leave this simple checking to the reader.
\begin{prop} \label{prop:grafting_from_NAP}
Let $\gamma: \taa \circ \taa \to \taa$ be the operad on trees. Consider the map
\[\zeta: \taa \to \taa'\]
such that any tree $t$ is sent to a new tree obtained by attaching a new root to the root of $t$, labeled by a ghost. Then, the grafting map $\Bb: \wE\circ\taa\to \taa'$ of \S \ref{subsubsec:grafting_trees} is $\Bb_\zeta$. Moreover, $\Bb$ is unique in the sense of the previous theorem.
\end{prop}
In other words, the nature of the classical map $\Bb$  can be understood from the operad NAP and a simple up operator.

\section{Up/down operators and branching maps} \label{sec:up_down}

\subsection{Up and down operators induced from grafting and cografting} \label{subsec:up_down_grafting}

Grafting and cografting maps are closely related to up and down operators, as we now detail.
\begin{definition} \label{def:up_dow_ensembliste}
Let $\tp$ be a species with grafting map $\Bb: \wE \circ \tp \to \tp'$. The \Def{up operator associated to $\Bb$} is the up operator $u_{\Bb}: \wE \circ \tp \to (\wE \circ \tp)'$ given by the composition of the following maps
\begin{equation}  \label{def:up_down_ensemblistes}
    u_{\Bb}: \wE \circ \tp  \xrightarrow{\Bb} \tp' \hookrightarrow  (\wE \circ \tp)'.
\end{equation}
The second arrow refers to the derivative of the natural inclusion $\tp \hookrightarrow  (\wE \circ \tp)$. 

\medbreak

Let $\tp$ be a species with cografting map $\bB:  \tp' \to \wE \circ \tp$. The \Def{down operator associated to $\bB$} is the down operator $d_{\tiny\bB}: (\wE \circ \tp)' \to \wE \circ \tp$ given by the composition of the following maps
\begin{equation} \label{def:down_ensemblistes}
    d_{\tiny\bB}: (\wE \circ \tp)'  	\twoheadrightarrow  \tp' \xrightarrow{\tiny\bB} \wE \circ \tp.
\end{equation}
The first arrow above is the derivative of the natural projection $\wE \circ \tp \twoheadrightarrow \tp$.
\end{definition}
The structures of $\wE\circ\tp$ allow us to build other up and down operators from grafting and cografting maps. These other operators are less convenient to manipulate but they have more algebraic properties.
\begin{prop} \label{prop:alg_up_coderiv}  
    Let $(\tq, \Bb)$ be a species with grafting map. The map
    \begin{equation}
        u: \wE\circ\tq \xlongrightarrow{\Delta} \wE\circ\tq \cdot \wE\circ\tq \xlongrightarrow{\text{id} \, \cdot \,\Bb 
         } \wE\circ\tq \cdot  \tq'= (\wE\circ\tq)'
    \end{equation}
    endows the comonoid $(\wE\circ\tq, \Delta)$ with an up coderivation structure.
\end{prop}
\begin{proof}
 We need to show that Diagram \ref{diag:up_coderivation} commutes.
 Then recall that for any finite set $I$, any $\pi\vdash I$ and any
 \begin{equation*}
  s=s_1\cdots s_k\in \tq(\pi)\subseteq \wE\circ\tq[I]
 \end{equation*}
 we have
 \begin{equation*}
  \Delta_I(s)=\sum_{J\subseteq[k]}s_J\otimes s_{[k]\setminus J}
 \end{equation*}
 with $s_J:=\otimes_{j\in J}s_j$ (recall that in this case, we use $\otimes$ to denote a commutative tensor product, as explained in \S \ref{subsubsec:comp_part}). Using
 \begin{equation*}
  u(s_J)=\sum_{\tilde J\subseteq J}s_{\tilde J}\otimes \Bb(s_{J\setminus \tilde J})
 \end{equation*}
 we obtain
 \begin{equation} \label{eq:half_diag_preuve}
  (\text{id}\otimes u+u\otimes\text{id})\Delta_I(s) = \sum_{J\subseteq[k]}\left[\left(
  \sum_{\tilde J\subseteq J}s_{\tilde J}\otimes \Bb(s_{J\setminus \tilde J})
  \right)
  \otimes s_{[k]\setminus J} + s_J\otimes\sum_{\tilde J\subseteq [k]\setminus J}s_{\tilde J}\otimes \Bb(s_{([k]\setminus J)\setminus \tilde J})
  \right].
 \end{equation}
 On the other hand we have
 \begin{align*}
  \Delta'_{J}(s_J\otimes x_*) & = \Delta_{J\cup\{*\}}(s_J\odot x_*) \\
  & = \sum_{\tilde J\subseteq J\cup\{*\}}s_{\tilde J}\otimes s_{(J\cup\{*\})\setminus J} \\
  & = \sum_{\tilde J\subseteq J}\left(s_{\tilde J}\otimes s_{(J\setminus\tilde J)\cup\{*\})} + s_{\tilde J\cup\{*\}}\otimes s_{J\setminus\tilde J}\right) \\
  & = \sum_{\tilde J\subseteq J}\left(s_{\tilde J}\otimes s_{J\setminus\tilde J}\otimes x_* + s_{\tilde J}\otimes x_*\otimes s_{J\setminus\tilde J}\right)
 \end{align*}
 where we used the obvious notation $s_{\tilde J\cup\{*\}}:=s_{\tilde J}\otimes x_*$ for any $\tilde J\subseteq J$. Then 
 using $u(s)=\sum_{J\subseteq [k]}s_{J}\otimes \Bb(s_{[k]\setminus J})$ we obtain by setting $x_*=\Bb(s_{[k]\setminus J})$:
 \begin{equation} \label{eq:half_diag_preuve2}
  \Delta'_I(u(s)) = \sum_{\tilde J\subseteq J\subseteq[k]}\left(s_{\tilde J}\otimes s_{J\setminus\tilde J}\otimes \Bb(s_{[k]\setminus J}) + s_{\tilde J}\otimes \Bb(s_{[k]\setminus J})\otimes s_{J\setminus\tilde J}\right).
 \end{equation}
 This is actually equal to the R.H.S. of Equation \eqref{eq:half_diag_preuve}. To see this, let us define relabel the sums in the second term of the  R.H.S. of Equation \eqref{eq:half_diag_preuve}. Let us set $J':=\tilde J\cup J$ and $\tilde J':=J\subseteq J'$. Then the second term of the R.H.S. of Equation \eqref{eq:half_diag_preuve} becomes exactly the first term of the R.H.S. of Equation \eqref{eq:half_diag_preuve2}.
 
 Finally, let us relabel the sums in the first term of the  R.H.S. of Equation \eqref{eq:half_diag_preuve}. Let us set $J':=([k]\setminus J)\cup\tilde J$ and $\tilde J':=\tilde J\subseteq J'$. Then the first term of the R.H.S. of Equation \eqref{eq:half_diag_preuve} becomes exactly the second term of the R.H.S. of Equation \eqref{eq:half_diag_preuve2}.
\end{proof}
With a very similar proof which we omit we find the following dual result. 
\begin{prop} \label{prop:alg_down_deriv}  
    Let $(\tp,\bB)$ be a species with a cografting. The map
    \begin{equation*}
        d: (\wE\circ\tq)'=\wE\circ\tq\cdot\tq' \xlongrightarrow{\text{Id}\cdot\tiny{\bB}} \wE\circ\tq \cdot \wE\circ\tq \xlongrightarrow{\mu}\wE\circ\tq
    \end{equation*}
    endows the monoid $(\wE\circ\tq, \mu)$ with an down derivation structure.
\end{prop}

\subsection{Up operators and graftings: isomorphism and adjointness} \label{subsec:up_down_cat} 

Let $\Moup$ be the category of commutative\footnote{Essentially the same discussion would hold for ordered grafting by working with (noncommutative) monoids.} monoids with up operators. Its objects are triplet $(\tp,\mu,u)$ where $(\tp,\mu)$ is a monoid and $u:\tp\longrightarrow\tp'$ is an up operator on $\tp$. Its morphisms are morphisms of monoids that are also morphisms of up operators. Now, recall from remark~\ref{rk:dec} that $\Spgrdec$ is the category of species with graftings whose morphisms are decomposable, where $f:\wE\circ\tp\longrightarrow\wE\circ\tq$ is decomposable if it exists $\bar f:\tp\longrightarrow\tq$ such that $f=\text{Id}\circ\bar f$.

\medbreak 

Then the construction above defines a functor
\begin{equation*}
    \fH:\Spgrdec\longrightarrow\Moup
\end{equation*}
characterised by $\fH(\tp,\Bb)=(\wE\circ\tp,\mu,u_\Bb)$ (where $\mu$ is the concatenation product of $\wE\circ\tp$) and $\fH(f)=f$; which is a morphism of monoids for the concatenation products since $f$ is decomposable.

\medskip

Let us explore more in details the link between categories of species with up operators and graftings. First observe that if $(\tp,\mu,u)$ is a monoid with up operator, then there is a natural associated branching map $\Bb_u:\wE\circ\tp\longrightarrow\tp'$ on $\tp$ given by, for any finite set $I$ and any $\pi\vdash I$:
\begin{align*}
    (\Bb_u)_\pi: \tp(\pi) & \longrightarrow \tq'[I] \\
    x_1\cdots x_k & \longmapsto u_I(\mu_\pi(x_1\cdots x_k)).
\end{align*}
Furthermore, for $f:(\tp,\mu,u)\longrightarrow(\tq,\nu,v)$ morphism of monoids with up operators, let us define $\fF(f):\wE\circ\tp\longrightarrow\wE\circ\tq$, by setting for any finite set $I$ and any $\pi=\{B_1,\cdots,B_k\}\vdash I$
\begin{align*}
    \fF(f)_\pi:\tp(\pi) & \longrightarrow\tq(\pi) \\
    x_1\dots x_k & \longmapsto f_{B_1}(x_1)\cdots f_{B_k}(x_k).
\end{align*}
Finally, let $\Mogr$ be the category of monoids with up operators, and let $\fU:\Mogr\longrightarrow\Spgr$ be the forgetful functor that forgets the monoid structure. Then we have the following results. The second point explains the fact that our branching maps generalize up operators.
\begin{theorem} \label{thm:comparison_up_grafting_NEW}
    \begin{enumerate}
        \item[(i)] Setting $\fF(\tp,\mu,u):=(\tp,\mu,\Bb_u)$, the construction above defines a functor 
        \begin{equation*}
            \fF:\Moup\longrightarrow\Mogr.
        \end{equation*}
        \item[(ii)] The category $\Moup$ is isomorphic, as a category, to a subcategory of $\Mogr$ whose objects are monoids with graftings $(\tp,\mu,\Bb)$ such that 
        \begin{equation} \label{eq:prop_subcategory_objects}
         \forall I, ~\forall\pi\vdash I,~\Bb_\pi=\B_{\{I\}}\circ\mu_\pi;
        \end{equation}
        and whose morphisms are decomposable morphisms of species with graftings $f:\wE\circ\tp\longrightarrow\wE\circ\tq$ which are morphisms of monoids in the following sense:
        \begin{equation} \label{eq:prop_subcat_morphisms}
            \for I,~\forall S\sqcup T=I,~f_{\{I\}}\circ\mu_{S,T} = \nu_{S,T}\circ(f_{\{S\}}\otimes f_{\{T\}})
        \end{equation}
        where $\mu$ is the monoid product of $\tp$ and $\nu$ is the monoid product of $\tq$.
        \item[(iii)] $\fH$ is left adjoint to the functor $\fU\fF:\Moup\longrightarrow\Spgrdec$.
    \end{enumerate}  
\end{theorem}
\begin{proof}
    For the first point, for $f:(\tp,\mu,u)\longrightarrow(\tq,\nu,v)$, it is trivial from the definition that $\fF(f):(\wE\circ\tp,\Bb_u)\longrightarrow(\wE\circ\tq,\Bb_v)$ is a morphism of monoids for the concatenation products. Thus we only need to show that $\fF(f)$ is a morphism of species with up operators.


    For any finite set $I$ and any $\pi=\{B_1,\cdots,B_k\}\vdash I$ we have for any $x_1\cdots x_k\in\tp(\pi)$
    \begin{align*}
        \fF(f)_{\{I^+\}}\left((\Bb_u)_\pi(x_1\cdots x_k)\right) & = \fF(f)_{\{I^+\}}\left(u_I(\mu_\pi(x_1\cdots x_k))\right) \tag{by definition of $\Bb_u$} \\
        & = f_{I^+}\left(u_I(\mu_\pi(x_1\cdots x_k))\right)\tag{by definition of $\fF(f)$} \\
        & = v_I\left(f_{I}(\mu_\pi(x_1\cdots x_k))\right)\tag{since $f$ is a morphism of up operators} \\
        & =  v_I\left(\nu_\pi(f_{B_1}(x_1)\cdots f_{B_k}(x_k))\right)\tag{since $f$ is a morphism of monoids} \\
        & = v_I\left(\nu_\pi(\fF_\pi(x_1\cdots x_k))\right)\tag{by definition of $\fF(f)$} \\
        & = (\Bb_v)_I(\fF_\pi(x_1\cdots x_k))
    \end{align*}
    by definition of $\Bb_v$. Thus $\fF$ is a functor as claimed.

    \medskip 

    For the second point, let $\mathcal{C}$ be the subcategory of $\Mogr$ described in the theorem. Let us define a functor $\fG:\mathcal{C}\longrightarrow\Moup$ such that $\fG\fF=\text{Id}_{\Moup}$ and $\fF\fG=\text{Id}_\mathcal{C}$. For any $(\tq,\mu,\Bb)$ object of $\mathcal{C}$, set $\fG(\tq,\mu,\Bb)=(\tq,\mu,\fG(\Bb))$ with $\fG(\Bb):\tq\longrightarrow\tq'$ is defined for any finite set $I$ by
    \begin{equation*}
        \fG(\Bb)_I:=\Bb_{\{I\}}:\tq(\{I\})=\tq[I]\longrightarrow\tq'[I].
    \end{equation*}
    Furthermore, for any $f:\wE\circ\tp\longrightarrow\wE\circ\tq$ morphism of $\mathcal{C}$, define $\fG(f):\tp\longrightarrow\tq$ by setting, for any finite set $I$
    \begin{equation*}
        \fG(f)_I:=f_{\{I\}}:\tp(\{I\})=\tp[I]\longrightarrow\tq(\{I\})=\tq[I].
    \end{equation*}
    Then $\fG(f)$ is a morphism of monoids since $f$ has property \eqref{eq:prop_subcat_morphisms} by hypothesis. Furthermore, for any finite set $I$ and any $x\in\tp[I]$ we have
    \begin{align*}
        \fG(f)_{I^+}\left(\fG(\Bb)_I(x)\right) & = f_{\{I^+\}}\left(\Bb_{\{I\}}(x)\right)\tag{by definition of $\fG$} \\
        & = \Cc_{\{I\}}\left(f_{\{I\}}(x)\right)\tag{since $f$ is a morphism of species with up operators} \\
        & = \fG(\Cc)_{I}\left(\fG(f)_I(x)\right)
    \end{align*}
    again by definition of $\fG$. Thus $\fG(f)$ is a morphism of species with up operators and $\fG$ is a functor between $\mathcal{C}$ and $\Moup$.

    Now, for any object $(\tp,\mu,\Bb)$ of $\mathcal{C}$ we have 
    \begin{equation*}
        \fF\fG(\tp,\mu,\Bb) = (\tp,\mu,\Bb_{\fG(\Bb)}).
    \end{equation*}
    For any finite set $I$ and any $\pi\vdash I$, by definition of $\fG(\Bb)$ and $\fF$ we have for any $x_1\cdots x_k\in\wE\circ\tp$ 
    \begin{equation*}
        (\Bb_{\fG(\Bb)})_\pi(x_1\cdots x_k) = \fG(\Bb)_I(\mu_\pi(x_1\cdots x_k))=\Bb_{\{I\}}(\mu_\pi(x_1\cdots x_k)) = \Bb_\pi(x_1\cdots x_k)
    \end{equation*}
    since $\Bb$ has property~\eqref{eq:prop_subcategory_objects} by hypothesis. Therefore $\fF\fG(\tp,\mu,\Bb)=(\tp,\mu,\Bb)$ as required. 

    Next, for any morphism $f:\wE\circ\tp\longrightarrow\wE\circ\tq$ of $\mathcal{C}$, we have by construction $\fF\fG(f):\wE\circ\tp\longrightarrow\wE\circ\tq$ and for any finite set $I$, any $\pi=\{B_1,\cdots,B_k\}\vdash I$, and any $x_1\cdots x_k\in\tp(\pi)$:
    \begin{align*}
        \fF\fG(f)(x_1\cdots x_k) & = \fG(f)_{B_1}(x_1)\cdots\fG(f)_{B_k}(x_k) \tag{by definition of $\fF$} \\
        & = f_{\{B_1\}}(x_1)\cdots f_{\{B_k\}}(x_k) \tag{by definition of $\fG$} \\
        & = f_\pi(x_1\cdots x_k)
    \end{align*}
    since $f$ is decomposable by hypothesis. So $\fF\fG(f)=f$ as required and $\fF\fG=\text{Id}_\mathcal{C}$.

    For the other composition let $(\tp,\mu,u)$ be an object of $\Moup$. Then $\fG\fF(\tp,\mu,u)=(\tp,\mu,\fG(\Bb_u))$ and for any finite set $I$ and any $x\in\tp[I]$ we have by definition of $\fG$
    \begin{equation*}
        \fG(\Bb_u)_I(x) = (\Bb_u)_{\{I\}}(x) = u_I(x)
    \end{equation*}
    with the second equality being by definition of $\Bb_u$, thus of $\fF$. Therefore $\fG\fF(\tp,\mu,u)=(\tp,\mu,u)$ as required.

    Finally, let $f:(\tp,\mu,u)\longrightarrow(\tq,\nu,v)$  be a morphism of $\Moup$. Then for any finite set $I$ and any $x\in\tp[I]$ we have by definition of $\fG$
    \begin{equation*}
        \fG\fF(f)_I(x) = \fF(f)_{\{I\}}(x)=f_{I}(x)
    \end{equation*}
    by definition of $\fF$, since $x\in\tp[I]=\tp(\{I\})\subseteq\wE\circ\tp[I]$. Thus $\fG\fF(f)=f$, and $\fG\fF=\text{Id}_{\Moup}$: $\mathcal{C}$ and $\Moup$ are isomorphic as categories as claimed.

    \medskip
    
    For the third point, we use the characterization of adjoint in terms of units and counits (see for example \cite{ML71}). Let us build to natural transformations $\epsilon:\fH\fU\fF\Longrightarrow\text{Id}_{\Moup}$ and $\eta:\text{Id}_{\Spgr_{\rm dec}}\Longrightarrow\fU\fF\fH$. 
    
    For any monoid with up operator $(\tp,\mu,u)$, let us use the short-hand notation $\bar u:=u_{\Bb_u}$. Then $\fH\fU\fF(\tp,\mu,u)=(\wE\circ\tp,\cdot,\bar u)$, with $\cdot$ the concatenation product. Then we define $\epsilon_{(\tp,\mu,u)}:(\wE\circ\tp,\cdot,\bar u)\longrightarrow(\tp,\mu,u)$ as follows: for any finite set $I$, any $\pi\vdash I$ we set
    \begin{equation*}
        (\epsilon_{(\tp,\mu,u)})_\pi(x_1\cdots x_k):=\mu_\pi(x_1\cdots x_k)\quad\forall x_1\cdots x_k\in\tp(\pi).
    \end{equation*}
    Let us check that $\epsilon$ is a morphism of monoids with up operators. First, $\epsilon_{(\tp,\mu,u)}$ is by definition a morphism of monoids. Second, for any finite set $I$, any $\pi\vdash I$, and any $x_1\cdots x_k\in\tp(\pi)$ we have
    \begin{align*}
        (\epsilon_{(\tp,\mu,u)})_{\{I^+\}}\left((\bar u)_\pi (x_1\cdots x_k)\right) & = (\epsilon_{(\tp,\mu,u)})_{\{I^+\}}\left((\Bb_u)_\pi(x_1\cdots x_k)\right)\tag{by definition of $\bar u$} \\
        & = (\Bb_u)_\pi(x_1\cdots x_k)\tag{by definition of $\epsilon$} \\
        & = u_I(\mu_\pi(x_1\cdots x_k))\tag{by definition of $\Bb_u$} \\
        & = u_I((\epsilon_{(\tp,\mu,u)})_\pi(x_1\cdots x_k))
    \end{align*}
    again by definition of $\epsilon$. Therefore $\epsilon_{(\tp,\mu,u)}$ is also a morphism of up operators: it is a morphism of $\Moup$.

    Now, let us show that $\epsilon$ is a natural transformation. Let $f:(\tp,\mu,u)\longrightarrow(\tq,\nu,v)$ be any morphism of $\Moup$. Then by definition of the functors we have
    \begin{equation*}
        \fH\fU\fF(f)=\text{Id}\circ f:(\wE\circ\tp,\cdot,\bar u)\longrightarrow(\wE\circ\tq,\cdot,\bar v).
    \end{equation*}
    Now, for any finite $I$, any $\pi=\{B_1,\cdots,B_k\}\vdash I$ and any $x_1\cdots x_k\in\tp(\pi)$ we have, using the definition of the functors
    \begin{align*}
        (\epsilon_{(\tq,\nu,v)})_\pi\left(\fH\fU\fF(f)_\pi(x_1\cdots x_k)\right) & = (\epsilon_{(\tq,\nu,v)})_\pi\left(f_{B_1}(x_1)\cdots f_{B_k}(x_k)\right) \\
        & = \mu_\pi\left(f_{B_1}(x_1)\cdots f_{B_k}(x_k)\right) \tag{by definition of $\epsilon$} \\
        & = f_I(\mu_\pi(x_1\cdots x_k))
    \end{align*}
    since $f$ is a morphism of $\Moup$ and thus a morphism of monoids. Thus $\epsilon$ is natural as claimed.

    \medskip

    For $\eta$, observe that for any species with up operator $(\tp,\Bb)$ we have $\fU\fF\fH(\tp,\Bb)=(\wE\circ\tp,\bar B)$ with $\bar B:\wE\circ(\wE\circ\tp)\longrightarrow(\wE\circ\tp)'$. One can explicitly write the action of $\bar\Bb$ but we will not need these details here. Instead, notice that for any finite set $I$ we have $\wE\circ\tp[I]=\wE\circ\tp(\{I\})\subseteq\wE\circ\wE\circ\tp[I]$. Then we can define $\eta$ to be this inclusion of $\wE\circ\tp$ into $\wE\circ\wE\circ\tp$: for any finite set $I$, for any $\pi\vdash I$, define
    \begin{align*}
        (\eta_{(\tp,\Bb)})_I:\wE\circ\tp[I] & \longrightarrow\wE\circ\wE\circ\tp[I] \\
        x_1\cdots x_k & \longmapsto X= x_1\cdots x_k
    \end{align*}
    where in the right-hand-side, $x_1\cdots x_k$ is seen as an element of $\wE\circ\tp(\{I\})\subseteq\wE\circ\wE\circ\tp[I]$. Then $\eta$ is trivially a morphism of species with graftings.

    Let us now check that $\eta$ is natural. For any $f:(\tp,\Bb)\longrightarrow(\tq,\Cc)$ we have $f:\wE\circ\tp\longrightarrow\wE\circ\tq$ as a map of species. Then by definition of the functors, for any finite set $I$, any $\pi=\{B_1,\cdots,B_k\}\vdash I$, and any $X_1\cdots X_k\in(\wE\circ\tp)(\pi)$ (i.e. $X_i\in\wE\circ\tp[B_i]$) we have
    \begin{equation*}
        \fU\fF\fH(f)_\pi(X_1\cdots X_k) = f_{B_1}(X_1)\cdots f_{B_k}(X_k)\in\wE\circ\wE\circ\tq[I].
    \end{equation*}
    Then, by definition of $\eta$ we have, for any $x\in\wE\circ\tp[I]$
    \begin{equation*}
        \fU\fF\fH(f)_{\{I\}}\left((\eta_{(\tp,\Bb)})_I(X)\right) = f_I(X) = (\eta_{(\tq,\Cc)})_I(f_I(X)).
    \end{equation*}
    Thus $\eta$ is natural as claimed.

    \medskip
    
    Finally, we need to check the compatibility of $\epsilon$ and $\eta$ with $\fH$ and $\fU\fF$. First, for any finite set $I$ and any $X=x_1\cdots x_k\in\wE\circ\tp[I]=\wE\circ\tp(\{I\})\subseteq\wE\circ\wE\circ\tp[I]$ we have $\fH((\eta_{(\tp,\Bb)})_I)(X)=X$ by definition of $\fH$ and $\eta$. Therefore
    \begin{equation*}
        (\epsilon_{\fH(\tp,\Bb)})_I\circ\fH((\eta_{(\tp,\Bb)})_I)(X) = (\epsilon_{\fH(\tp,\Bb)})_I(X) = X.
    \end{equation*}
    The last equality comes from the fact that $X$ is an element of $\wE\circ\wE\circ\tp[I]$ with only one block so that the action of $\epsilon$ is trivial on this object. Then $\epsilon_{\fH(\tp,\Bb)}\circ\fH(\eta_{(\tp,\Bb)})=\text{Id}_{\fH(\tp,\Bb)}$ as necessary.

    Similarly, $\eta_{\fU\fF(\tp,\mu,u)}:\wE\circ\tp\longrightarrow\wE\circ\wE\circ\tp$ maps an element $X\in\wE\circ\tp$ to a one-block element of $\wE\circ\wE\circ\tp$. Therefore, by definition of $\fF$ and $\epsilon$ we have that the action of $\fU\fF(\epsilon_{(\tp,\mu,u)})$ on such an element is trivial and $\fU\fF(\epsilon_{(\tp,\mu,u)})\circ\eta_{\fU\fF(\tp,\mu,u)}=\text{Id}_{\fF(\tp,\mu,u)}$ as needed.

    All these points together show that $\fH$ is left adjoint to $\fU\fF$ as claimed. 
\end{proof}

\medbreak

\subsection{Balanced up and down operators} \label{subsec:balanced}


Recall that a permutation $\sigma\in \mathfrak{S}_n$ is a bijection $\sigma:[n]\to[n]$. Therefore, for any species $\tp$ we have a linear map $\tp[\sigma]:\tp[n]\to\tp[n]$ by the functoriality of $\tp$. Then for $x\in\tp[n]$ we write $\sigma.x$ instead of $\tp[\sigma](x)$. Furthermore, if $\sigma\in\mathfrak{S}_n$ is the transposition that exchanges $i$ and $j$, we write $\sigma=(i,j)$. We also write $\sigma=(n,n-1,\cdots,1)$ for the cyclic permutation $\sigma\in\mathfrak{S}_n$ given by $\sigma^{-1}(k)=(k+1)\,\text{mod}\, n$ for any $k\in[n]$. These notations were introduced in \cite[Notation 19.23]{AM2010}, reformulating \cite{gucta2002symmetric}.

\subsubsection{Balanced operators from graftings and cograftings}

We can now introduce the notion we will be investigating, from \cite[Definition 19.23]{AM2010}.
\begin{definition} \label{def:up_down_bal}
    Let $\tp$ be a species with up operator $u: \tp \to \tp'$ and down operator $d: \tp' \to \tp$. Then $(u,d)$ is a \Def{pair of balanced up and down operators} if, for any $n\in\N$, the following relations hold
    \begin{align}
   (1,2)\cdot u^2(x) & = u^2(x)\quad\forall x\in\tp[n] \label{eq:up_down1} \\
   d^2((1,2)\cdot x) & = d^2(x)\quad\forall x\in\tp[n+2] \label{eq:up_down2} \\
   d\circ u=\lambda_n & \Id\quad\text{ for some  }\lambda_n\in\mathbb{K} \label{eq:up_down3} \\
   d((k+1,1)\cdot u(x)) & =(k,1)\cdot u(d((k,1)\cdot x))\quad\forall x\in\tp[n],~1\leq k<n.  \label{eq:up_down4}
  \end{align}
\end{definition}
The introduction of these objects is motivated by the observation that together with the Fock functor, one can build pairs of creation-annihilation operators from pairs of balanced up and down operators. Then the canonical commutation relations hold thanks to equations \eqref{eq:up_down1}-\eqref{eq:up_down4}. Notice that a twisted version of the Fock functor gives a fermionic version of creation-annihilation operators. We do not aim here to enter these details and instead refer the reader to \cite[section 19.4]{AM2010}.

\medbreak

Instead, we are interested in producing examples of pairs of balanced up and down operators. Indeed, such pairs seem to be difficult to construct in full without any constraints and \cite{AM2010} contains only essentially trivial examples. Our result follows from the next simple proposition.
\begin{prop} \label{prop:balanced_grafting_cografting}
    Let $\tp$ be a species with grafting map $\Bb: \wE \circ \tp \to \tp'$ and with cografting $\bB: \tq' \to \wE\circ \tq$. Let $u_{\Bb}: \wE \circ \tp \to (\wE \circ \tp)'$ be the up operator on $\wE \circ \tp$ associated to $\Bb$ as in definition \ref{def:up_down_ensemblistes}; and $d_{\tiny\bB}: (\wE \circ \tp)' \to \wE \circ \tp$ be the down operator on $\wE \circ \tp$ associated to $\bB$ as in definition \ref{def:down_ensemblistes}. Then $(u_{\Bb},d_{\tiny\bB})$ is a pair of balanced up and down operators on $\wE\circ\tp$ if, and only if, for any $n\in\N$
     \begin{align}
   (1,2)\cdot \Bb^2(x) & = \Bb^2(x)\quad\forall x\in\tp[n] \label{eq:balanced_1} \\
   \bB^2((1,2)\cdot x) & = \bB^2(x)\quad\forall x\in\tp[n+2] \label{eq:balanced_2}\\
   \bB\circ \Bb=\lambda_n & \Id\quad\text{for some }\lambda_n\in\mathbb{K}\label{eq:balanced_3} \\
   \bB((k+1,1)\cdot \Bb(x)) & =(k,1)\cdot\Bb(\bB((k,1)\cdot x))\quad\forall x\in\tp[n],~1\leq k<n \label{eq:balanced_4}.
  \end{align}
\end{prop}
\begin{proof}
    This follows directly from the definitions \ref{def:up_down_ensemblistes} and \ref{def:down_ensemblistes} of $u_{\Bb}$ and $d_{\tiny\bB}$.
\end{proof}
Notice that in this proposition above we use the same abuse of notations that in definition~\ref{def:up_down_bal}. For example, $\Bb^2(x)$ should be read as $\Bb_{[n+1]}(\Bb_{[n]}(x))$ for $x\in\tp[n]$.

We can now easily built examples of balanced pairs of up and down operators build from grafting and cografting maps in the fashion of definition \ref{def:up_dow_ensembliste}. To illustrate this, we give two such examples below. Looking for examples relevant to quantum field theories is left for future investigations.

\subsubsection{An example on graphs}

Our first example in on the species of connected graphs.
\begin{example} \label{ex:balanced_graphs}
    Recall from \S \ref{subsubsec:conn_graph} that $\tcG$ is the species of connected graphs, and $\wE\circ\tcG$ is the species of non-connected graphs. More specifically, for a any finite $I$, $\tcG[I]$ is spanned by elements of the form $G=(V(G),E(G))$, with $V(G)=I$ and $E(G)$ a finite set.
    
    Further recall (from \S\ref{subsubsec:conn_graph}) that $\Bb:\wE\circ\tcG\rightarrow\tcG'$ is the grafting map defined as follows: for any finite set $I$ and any 
    \begin{equation*}
        G=(V(G)=I, E(G)=\{e_1,\cdots,e_N\})\in\wE\circ\tcG[I]
    \end{equation*}
    we set
    \begin{equation*}
        \Bb(G)=(V(G)\cup\{*_I\},E(G)\cup_{i\in I}\{f_i \})\in\tcG'[I]
    \end{equation*}
    where $f_i$ is a new edge between the new vertex $*_I$ and the old vertex $i\in I=V(G)$. In words: the grafting map adds one vertex and attach this new vertex to every other vertices of the graph.

    Further set $\bB:\tcG'\longrightarrow\wE\circ\tcG$ to be the cografting map that on a graph in $\tcG'[I]$ removes the vertex $*_I$ and all the edges that have an end at this vertex. Then $(u_{\Bb},d_{\tiny\bB})$ is a pair of balanced up and down operators on $\wE\circ\tcG$ with $\lambda_n=1$ for any $n\in\N$.
\end{example}
\begin{proof}
    By proposition \ref{prop:balanced_grafting_cografting}, it is enough to prove that the equations \eqref{eq:balanced_1}-\eqref{eq:balanced_4} hold.
 \begin{enumerate}
        \item Equation \eqref{eq:balanced_1} trivially holds.
        \item Let $n\in\N$ and let $G\in\tcG[n]$. Then it is clear that $\bB^2(G)=G_{1,2}$. Since removing the vertex 1 and then the vertex 2 or removing first the vertex 2 and then the vertex 1 lead to the same result, it is clear that we also have $\bB^2((1,2).G)=G_{1,2}$. Then Equation \eqref{eq:balanced_2} holds.
        \item Equation \eqref{eq:balanced_3} also trivially holds with $\lambda_n=1$ for all $n\in\N$.
        \item For Equation \eqref{eq:balanced_4}, we use the following graphical conventions rather than a lengthy technical proof.
    \begin{itemize}
        \item For a graph $G\in\tcG[n]$ and $I=\{i_1,\cdots,i_k\}\subseteq[n]$, we write $G_{i_1,\cdots,i_k}$ (eventually in a circle) for the graph obtained by removing all the elements of $I$ to the vertices of $G$, and all the edges attached at either end to these vertices. 
        \item We use the standard notation that a solid edge between two vertices represents an edge of the graph.
        \item A dotted edge between two vertices indicates that there is possibly an edge between these two vertices.
        \item A double edge between a vertex and a subgraph indicates that this vertex is possibly related to some vertices of the subgraph. The number of edges between this vertex and the subgraph can be zero, and can also be the number of vertices of the subgraph.
        \item A triple edge between a vertex and a subgraph indicates that this vertices is related to each of the vertices of the subgraph.
    \end{itemize}
     Then for $k\in[n]$, $k\neq1$, write 
        
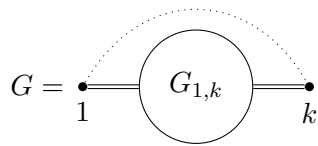
\begin{equation*}
            \begin{tikzpicture}[every node/.style={font=\small}]

\node[label={$G=$}] (G) at (-2.1,-0.4) {};
  \node[draw, circle, minimum size=1.5cm] (A) at (0,0) {$G_{1,k}$};

  \node[circle, fill, inner sep=1.2pt, label=below:{$1$}] (L) at (-1.5,0) {};
  \node[circle, fill, inner sep=1.2pt, label=below:{$k$}] (R) at (1.5,0) {};

  \draw[line width=0.4pt, double, double distance=1pt] (L) -- (A);
  \draw[line width=0.4pt, double, double distance=1pt] (A) -- (R);

  \draw[dotted]
    (L) to[out=60, in=120, looseness=1.3] (R);

\end{tikzpicture}
        \end{equation*}
        . Then we have 
        
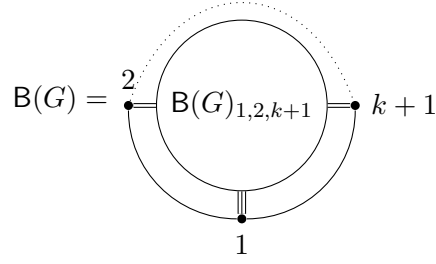
\begin{equation*}
        \begin{tikzpicture}[every node/.style={font=\small}]

\node[label={$\Bb(G)=$}] (G) at (-2.4,-0.4) {};

  \node[draw, circle, minimum size=1.5cm] (A) at (0,0) {$\Bb(G)_{1,2,k+1}$};

  \node[circle, fill, inner sep=1.2pt, label=above:{$2$}] (L) at (-1.5,0) {};
  \node[circle, fill, inner sep=1.2pt, label=right:{$k+1$}] (R) at (1.5,0) {};
  \node[circle, fill, inner sep=1.2pt, label=below:{$1$}] (D) at (0,-1.5) {};

  \draw[line width=0.4pt, double, double distance=1pt] (L) -- (A);
  \draw[line width=0.4pt, double, double distance=1pt] (A) -- (R);
  
  \draw[line width=0.4pt, double, double distance=2.4pt] (A) -- (D);
  \draw[line width=0.4pt] (A) -- (D);
  
  \draw[line width=0.4pt] (D) to[out=0, in=270] (R);
  \draw[line width=0.4pt] (L) to[out=270, in=180] (D);

  \draw[dotted]
    (L) to[out=70, in=110, looseness=1.6] (R);

\end{tikzpicture}
  \end{equation*}
  . Applying the permutation $(1,k+1)$ and then $\bB$ we obtain the LHS of Equation \eqref{eq:balanced_4} to be
  
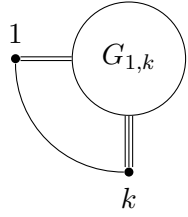
\begin{equation*}
      \begin{tikzpicture}[every node/.style={font=\small}]


  \node[draw, circle, minimum size=1.5cm] (A) at (0,0) {$G_{1,k}$};

  \node[circle, fill, inner sep=1.2pt, label=above:{$1$}] (L) at (-1.5,0) {};
  \node[circle, fill, inner sep=1.2pt, label=below:{$k$}] (D) at (0,-1.5) {};

  \draw[line width=0.4pt, double, double distance=1pt] (L) -- (A);
  
  \draw[line width=0.4pt, double, double distance=2.4pt] (A) -- (D);
  \draw[line width=0.4pt] (A) -- (D);
  
  \draw[line width=0.4pt] (L) to[out=270, in=180] (D);


\end{tikzpicture}
  \end{equation*}
  . For the RHS, starting from the same representation of $G$ and applying the permutation $(1,k)$ and then $\bB$ we obtain:
  
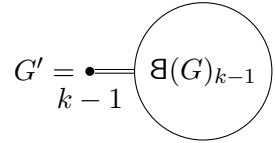
\begin{equation*}
      \begin{tikzpicture}[every node/.style={font=\small}]

\node[label={$G'=$}] (G) at (-2.1,-0.4) {};

  \node[draw, circle, minimum size=1.5cm] (A) at (0,0) {$\bB(G)_{k-1}$};

  \node[circle, fill, inner sep=1.2pt, label=below:{$k-1$}] (L) at (-1.5,0) {};

  \draw[line width=0.4pt, double, double distance=1pt] (L) -- (A);


\end{tikzpicture}
  \end{equation*}
  for $G':=\bB((1,k).G)$. We then have
  
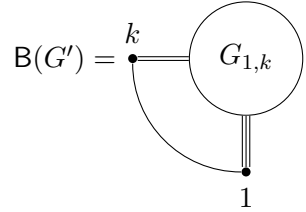
\begin{equation*}
      \begin{tikzpicture}[every node/.style={font=\small}]

\node[label={$\Bb(G')=$}] (G) at (-2.4,-0.5) {};

  \node[draw, circle, minimum size=1.5cm] (A) at (0,0) {$G_{1,k}$};

  \node[circle, fill, inner sep=1.2pt, label=above:{$k$}] (L) at (-1.5,0) {};
  \node[circle, fill, inner sep=1.2pt, label=below:{$1$}] (D) at (0,-1.5) {};

  \draw[line width=0.4pt, double, double distance=1pt] (L) -- (A);

  \draw[line width=0.4pt, double, double distance=2.4pt] (A) -- (D);
  \draw[line width=0.4pt] (A) -- (D);
  
  \draw[line width=0.4pt] (L) to[out=270, in=180] (D);


\end{tikzpicture}
  \end{equation*}
  . Exchanging the vertices $1$ and $k$ with the permutation $(1,k)$, we see that the RHS coincides with the LHS. Thus Equation \eqref{eq:balanced_4} holds and we have an example of a pair of balanced operators as claimed.
    \end{enumerate}
\end{proof}

\subsubsection{An example on rooted trees}

Our last example is on the species of rooted trees.
\begin{example}\label{ex:balanced_rooted_trees}
Let $\taa$ be the species of (non-planar) rooted trees. Gu\c{t}\v{a} and Maasen defined in \cite{gucta2002symmetric} (see Section 4.1) a simple pair $(u,d)$ of up and down operators on $\taa$, which is not balanced and therefore does not leads to creation/annihilation operators with commutation relations (for a detailed description of why their failed to be balanced, see \cite{AM2010}, Example 19.32). We define here an up operator $u:\taa\to\taa'$ and a down operator $d:\taa'\to\taa$ who does satisfy the balanced property. 

\medbreak

For a finite nonempty set $I$ and a rooted tree $t\in\taa[I]$, let $u_I(t)$ be the rooted tree on $I^+=I\sqcup\{*_I\}$ obtained
by adjoining the \emph{ghost vertex} $*_I$ as a new leaf attached to the root of $t$.

\medbreak
Conversely, let $t\in\taa'[I]=\taa[I^+]$. If the ghost vertex $*_I$ is the root of $t$, set $d_I(t)=0$. Otherwise, delete $*_I$ and all its incident edges, and graft each of its former children directly to the root of $t$. All other edges and the root are left unchanged. Extending both constructions linearly gives morphisms of species $u$ and $d$.

\medbreak

For instance, the actions of $u_I$ and $d_I$ are illustrated side by side below. The enlarged symbol $\star$ represents the
ghost vertex $*_I$. The two illustrations have independent input trees, each denoted by $t$.
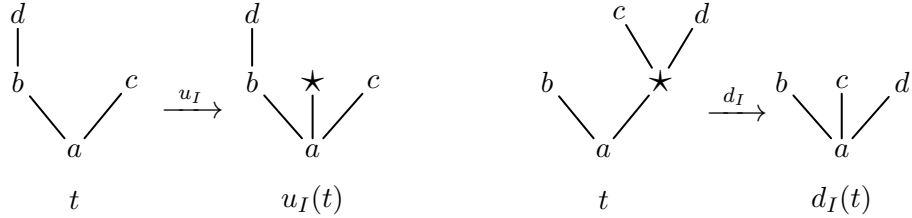
\begin{center}
\begin{tikzpicture}[
        x=1cm, y=0.91cm,
        vertex/.style={font=\small, inner sep=2pt},
        ghost/.style={font=\Large, inner sep=1pt},
        treeedge/.style={draw=black, line width=0.75pt,
                         line cap=round},
        operation/.style={font=\normalsize, inner sep=1pt}
    ]

       \node[vertex] (ua) at (0,0) {$a$};
       \node[vertex] (ub) at (-0.75,1) {$b$};
       \node[vertex] (uc) at (0.75,1) {$c$};
       \node[vertex] (ud) at (-0.75,2) {$d$};
       \draw[treeedge] (ua)--(ub) (ua)--(uc) (ub)--(ud);

       \node[vertex] (uva) at (3.15,0) {$a$};
       \node[vertex] (uvb) at (2.35,1) {$b$};
       \node[ghost]  (uvs) at (3.15,1) {$\star$};
       \node[vertex] (uvc) at (3.95,1) {$c$};
       \node[vertex] (uvd) at (2.35,2) {$d$};
       \draw[treeedge] (uva)--(uvb) (uva)--(uvs)
                       (uva)--(uvc) (uvb)--(uvd);
       \node[operation] at (1.56,0.70) {$\xrightarrow{\ u_I\ }$};
       \node[vertex] at (0,-0.72) {$t$};
       \node[vertex] at (3.15,-0.72) {$u_I(t)$};

       \node[vertex] (da) at (7.00,0) {$a$};
       \node[vertex] (db) at (6.25,1) {$b$};
       \node[ghost]  (ds) at (7.75,1) {$\star$};
       \node[vertex] (dc) at (7.20,2) {$c$};
       \node[vertex] (dd) at (8.30,2) {$d$};
       \draw[treeedge] (da)--(db) (da)--(ds)
                       (ds)--(dc) (ds)--(dd);

       \node[vertex] (dva) at (10.15,0) {$a$};
       \node[vertex] (dvb) at (9.35,1) {$b$};
       \node[vertex] (dvc) at (10.15,1) {$c$};
       \node[vertex] (dvd) at (10.95,1) {$d$};
       \draw[treeedge] (dva)--(dvb) (dva)--(dvc)
                       (dva)--(dvd);
       \node[operation] at (8.76,0.70) {$\xrightarrow{\ d_I\ }$};
       \node[vertex] at (7.00,-0.72) {$t$};
       \node[vertex] at (10.15,-0.72) {$d_I(t)$};
\end{tikzpicture}
\end{center}

In the right-hand illustration, the children $c$ and $d$ of the ghost vertex become children of the root $a$ after the ghost is deleted. Then $(u,d)$ is a pair of balanced up and down operators on $\taa$, with $\lambda_n=1$ for every $n\geq1$ (and with the degree-zero identity understood on $\taa[\emptyset]=0$).
\end{example}

\begin{proof}
We verify equations \eqref{eq:up_down1}--\eqref{eq:up_down4}.
Throughout, if a vertex to be removed is the root, the result is understood to be zero.
\begin{enumerate}[label=(\roman*)]
    \item Applying $u$ twice adds two new leaves, both adjacent
    to the original root. Exchanging their labels does not
    change the resulting rooted tree. Hence
    $(1,2)\cdot u^2(t)=u^2(t)$.

    \item Removing two distinct vertices in either order gives the same result: if either is the root, both compositions
    vanish; otherwise, each surviving child whose parent is one of the removed vertices is reattached to the root, while all other parent--child relations are unchanged. Thus $d^2((1,2)\cdot t)=d^2(t)$.

    \item The ghost vertex in $u_I(t)$ is a leaf. Removing it therefore recovers $t$, so $d_Iu_I(t)=t$. Consequently
    $d\circ u=\mathrm{id}_{\taa}$ and $\lambda_n=1$.

    \item Finally, consider the mixed relation \eqref{eq:up_down4}. We use the convention implicit in
    that relation: on $[n]$, the up operator inserts the new vertex with label $1$, shifting the old labels to $2,\ldots,n+1$, whereas the down operator deletes label $1$ and shifts the remaining labels down by one. Fix $1\leq k<n$ and $t\in\taa[n]$. On the left-hand side,
    \[
        d\bigl((k+1,1)\cdot u(t)\bigr),
    \]
    the transposition exchanges the newly created leaf with the vertex of $t$ originally labeled $k$. The down operator then removes that original vertex and grafts each of its children to the root. Indeed, the newly created leaf remains a child of the root and carries label $k$ after relabeling.
     On the right-hand side,
    \[(k,1)\cdot u\bigl(d((k,1)\cdot t)\bigr),\]
    the same original vertex is removed first, its children are grafted to the root, and a new leaf is subsequently attached to the root and relabeled $k$. The resulting rooted trees coincide. If the vertex originally labeled $k$ is the root, both sides are zero. This proves
    \eqref{eq:up_down4}.
    \end{enumerate}
    The four relations show that $(u,d)$ is balanced.
\end{proof}
Notice that the operators $u$ and $d$ of example ~\ref{ex:balanced_rooted_trees} are from grafting and cografting maps in the sense of proposition~\ref{prop:balanced_grafting_cografting}. Set $\bB$ to be like $d$, and $\Bb$ to be zero on partitions of length at least 2 and defined as $u$ in partitions of length one.

\medbreak

\section{Universal properties of the species of rooted trees} \label{sec:trees}

\subsection{Universality of the branching of rooted trees}

We present in this subsection an analogue of the universality result of the celebrated Connes-Kreimer Hopf algebra with its grafting operator $(\mathcal{H}_{CK}, B_+)$ introduced in \cite{connes1999hopf}. These properties use the branching and cobranching operators defined in \S \ref{subsubsec:grafting_trees} and \ref{ex:cografting_trees}, and is formulated in a category of monoids (resp. comonoids) with up (resp. down) operators. These operators are the one obtained from the branching and cobranching maps as explained in definition \ref{def:up_dow_ensembliste}.

\begin{theorem}\label{Universality}
\begin{enumerate}
    \item[i.] 
Let $\mathbf{a}$ be the species of trees and $\Bb: \wE\circ\taa\to \mathbf{a}'$ its grafting map defined in \S \ref{subsubsec:grafting_trees}. For every commutative monoid $(\tq, \nu, \kappa)$ with up operator $u: \tq \to \tq'$, there exists a unique map of monoids with up operators $\phi: (\wE\circ\taa, \mu, \iota, u_{\Bb}) \to (\tq, \nu, \kappa, u)$. In particular, the following diagram commutes:
    \begin{equation}\label{univdiagram}
    \xymatrix{
    \wE\circ\taa \ar[r]^-{\phi} \ar[d]_{u_{\Bb}} & \tq \ar[d]^u\\
    \wE\circ\taa' \ar[r]^-{\phi'}  & \tq' 
    }
\end{equation}
\item[ii.] Let $\bB:\taa'\longrightarrow\wE\circ\taa$ be the cografting map on trees defined in in \S \ref{subsubsec:grafting_trees} and $(\tp,\Delta,d)$ be a finite dimensional comonoid (i.e. such that dim$(\tp[I])<+\infty$ for any finite set $I$) with down operator. Then it exists a unique morphism of comonoids with down operators $\psi^{\rm co}:(\tp,\Delta,d)\longrightarrow(\wE\circ\taa,\Delta^*,d_{\tiny\bB})$, where $\Delta^*$ is the free coproduct on $\wE\circ\taa$.
\end{enumerate}
\end{theorem}

\begin{proof}
\begin{enumerate}
    \item[i.]
We prove first the existence of the map $\phi: \wE\circ\taa\to \tq$. Let $1_{\mathbf{F}}:= \iota_\emptyset(1_{\mathbb{K}}) \in \wE\circ\taa[\emptyset]$ and $1_{\tq}:= \kappa_\emptyset(1_{\mathbb{K}}) \in \tq[\emptyset]$ be the unit elements of $\wE\circ\taa$ and $\tq$, respectively. We put
\begin{equation}\label{phi0}
    \phi_{\emptyset}(1_{\mathbf{F}}):=1_{\tq}.
\end{equation}
Let $I$ be a non-empty finite set. Recall that the vector space $\wE\circ\taa[I]$ is given by
\[\wE\circ\taa[I] := \bigoplus_{\pi \, \vdash \, I} \taa(\pi),\]
where, for every partition $\pi$ of $I$, th space $\taa(\pi)$ is defined as
\[\taa(\pi):=\bigotimes_{B \in \pi}\taa[B].\]
We treat two cases separately.
\begin{itemize}
    \item Let $\pi=\{B_1, \hdots, B_k\} \vdash I$ be a partition with more than one block. For every $1 \leq i \leq k$, let $t_i \in \taa[B_i]$ be a rooted treed labelled with the set $B_i$. If 
    \begin{equation}\label{F}
    F=t_1  \cdots t_k \in \taa(\pi),
    \end{equation}
    where, as before, we do not indicate the symmetric tensor products between the trees, we let
    \begin{equation}\label{phi1}
        \phi_I(F):= \phi_{B_1}(t_1)\cdot\,  \cdots \, \cdot \phi_{B_k}(t_k) =\nu_{B_1,\cdots,B_k}(\phi_{B_1}(t_1)\otimes  \cdots \otimes \phi_{B_k})\in \tq[I].
    \end{equation}
    Notice that $\nu_{B_1,\cdots,B_k}$, the iteration of the product on $\tq$, is defined and its image is contained in $\tq[I]$ since $\nu$ endows $\tq$ with a monoid structure.

\

    \item Suppose that $\pi = \{I\}$, the (unique) partition of $I$ with one block. 

    For every $i \in I$, let $I^{-i}:=I\setminus \{i\}$. Consider the bijection
    \[\text{id}_i: I^{-i} \sqcup \{*_{I^{-i}}\} \to I,\]
    such that $*_{I^{-i}} \mapsto i$ and $j \mapsto j$ for every $j \in I^{-i}$. Then we have an isomorphism of vector spaces $q_i: \tq'[I^{-i}] \xrightarrow{\cong} \tq[I]$ given by
    \[q_i: \tq'[I^{-i}] = \tq[I\setminus \{i\} \sqcup \{*_{I^{-i}}\}] \xrightarrow{\tq[\text{id}_i]} \tq[I].\]

    Let $t \in \taa(\pi)=\taa[I]$ be a rooted tree. If $i \in I$ is the label of the root of $t$, then there exist a decomposition $I^{-i}=S_1 \sqcup \cdots \sqcup S_n$ and rooted trees $t_j \in \taa[S_j]$, for $1 \leq j \leq n$, such that $t =B_+^i( t_1 \hdots, t_n)$, with $B_+$ the usual combinatorial grafting operator of forests. We define
    \begin{equation}\label{phi2}
        \phi_I(t):=q_i \circ u_{I^{-i}}\left(\phi_{S_1}(t_1)\cdot \, \cdots \, \cdot \phi_{S_n}(t_n) \right).
    \end{equation}    
\end{itemize}

\

By \eqref{phi1}, $\phi$ is a monoid map. If $F=t_1 \cdots  t_k$ as in \eqref{F}, we have
\begin{align*}
    \phi'_I \circ u_{\Bb}(F)=\phi'_I(\Bb_I(t_1 \hdots t_k)) &= \phi_{I^+}(\Bb_I(t_1 \hdots t_k))\\
    &=q_{*_I}\circ u_{I^+\setminus \{*_I\}}\left(\phi_{B_1}(t_1)\cdot\,  \cdots \, \cdot \phi_{B_k}(t_k)\right)\\
    &= u_I\left(\phi_{B_1}(t_1)\cdot\,  \cdots \, \cdot \phi_{B_k}(t_k)\right) \tag{\footnotesize since $q_{*_I} = \tq[\text{id}_{*_I}] = \tq[\text{id}_{I^+}]=\text{id}_{\tq[I^+]}=\text{id}_{\tq'[I]}$ }\\
    &=u_I \circ \phi_I(F).
\end{align*}

Therefore, the diagram \eqref{univdiagram} commutes.

\

We still need to verify the uniqueness of $\phi = (\phi_I)_I$. We proceed by induction on the size of the finite sets $I$. Let $\psi:\wE\circ\taa\to \tq$ be a monoid map such that \eqref{univdiagram} commutes. Since $\psi$ is a monoid map, we must have
\[\psi_{\emptyset}(1_{\mathbf{F}})=1_{\tq} =\phi_{\emptyset}(1_{\mathbf{F}}).\]
For the induction step, we will once again separate two cases.
\begin{itemize}
    \item Let $I$ be a non-empty finite set. Let $\pi =\{I\}$ be the partition of $I$ with only one block. Let $t \in \taa(\pi)=\taa[I]$ be a rooted tree. If $i \in I$ is the label of the root of $t$, then there exist a decomposition $I^{-i}=S_1 \sqcup \cdots \sqcup S_n$ and rooted trees $t_j \in \taa[S_j]$, for $1 \leq j \leq n$, such that $t = B_+^i(t_1 \cdots t_n)$.

\

Consider again the bijection
\[\text{id}_i: I^{-i} \sqcup \{*_{I^{-i}}\} \to I,\]
such that $*_{I^{-i}} \mapsto i$ and $j \mapsto j$ for every $j \in I^{-i}$. We have an isomorphism of vector spaces $a_i: \wE\circ\taa'[I^{-i}] \xrightarrow{\cong} \wE\circ\taa[I]$ given by
\[a_i:  \wE\circ\taa'[I^{-i}] = \wE\circ\taa[I\setminus \{i\} \sqcup \{*_{I^{-i}}\}] \xrightarrow{\wE\circ\taa[\text{id}_i]} \wE\circ\taa[I].\]

\

Let $s:= a_i^{-1}(t)$. Under the action of $a_i^{-1}$, all labels of the tree $t$ remain fixed, except the label of the root, which is changed from $i$ to $*_{I^{-i}}$. Therefore, 
\[s = B_+^{*_{I^{-i}}}(t_1\cdots t_n) = u_{\Bb}(t_1 \cdots  t_n),\]
where $t_1  \cdots t_n \in \wE\circ\taa[I^{-i}]$. Hence, we have
\begin{align*}
    \psi'_{I^{-i}}(s)&= \psi'_{I^{-i}}\circ u_B(t_1  \cdots t_n)\\
    &=u_{I^{-i}}\circ \psi_{I^{-i}}(t_1  \cdots  t_n) \tag{\footnotesize $\psi$ satisfies \eqref{univdiagram}}\\
    &= u_{I^{-i}}(\psi_{S_1}(t_1)\cdot \, \cdots \, \cdot \psi_{S_n}(t_n)) \tag{\footnotesize $\psi$ is a monoid morphism}\\
    &= u_{I^{-i}}(\phi_{S_1}(t_1)\cdot \, \cdots \, \cdot \phi_{S_n}(t_n)) \tag{\footnotesize by the induction hypothesis}\\
    &=u_{I^{-i}}\circ \phi_{I^{-i}}(t_1  \cdots  t_n) \tag{\footnotesize $\phi$ is a monoid morphism}\\
    &=\phi'_{I^{-i}}\circ u_B(t_1  \cdots  t_n)\tag{\footnotesize $\phi$ satisfies \eqref{univdiagram}}\\
    &=\phi'_{I^{-i}}(s).
\end{align*}

\

Since $a_i^{-1}= \wE\circ\taa[\text{id}_i]^{-1}=\wE\circ\taa[\text{id}_i^{-1}]$ (by functoriality of $\wE\circ\taa$), then $s= \wE\circ\taa[\text{id}_i^{-1}](t)$, so
\[\psi'_{I^{-i}}\circ (\wE\circ\taa)[\text{id}_i^{-1}](t) = \phi'_{I^{-i}}\circ (\wE\circ\taa)[\text{id}_i^{-1}](t) \]

\

Finally, by naturality of the species map $\phi$ and $\psi$, we obtain
\[\tq[\text{id}_i^{-1}]\circ \psi_I (t) = \tq[\text{id}_i^{-1}]\circ \phi_I(t).\]

As $\tq[\text{id}_i^{-1}]$ is an isomorphism, we must have $\psi_I (t) = \phi_I(t)$.

\item For the second case, let $\pi\vdash I$ be a partition of $I$ with more than one block and consider a forest $F=t_1\cdots t_k \in \taa(\pi)$ as in \eqref{F}, with $k >1$. Applying the recurrence hypothesis to each element $t_1, \hdots, t_k$, we obtain
\begin{align*}
\psi_I(F)= \psi_I(t_1 \cdots t_k ) &= \psi_{B_1}(t_1)\cdot \, \cdots \, \cdot \psi_{B_k}(t_k)\\
&=\phi_{B_1}(t_1)\cdot \, \cdots \, \cdot \phi_{B_k}(t_k) = \phi_I(t_1 \cdots t_k ) = \phi_I(F).
\end{align*}
\end{itemize}
These are the two possible cases. So $\psi_I=\phi_I$, and by induction $\psi=\phi$.
\item[ii.]
     Let $(\tp,\Delta,d)$ be a comonoid with down operator such that dim$(\tp[I])<+\infty$ for any finite set $I$. Then $(\tp^*,\Delta^*,d^*)$ is a monoid with up operator. Thus by the first point it exists a unique morphism of monoids with up operator 
     $\phi:(\wE\circ\taa,\mu,u_\Bb)\longrightarrow(\tp^*,\Delta^*,d^*)$. Taking the dual we have a map
     \begin{equation*}
         \phi^*:(\tp^*)^*=\tp\longrightarrow(\wE\circ\taa)^*=\wE\circ\taa
     \end{equation*}
     (where we have used again that dim$(\tp[I])<+\infty$). Now since $\phi$ is a morphism of monoids with up operators, $\phi^*:(\tp,\Delta,d)\longrightarrow(\wE\circ\taa,\mu^*,u_\Bb^*)$ is a morphism of comonoids with down operators. The result then follows from the fact that the free coproduct $\Delta^*$ on $\wE\circ\taa$ is exactly the dual of $\mu$ and the fact that $u_\Bb^*=d_{\tiny \bB}$.

     To see this last point, take $x$ and $f$ two forests in $\wE\circ\taa[I]$. Then by definition of $u_\Bb^*$ we have 
     \begin{equation*}
         \langle u_\Bb^*(x),f\rangle  = \langle x,u_\Bb(f)\rangle  =
         \begin{cases}
             1 \quad\text{ if }x=B_+^{*_I}(f) \\
             0\quad\text{ otherwise}.
         \end{cases} 
     \end{equation*}
     This is exactly $\langle d_{\tiny \bB}(x),f\rangle$ so we have $u_\Bb^*=d_{\tiny \bB}$ as claimed and the theorem is proven. 
\end{enumerate}
\end{proof}
The second point of this theorem could be proved without the assumption of finite dimensionality, by essentially dualizing the proof of the first point. We omit this proof since it does not bring much to the discussion and the level of generality it brings it not of much use in practice.

\subsection{Applications of the universal property}



Theorem~\ref{Universality} is formulated in the category of species with up operators. However, the categorical constructions of \S\ref{subsec:up_down_cat} allow us to formulate it in the category $\Spgrdec$ of species with grafting and decomposable morphisms (see remark~\ref{rk:dec}).
\begin{cor} \label{coro:univ_grafting}
 $(\wE\circ\taa,\Bb)$ is the initial object of the category $\Spgrdec$.
\end{cor}
\begin{proof}
    Let $(\wE\circ\tp,\Cc)$ be a species with grafting. Applying the functor $\fH:\Spgrdec\longrightarrow\Moup$ of \S\ref{subsec:up_down_cat} we obtain a monoid with up operator $\fH(\tp,\Cc)$. Then by theorem~\ref{Universality} there is a unique morphism of monoids with up morphism $f:(\wE\circ\taa,\mu,u_\Bb)\longrightarrow\fH(\tp,\Cc)$, where $\mu$ is the concatenation of $\wE\circ\taa$.

    By the third point of theorem~\ref{thm:comparison_up_grafting_NEW} $\fH$ is left adjoint of $\fU\fF$ so there is a unique decomposable morphism of species with graftings $g:\fU\fF(\wE\circ\taa,\mu,u_\Bb)\longrightarrow(\tp,\Cc)$. By definition of the functors $\fU$ and $\fF$ we have $\fU\fF(\wE\circ\taa,\mu,u_\Bb)=(\wE\circ\taa,\Bb)$ which is then the initial object of $\Spgrdec$ as claimed.
\end{proof}

\medskip


Let $\epsilon:\wE\circ\taa{\longrightarrow\bf 1}_{\mathbb K}$ be the morphism of species defined by $\epsilon_\emptyset(1_\taa)=1$ and  $\epsilon_I(t)=0$ for any finite set $I$ with $|I|>0$. Here $1_\taa$ is the empty tree, which is the unit of $\wE\circ\taa$. Then, let us define the following map of species
\begin{align*}
    L:(\wE\circ\taa)\cdot(\wE\circ\taa) & \longrightarrow (\wE\circ\taa)'\cdot(\wE\circ\taa) + (\wE\circ\taa)\cdot(\wE\circ\taa)' \\
    x\otimes y& \longmapsto \Bb(x)\otimes\epsilon(y)1_\taa + x\otimes\Bb(y).
\end{align*}
We then have the following.
\begin{cor} \label{coro:CK_Hopf}
    There exists a map $\Delta:\wE\circ\taa\longrightarrow(\wE\circ\taa)\cdot(\wE\circ\taa)$ such that $(\wE\circ\taa,\mu,\Delta,\iota,\epsilon)$ is a commutative Hopf monoid. 

    Furthermore, let $(\tp,\nu,\Delta^\tp,\kappa,\varepsilon)$ be a commutative Hopf monoid  with up operator $u:\tp\longrightarrow\tp'$ such that, for any finite set $I$, any $S\sqcup T=I$, and any $x\in\tp[I]$ we have
    \begin{equation*}
        \Delta^\tp_{S,T+}\circ u_I(x) = (\text{id}_S\otimes u_T)\circ\Delta_{S,T}(x),\quad \Delta^\tp_{I+,\emptyset}\circ u_I(x) = u_I(x)\otimes 1_\tp.
    \end{equation*}
    Then there is a unique morphism of Hopf monoids $\phi:\wE\circ\taa\longrightarrow\tp$ which is a morphism of species with up operator.
\end{cor}
\begin{proof}
    Since $(\wE\circ\taa)'\cdot(\wE\circ\taa) + (\wE\circ\taa)\cdot(\wE\circ\taa)'=\left((\wE\circ\taa)\cdot(\wE\circ\taa)\right)'$, the map $L$ is an up operator on $(\wE\circ\taa)\cdot(\wE\circ\taa)$. Thus by theorem \ref{Universality} we have a morphism of monoids with up operators $\Delta:\wE\circ\taa\longrightarrow(\wE\circ\taa)\cdot(\wE\circ\taa)$.

    The coassociativity of $\Delta$ can be proven by adapting the standard proof (for example in \cite{FoissyBerlin}) to the language of species, and we omit it here. Then we have a connected bimonoid, which is a Hopf monoid by \cite[Proposition 6]{AM2013}.
    
    For the second point, theorem \ref{Universality} gives the existence and unicity of a morphism of monoids with up operators $\phi:\wE\circ\taa\longrightarrow\tp$. On can once again prove by adapting the standard proof (which can also be found in \cite{FoissyBerlin}, among others) to the language of species that this map is also a morphism of comonoids as needed.
\end{proof}

\medbreak

\section{Recovering the classical Hopf algebras of rooted forests}
\label{sec:recovering_CK_Foissy}

The universal properties of the preceding section can be transported from
species to graded vector spaces by means of Fock functors. We explain how
this recovers three familiar, but different, Hopf algebras: the
commutative Connes--Kreimer Hopf algebra of rooted forests, its
noncommutative counterpart on planar rooted forests, and Foissy's
infinitesimal Hopf algebra of planar rooted forests. All three share the
same grafting operator and a common Hochschild-type cocycle identity.
The distinction between the two planar versions lies in the
the compatibility of the coproduct with concatenation. We also
construct an up operator from ordered grafting and unshuffle, and
explain how the $q$-deformed Fock functor interpolates between the
ordinary and infinitesimal planar cut coproducts.

\subsection{The bosonic Fock functor and the branching operator}
\label{subsec:Fock_and_branching}

For a vector species $\tp$, its \emph{bosonic Fock space} is the graded
vector space
\begin{equation}\label{eq:Fock_bosonic_CK}
  \mathcal K(\tp):=\bigoplus_{n\geq0}\tp[n]_{\mathfrak S_n},
  \qquad
  \tp[n]_{\mathfrak S_n}:=
  \tp[n]/\langle\sigma x-x:\sigma\in\mathfrak S_n,\ x\in\tp[n]\rangle.
\end{equation}
Here $[n]=\{1,\ldots,n\}$, and a class in the coinvariant space is
denoted by $[x]$. The Fock functor is strong monoidal for the Cauchy
product and takes Hopf monoids in species to graded Hopf algebras; see
\cite[Chapter~15]{AM2010}. More concretely, the product of two classes
is obtained by putting their label sets in disjoint blocks and applying
the species product. If $\tp$ has a comonoid structure with coproduct $\Delta$, the coproduct of a class $[x]$, with
$x\in\tp[n]$, is induced by all decompositions $[n]=S\sqcup T$:
\begin{equation}\label{eq:Fock_coproduct_CK}
  \Delta_{\mathcal K}([x])=
  \sum_{S\sqcup T=[n]}
  [x_{(1),S}]\otimes[x_{(2),T}],
  \qquad
  \Delta_{S,T}(x)=\sum_{S\sqcup T=[n]} x_{(1),S}\otimes x_{(2),T}.
\end{equation}
In each summand the labels of $S$ and $T$ are standardized. In
particular, distinct labelled decompositions may have identical images
and therefore contribute multiplicities.

\medbreak

Suppose that a species $\tp$ has an up operator
$u:\tp\to\tp'$. Relabelling the additional ghost vertex as a new
ordinary vertex, and then passing to coinvariants, gives an operator
of degree $+1$,
\begin{equation}\label{eq:Fock_up_Bplus}
  B_u:\mathcal K(\tp)_n\longrightarrow
       \mathcal K(\tp)_{n+1},\qquad
  [x]\longmapsto[u_{[n]}(x)].
\end{equation}
In this formula the displayed class in degree $n+1$ is taken modulo
$\mathfrak S_{n+1}$, rather than merely modulo the permutations of
the original labels. This distinction is needed because
$\mathcal K(\tp')_n$ is not the $(n+1)$-degree part of
$\mathcal K(\tp)$. For the branching operator $u_{\Bb}$ on rooted
forests, the operator $B_{u_{\Bb}}$ is the usual $B_+$: it adjoins a
new root below all components of the forest.

\subsection{The commutative Connes--Kreimer Hopf algebra}
\label{subsec:Fock_classical_CK}

Recall that $\mathbf F=\wE\circ\taa$ is the species of non-planar
rooted forests, with disjoint union as its commutative product, and
that Corollary~\ref{coro:CK_Hopf} equips it with a Hopf monoid
coproduct $\Delta_{\rm CK}$. Its bosonic Fock image is
\begin{equation}\label{eq:recover_commutative_CK}
  \mathcal H_{\rm CK}:=\mathcal K(\mathbf F)
  =\bigoplus_{n\geq0}\mathbf F[n]_{\mathfrak S_n}.
\end{equation}
The degree-$n$ basis is indexed by unlabelled non-planar rooted
forests with $n$ vertices. Multiplication is disjoint union, so
$\mathcal H_{\rm CK}$ is the free commutative algebra on non-planar
rooted trees. Its unit is the empty forest $1$.

\medbreak

To describe the coproduct, an \Def{admissible cut} of a forest $F$ is
a collection of edges meeting each root-to-leaf path at most once.
Write $P^c(F)$ for the forest pruned above the cut and $R^c(F)$
for the forest remaining below it. We also include the empty and
full cuts, which contribute $1\otimes F$ and $F\otimes1$.
The coproduct furnished by the Fock functor is
\begin{equation}\label{eq:recover_CK_cuts}
  \Delta_{\rm CK}(F)=
    \sum_{c\in\operatorname{Adm}(F)}P^c(F)\otimes R^c(F).
\end{equation}
Here and below, a cut of a forest is an independent choice of cuts
in its component trees, and terms with equal unlabelled images are
counted with their multiplicities.

\begin{prop}\label{prop:recover_CK}
The image through the Fock functor $\mathcal{K}$ of the Hopf monoid of
Corollary~\ref{coro:CK_Hopf} is the classical Connes--Kreimer Hopf
algebra. Under the identification~\eqref{eq:recover_commutative_CK},
the up operator $u_{\Bb}$ becomes $B_+$ and
\begin{equation}\label{eq:CK_Hochschild}
  \Delta_{\rm CK}(B_+(F))
   =B_+(F)\otimes1+
     (\mathrm{id}\otimes B_+)\Delta_{\rm CK}(F).
\end{equation}
\end{prop}
\begin{proof}
The coinvariants identify precisely the labelled forests that differ
only by a relabelling. The species product becomes disjoint union.
At the labelled level, the coproduct of
Corollary~\ref{coro:CK_Hopf} splits a rooted forest into a pruned
upper part and a complementary rooted lower part. Such a split is
uniquely specified by an admissible cut; applying
\eqref{eq:Fock_coproduct_CK} sums over all label-set choices and
gives~\eqref{eq:recover_CK_cuts}. Finally, the defining compatibility
of $u_{\Bb}$ with the coproduct becomes
\eqref{eq:CK_Hochschild}. Since the species Hopf monoid is connected,
its bosonic Fock image is a graded connected Hopf algebra.
\end{proof}

\medbreak

\subsection{Planar forests and the noncommutative Connes--Kreimer Hopf algebra}
\label{subsec:Fock_planar_CK}

Let $\taa_{\rm pl}$ denote the positive species of \Def{planar
rooted trees}; at each vertex, the children are linearly ordered.
The species of \Def{planar rooted forests} is
\begin{equation}\label{eq:planar_forests_species}
  \mathbf F_{\rm pl}:=\wL\circ\taa_{\rm pl}.
\end{equation}
Thus a basis element is an \emph{ordered} list
$F=t_1|\cdots|t_r$ of planar rooted trees, and the product is
concatenation $FG$. It is the (noncommutative)
free monoid in species on $\taa_{\rm pl}$.

\medbreak

The ordered grafting map
\begin{equation*}
  \overrightarrow{\Bb}:\mathbf F_{\rm pl}
       \longrightarrow\taa_{\rm pl}'
\end{equation*}
adjoins a ghost-labelled root below the trees $t_1,\ldots,t_r$,
retaining their order as its list of children. Together with the
canonical inclusion $\taa_{\rm pl}'\hookrightarrow
\mathbf F_{\rm pl}'$, this gives an up operator
$v:\mathbf F_{\rm pl}\to\mathbf F_{\rm pl}'$. Its Fock image is the
usual grafting $B_+$ on planar rooted forests.

\medbreak

The next result is the planar version of the first point of theorem~\ref{Universality}. It is proved in the same way and we only outline the proof.
\begin{prop}\label{prop:planar_universal_up}
For every (not necessarily commutative) monoid $(\tp,\mu,\iota)$
with an up operator $u:\tp\to\tp'$, there is a unique morphism
of monoids with up operators
\begin{equation*}
  \phi:(\mathbf F_{\rm pl},v)\longrightarrow(\tp,u).
\end{equation*}
\end{prop}
\begin{proof}
Every planar rooted forest admits a unique factorization as an
ordered concatenation of planar rooted trees. Every planar rooted
tree with root labelled $i$ is uniquely of the form
$B_+^i(t_1|\cdots|t_r)$, for a uniquely determined ordered list of
smaller trees. Starting with $\phi(1)=1_{\tp}$, these two facts
recursively force $\phi(FG)=\phi(F)\phi(G)$ and
\begin{equation*}
  \phi(B_+^i(F))
   =\tp[\sigma_i]\bigl(u(\phi(F))\bigr),
\end{equation*}
where $\sigma_i$ identifies the newly adjoined ghost label with $i$
and fixes all other labels. Conversely, these rules define a natural
monoid morphism commuting with the up operators. Induction on the
number of vertices proves uniqueness.
\end{proof}

Equip $\mathbf F_{\rm pl}$ with the \Def{planar cut
coproduct} $\Delta_1$, defined as in
\eqref{eq:recover_CK_cuts}, but retaining the planar order on each
pruned and remaining forest. The coproduct is multiplicative with
respect to concatenation, with the usual multiplication on the
tensor product. Its naturality and coassociativity follow from the
same nested-cut argument as in the non-planar case. We therefore
obtain a connected Hopf monoid in species, and hence
\begin{equation}\label{eq:recover_NCK}
  \mathcal H_{\rm NCK}
       =\mathcal K(\mathbf F_{\rm pl},\Delta_1).
\end{equation}
This is the noncommutative Connes--Kreimer Hopf
algebra introduced by Foissy (see \cite{foissy2002algebres}): as an algebra it is freely generated by planar rooted trees,
and the coproduct is given by planar admissible cuts. In particular,
\begin{equation}\label{eq:NCK_Hochschild}
  \Delta_1(B_+(F))=B_+(F)\otimes1+
         (\mathrm{id}\otimes B_+)\Delta_1(F),
  \qquad \Delta_1(FG)=\Delta_1(F)\Delta_1(G).
\end{equation}
The second identity uses the ordinary tensor-product algebra
multiplication. Notice that it is the ordered, planar replacement
of $\wE\circ\taa$ by $\wL\circ\taa_{\rm pl}$ that makes the
resulting algebra noncommutative.

\subsection{Ordered grafting induced by the unshuffle coproduct}
\label{subsec:ordered_grafting_unshuffle}

We now construct an up operator from an ordered grafting by using the
\emph{unshuffle} (also called \emph{deshuffle}) coproduct on the free
monoid of a positive species. The construction uses arbitrary
subsequences of an ordered forest, rather than cuts at consecutive
positions. We shall distinguish the resulting up operator from the
ordinary branching operator $B_+$ of the Connes--Kreimer algebras.

\medbreak

Let $\tq$ be a positive species and put
$\mathcal T(\tq):=\wL\circ\tq$.
For an ordered word $w=x_1|\cdots|x_r$ of $\tq$-structures and a
subset $A=\{i_1<\cdots<i_s\}\subseteq[r]$, set
$w_A:=x_{i_1}|\cdots|x_{i_s}$, with $w_\emptyset=1$.
The \emph{unshuffle coproduct}
$\Delta_{\rm unsh}:\mathcal T(\tq)\to
\mathcal T(\tq)\cdot\mathcal T(\tq)$ is given, in the customary
notation for the sum of its components, by
\begin{equation}\label{eq:ordered_unshuffle_CK}
  \Delta_{\rm unsh}(w)
   =\sum_{A\subseteq[r]}
     w_A\otimes w_{[r]\setminus A}.
\end{equation}
At the species level, a component indexed by a decomposition of the
underlying label set vanishes unless each block carrying a
$\tq$-structure lies entirely in one part. Otherwise the component
retains the relative order of the blocks assigned to each part.
This coproduct is coassociative and cocommutative; together with
concatenation it makes $\mathcal T(\tq)$ a connected Hopf monoid.
In particular,
\begin{equation}\label{eq:ordered_unshuffle_mult}
 \Delta_{\rm unsh}(vw)
   =\Delta_{\rm unsh}(v)\Delta_{\rm unsh}(w),
\end{equation}
where multiplication on the right is the usual tensor-product
multiplication. This is the free Hopf monoid construction associated
with the trivial positive comonoid structure on $\tq$; see
\cite[Section~11.4]{AM2010}.

\medbreak

The twice-iterated coproduct
$\Delta_{\rm unsh}^{(2)}:=
(\Delta_{\rm unsh}\cdot\mathrm{id})\Delta_{\rm unsh}
=(\mathrm{id}\cdot\Delta_{\rm unsh})\Delta_{\rm unsh}$
assigns each block independently to one of three ordered tensor
factors. Consequently,
\begin{equation}\label{eq:ordered_triple_unshuffle_CK}
 \Delta_{\rm unsh}^{(2)}(w)
 =\sum_{A\sqcup B\sqcup C=[r]}
     w_A\otimes w_B\otimes w_C,
\end{equation}
where $A,B,C$ may be empty and the subwords $w_A,w_B,w_C$ need not
be consecutive. Thus an ordered word of length $r$ has $3^r$
terms in~\eqref{eq:ordered_triple_unshuffle_CK}, counted with
multiplicity.

\begin{prop}\label{prop:ordered_grafting_up_unshuffle}
Let
$\sBb:\wL\circ\tq\longrightarrow\tq'$
be an ordered grafting map. The composite
\begin{equation}\label{eq:ordered_grafting_threefold_map}
\begin{gathered}
  \widetilde u_{\sBb}:\quad
  \wL\circ\tq
   \xrightarrow{\ \Delta_{\rm unsh}^{(2)}\ }
       (\wL\circ\tq)^{\cdot3}
   \xrightarrow{\ \mathrm{id}\cdot\sBb\cdot\mathrm{id}\ }
       (\wL\circ\tq)\cdot\tq'\cdot(\wL\circ\tq)
  \\[4pt]
  \cong\ \tq'\cdot\bigl((\wL\cdot\wL)\circ\tq\bigr)
   \cong\ \tq'\cdot(\wL'\circ\tq)
   \cong\ (\wL\circ\tq)'
\end{gathered}
\end{equation}
defines a natural up operator on $\wL\circ\tq$.
If $\pi_{\tq'}:(\wL\circ\tq)'\to\tq'$ denotes the
projection onto the direct summand consisting of the single
$\tq'$-block containing the ghost, then
\begin{equation}\label{eq:ordered_unshuffle_projection}
   \pi_{\tq'}\circ\widetilde u_{\sBb}=\sBb.
\end{equation}
\end{prop}
\begin{proof}
The unshuffle coproduct, its iteration, and the ordered grafting
are natural in relabellings, so the first two arrows in
\eqref{eq:ordered_grafting_threefold_map} are morphisms of species.
The last isomorphisms are obtained from the composition and derivative rules
\begin{equation*}
(\tp\circ\tq)\cdot(\tp\circ\tq)\cong(\tp\cdot\tp)\circ\tq,
  \qquad
  \wL'\cong\wL\cdot\wL,
  \qquad
  (\wL\circ\tq)'
   \cong(\wL'\circ\tq)\cdot\tq'.
\end{equation*}
Indeed, if $I=S\sqcup T$, $\pi\vdash S$ and $\rho\vdash T$, then $\pi\sqcup\rho\vdash I$. Furthermore, distinguishing a ghost in a linear order separates the
blocks preceding and following the ghost-containing block.
The symmetry of the Cauchy product puts the $\tq'$-factor first
when desired.

\medbreak

For the projection identity, consider a summand indexed by
$A\sqcup B\sqcup C=[r]$ in
\eqref{eq:ordered_triple_unshuffle_CK}. Its image consists of
$w_A$, one distinguished $\tq'$-block $\sBb(w_B)$, and $w_C$.
It has no other blocks precisely when $A=C=\emptyset$, in which
case $B=[r]$. The surviving term is therefore $\sBb(w)$.
\end{proof}

Writing $[\sBb(v)]_*$ for the distinguished block containing
the ghost, the operator is explicitly
\begin{equation}\label{eq:ordered_grafting_unshuffle_formula}
 \widetilde u_{\sBb}(x_1|\cdots|x_r)
   =\sum_{A\sqcup B\sqcup C=[r]}
         w_A\,\big|\,[\sBb(w_B)]_*\,\big|\,w_C.
\end{equation}
The vertical bars on the right separate blocks of the resulting
pointed forest, with the ghost contained in the middle block.

\begin{example}\label{ex:ordered_grafting_two_blocks}
For a one-block word $x$ the three terms are
\begin{equation*}
 \widetilde u_{\sBb}(x)
 = [\sBb(1)]_*|x
   +[\sBb(x)]_*+x|[\sBb(1)]_*.
\end{equation*}
For a two-block word $x|y$ there are nine ordered assignments:
\begin{align*}
 \widetilde u_{\sBb}(x|y)
  ={}&[\sBb(1)]_*|x|y
       +[\sBb(x)]_*|y
       +[\sBb(y)]_*|x
       +[\sBb(x|y)]_*\\
    &+x|[\sBb(1)]_*|y
       +y|[\sBb(1)]_*|x\\
    &+x|[\sBb(y)]_*
       +y|[\sBb(x)]_*
       +x|y|[\sBb(1)]_*.
\end{align*}
Unlike an operation based on consecutive intervals, this sum
includes the terms $[\sBb(y)]_*|x$ and
$y|[\sBb(1)]_*|x$ arising from nonconsecutive assignments of the
original ordered blocks.
\end{example}

\begin{prop}\label{prop:unshuffle_grafting_and_NCK}
Take $\tq=\taa_{\rm pl}$ and let
$\sBb$ be the usual ordered root-grafting.
Let $\pi_{\rm tree}:\mathbf F_{\rm pl}'\to
\taa_{\rm pl}'$ be projection onto pointed forests consisting of
one rooted tree. Then
\begin{equation}\label{eq:unshuffle_grafting_connected_part}
 \pi_{\rm tree}\circ
 \widetilde u_{\overrightarrow\Bb}
 =\overrightarrow\Bb.
\end{equation}
After taking bosonic Fock functor, the right-hand side is
the usual Connes--Kreimer grafting operator $B_+$.
Moreover, the unshuffle coproduct on $\mathbf F_{\rm pl}$ is
exactly the sub-sum of the ordinary planar cut coproduct
$\Delta_1$ in which no edge internal to a component tree is cut.
\end{prop}
\begin{proof}
The first identity is the specialization of
\eqref{eq:ordered_unshuffle_projection}.
For the second statement write
$F=t_1|\cdots|t_r$. The ordinary cut coproduct is
multiplicative, and for each component $t_i$ it contains both
trivial terms $t_i\otimes1$ and $1\otimes t_i$.
Choosing only one of these two terms for each component gives
\begin{equation*}
 \prod_{i=1}^r(t_i\otimes1+1\otimes t_i)
   =\sum_{A\subseteq[r]}F_A\otimes F_{[r]\setminus A},
\end{equation*}
which is precisely $\Delta_{\rm unsh}(F)$. All other terms in
$\Delta_1(F)$ involve a nontrivial admissible cut in at least
one component tree.
\end{proof}

\medbreak

\begin{remark}\label{rem:deconcatenation_ordered_branching_outlook}
There is a second ordered branching construction that can be considered.  If
$\sBb:\wL\circ\tq\to\tq'$ is an ordered grafting, the
\emph{deconcatenation} coproduct on $\wL\circ\tq$ produces
\begin{equation*}
 (\mathrm{id}\cdot\sBb\cdot\mathrm{id})
 \circ\Delta_{\mathrm{dec}}^{(2)}:
 \wL\circ\tq\longrightarrow
 (\wL\circ\tq)\cdot\tq'\cdot(\wL\circ\tq)
 \cong(\wL\circ\tq)'.
\end{equation*}
For $w=x_1|\cdots|x_r$, this operator sums over
$0\leq i\leq j\leq r$, grafting the consecutive interval
$x_{i+1}|\cdots|x_j$ and leaving the prefix and suffix fixed.
The result is an up coderivation for \emph{deconcatenation}:
its coproduct splits the resulting word on either side of the
distinguished grafted block.  This property is different from the
unshuffle construction above.  We will develop this deconcatenation-based construction and
its further applications in a future work.
\end{remark}

\medbreak

\subsection{A $q$-infinitesimal interpolation via deformed Fock functors}
\label{subsec:qCK_deformed_Fock}

The ordinary and infinitesimal cut coproducts admit a common
$q$-deformation.   The construction can be read
directly from the universal noncommutative forest monoid and the
deformed Fock functor of Aguiar and Mahajan
\cite[Sections~13.3.5, 16.1, and 17.5.4]{AM2010}.

\medbreak

\paragraph{The deformed Fock construction.}
Recall the planar forest Hopf monoid
$\mathbf F_{\rm pl}=\wL\circ\taa_{\rm pl}$, with concatenation and
the ordinary cut coproduct $\Delta_1$ of
Subsection~\ref{subsec:Fock_planar_CK}.  We consider the
\Def{pruned forest in the first tensor factor}.  For $q\in\Bbbk$, consider the $q$-deformed full
Fock functor
\begin{equation}\label{eq:qCK_full_Fock}
 \mathcal K_q(\mathbf F_{\rm pl})
    =\bigoplus_{n\geq0}\mathbf F_{\rm pl}[n].
\end{equation}
Its product is concatenation followed by shifting the labels in
the second factor.  For a labelled forest $F$ on $[n]$, its coproduct
is
\begin{equation}\label{eq:qCK_labelled_Fock}
\begin{aligned}
 \Delta_q(F)
 &=\sum_{\substack{[n]=S\sqcup T\\S\text{ is pruned-admissible}}}
    q^{\operatorname{sch}_n(S)}
    \operatorname{std}(F|_S)\otimes\operatorname{std}(F|_T),\\
 \operatorname{sch}_n(S)
 &=\#\{(s,t)\in S\times T:s>t\},
\end{aligned}
\end{equation}
where $\operatorname{std}$ denotes the usual standardization map on the set of labels. The restriction to a pruned-admissible subset $S$ is the upper
forest $P^c(F)$ of the corresponding admissible cut, while the
complement is the remaining forest $R^c(F)$.
This follows by applying the colax transformation of
$\mathcal K_q$ to the cut coproduct of the Hopf monoid; see
\cite[Theorem~16.1]{AM2010}.  The resulting graded algebra and
coalgebra satisfy the $q$-braided bialgebra compatibility.

\medbreak

To recover the usual \emph{unlabelled} planar forest basis, label
any planar forest $F$ canonically in \emph{postorder}: visit the
components and the children from left to right, assigning a label
to a parent only after all its descendants have been visited.
Let $\mathrm{post}(F)\in\mathbf F_{\rm pl}[n]$ be this labelled
forest and put
\begin{equation}\label{eq:qCK_canonical_subspace}
 \mathcal H_q
   :=\operatorname{span}_{\Bbbk}
       \{\mathrm{post}(F):F\text{ a planar rooted forest}\}
       \ \subseteq\ \mathcal K_q(\mathbf F_{\rm pl}).
\end{equation}
The chosen labelling is compatible both with concatenation and
with the two planar forests obtained by an admissible cut
(restriction and standardization preserve postorder).  Hence
$\mathcal H_q$ is a graded $q$-Hopf subalgebra.  We identify
$\mathrm{post}(F)$ with the unlabelled forest $F$ from now on.
In this form, the same deformation may also be obtained by the
bosonic Fock functor applied to a weighted $q$-Hopf monoid of
planar forests.  This is the approach of
\cite[Section~13.3.5]{AM2010}, up to the ordering convention for
the two tensor factors.

\begin{theorem}\label{thm:qCK_Fock_interpolation}
Let $\mathcal H_q$ be defined by~\eqref{eq:qCK_canonical_subspace}.
For an admissible cut $c$ of a planar forest $F$, write
$P^c(F)$ and $R^c(F)$ for its pruned and remaining forests.
Define
\begin{equation}\label{eq:qCK_cut_statistic}
 \operatorname{inv}_F(c)
   :=\#\bigl\{(s,t)\in V(P^c(F))\times V(R^c(F)):
       \mathrm{post}_F(t)<\mathrm{post}_F(s)\bigr\}.
\end{equation}
Then the coproduct of $\mathcal H_q$ is
\begin{equation}\label{eq:qCK_weighted_cuts}
 \Delta_q(F)
  =\sum_{c\in\operatorname{Adm}(F)}
     q^{\operatorname{inv}_F(c)}
       P^c(F)\otimes R^c(F),
\end{equation}
where use the convention $0^0:=1$. For homogeneous tensors, set
\begin{equation}\label{eq:qCK_braided_tensor}
 (a\otimes b)\,\star_q\,(c\otimes d)
   :=q^{|b||c|}\,ac\otimes bd.
\end{equation}
The coproduct is coassociative and satisfies
\begin{equation}\label{eq:qCK_braided_compatibility}
 \begin{aligned}
 \Delta_q(1)&=1\otimes1,\\
 \Delta_q(FG)&=\Delta_q(F)\star_q\Delta_q(G),\\
 \Delta_q(B_+(F))&=B_+(F)\otimes1
      +(\mathrm{id}\otimes B_+)\Delta_q(F).
 \end{aligned}
\end{equation}
These three identities uniquely determine $\Delta_q$ on the free
associative operated algebra of planar rooted forests.  At $q=1$
one obtains the noncommutative Connes--Kreimer Hopf algebra;
at $q=0$ one obtains Foissy's unital infinitesimal Hopf algebra
with the pruned-first convention of
Subsection~\ref{subsec:Fock_Foissy_infin}.
\end{theorem}
\begin{proof}
A cut of $F$ corresponds precisely to a pruned-admissible subset
of its vertex set.  When the vertices carry their canonical
postorder labels, the Schubert exponent in
\eqref{eq:qCK_labelled_Fock} is exactly
$\operatorname{inv}_F(c)$.  Restriction and standardization yield
the canonical postorder labellings of $P^c(F)$ and $R^c(F)$;
this proves~\eqref{eq:qCK_weighted_cuts} and also shows that
$\mathcal H_q$ is preserved by the coproduct.  Concatenation
preserves postorder as well.  Since the deformed Fock functor
sends an ordinary Hopf monoid in species to a graded $q$-Hopf
algebra, its canonical subspace satisfies coassociativity and
the $q$-braided product identity in
\eqref{eq:qCK_braided_compatibility}.

\medbreak

We may also verify the grafting identity directly.  Apart from
the full cut, every cut of $B_+(F)$ leaves its new root in the
remaining forest and corresponds uniquely to a cut of $F$.
The new root is last in postorder, so it contributes no new
inversion to~\eqref{eq:qCK_cut_statistic}.  The right factor of
the corresponding term is $B_+(R^c(F))$.  The full cut accounts
for $B_+(F)\otimes1$.

\medbreak

Finally, every planar forest is a product of planar rooted trees,
and every tree can be written as $B_+(F)$ in a unique way for a smaller ordered forest.
Thus the unit, braided product rule, and cocycle recursion
determine the coproduct inductively on the number of vertices.
At $q=1$ every cut has coefficient~$1$.  At $q=0$ only cuts of
inversion number zero survive.  Their pruned vertices form an
initial interval of the postorder labelling, since
$\operatorname{sch}_n(S)=0$ exactly when $S=[|S|]$; conversely,
every postorder initial interval is pruned-admissible because
descendants precede their parents.  The result is the
left-admissible-cut coproduct of Foissy, with the stated choice
of left-to-right and tensor-factor conventions.
\end{proof}

\paragraph{The generalized infinitesimal identity.}
Writing $\Delta_q(F)=\sum F_{(1)}\otimes F_{(2)}$, the second
identity in~\eqref{eq:qCK_braided_compatibility} is equivalently to
\begin{equation}\label{eq:qCK_generalized_infinitesimal}
\begin{aligned}
\Delta_q(FG)={}&(F\otimes1)\Delta_q(G)
        +\Delta_q(F)(1\otimes G)-F\otimes G\\
&+\!\!\sum_{\substack{(F),(G)\\
                      |F_{(2)}|>0,\,|G_{(1)}|>0}}
 q^{|F_{(2)}||G_{(1)}|}
 F_{(1)}G_{(1)}\otimes F_{(2)}G_{(2)}.
\end{aligned}
\end{equation}
The extra summation vanishes at $q=0$; this is precisely the
unital infinitesimal compatibility.  At $q=1$ the sum restores
the ordinary multiplicativity of $\Delta_1$.

\medbreak

\begin{example}\label{ex:qCK_low_degrees}
Let $\bullet=B_+(1)$, let $\ell=B_+(\bullet)$ be the ladder
with two vertices, and let $t=B_+(\bullet\bullet)$ be the
three-vertex planar corolla.  The recursions give
\begin{align*}
 \Delta_q(\bullet)&=\bullet\otimes1+1\otimes\bullet,\\
 \Delta_q(\bullet\bullet)
  &=\bullet\bullet\otimes1
     +(1+q)\,\bullet\otimes\bullet
     +1\otimes\bullet\bullet,\\
 \Delta_q(\ell)
  &=\ell\otimes1+\bullet\otimes\bullet+1\otimes\ell,\\
 \Delta_q(t)
  &=t\otimes1+\bullet\bullet\otimes\bullet
    +(1+q)\,\bullet\otimes\ell+1\otimes t.
\end{align*}
In particular, the two nontrivial cuts pruning a single leaf
of the corolla contribute coefficients $1$ and $q$, respectively.
\end{example}

\medbreak

\subsection{Foissy's infinitesimal Hopf algebra}
\label{subsec:Fock_Foissy_infin}

The \emph{infinitesimal} coproduct $\Delta_0$ is another coproduct on the same
vector space of planar rooted forests. It differs from the ordinary
cut coproduct, and is the $q=0$ specialization of the family in
Subsection~\ref{subsec:qCK_deformed_Fock}. Its construction is particularly transparent from the $q=0$
specialization of the $q$-braided compatibility obtained by the
deformed Fock functor in the preceding subsection.

\medbreak

For homogeneous elements $x,y$ of degrees $m,n$, respectively,
consider the twisted interchange $\beta_q(x\otimes y)=
q^{mn}y\otimes x$, with $0^0=1$. At $q=0$, this interchange
vanishes when both degrees are positive. Consequently the induced
product on $\mathcal H\otimes\mathcal H$, for any connected graded
algebra $\mathcal H$ with augmentation $\varepsilon$, becomes
\begin{equation}\label{eq:Foissy_epsilon_product}
 (a\otimes b)\mathbin{\cdot_{\varepsilon}}(c\otimes d)
  =\varepsilon(b)\,ac\otimes d
  +\varepsilon(c)\,a\otimes bd
  -\varepsilon(b)\varepsilon(c)\,a\otimes d.
\end{equation}
In particular, requiring a coproduct
$\Delta_0:\mathcal H\to(\mathcal H\otimes\mathcal H,
\cdot_\varepsilon)$ to be an algebra morphism is equivalent to the
unital infinitesimal bialgebra identity
\begin{equation}\label{eq:Foissy_infinitesimal_identity}
  \Delta_0(FG)
  =(F\otimes1)\Delta_0(G)
   +\Delta_0(F)(1\otimes G)-F\otimes G.
\end{equation}
The multiplications on the right-hand side are the ordinary left
and right actions of $\mathcal H$ on $\mathcal H\otimes\mathcal H$.

\medbreak

On the free planar forest monoid, Proposition~\ref{prop:planar_universal_up}
constructs the unique coproduct that is multiplicative for
\eqref{eq:Foissy_epsilon_product} and satisfies the same cocycle
identity as in~\eqref{eq:NCK_Hochschild}. Equivalently, it is
uniquely determined by
\begin{equation}\label{eq:Foissy_recursive_delta}
  \begin{aligned}
  \Delta_0(1)&=1\otimes1,\\
  \Delta_0(FG)&=(F\otimes1)\Delta_0(G)
                +\Delta_0(F)(1\otimes G)-F\otimes G,\\
  \Delta_0(B_+(F))&=B_+(F)\otimes1+
                (\mathrm{id}\otimes B_+)\Delta_0(F).
  \end{aligned}
\end{equation}

There is an explicit combinatorial interpretation. Let $F$ be a
planar rooted forest with $n$ vertices. Following Foissy, order
its vertices as
$s_1>_{h,l}s_2>_{h,l}\cdots>_{h,l}s_n$, where the
total order $>_{h,l}$ combines the descendant relation with the
left-to-right order on incomparable vertices. Set
$S_k:=\{s_1,\ldots,s_k\}$ for $0\leq k\leq n$. Each $S_k$ is a
\Def{biideal}: it is closed upwards for the descendant order and
for the left order. Hence deleting the edges between $S_k$ and
its complement yields a planar pruned forest $F|_{S_k}$ and a
planar remaining forest $F|_{V(F)\setminus S_k}$. The formula is
\begin{equation}\label{eq:Foissy_left_cuts}
  \Delta_0(F)=\sum_{k=0}^{n}
       F|_{S_k}\otimes F|_{V(F)\setminus S_k}.
\end{equation}
Equivalently, the sum runs over \Def{left-admissible cuts} of $F$,
including the two trivial cuts; this is the description in
Foissy's infinitesimal Hopf algebra of planar forests (c.f. \cite{FoissyInf2009})

\medbreak

\begin{prop}\label{prop:recover_Foissy_infin}
The coproduct in~\eqref{eq:Foissy_left_cuts} is coassociative
and satisfies~\eqref{eq:Foissy_recursive_delta}. Together with
concatenation, it makes $\mathbf F_{\rm pl}$ a connected
$0$-bimonoid. Its bosonic Fock image
\begin{equation}\label{eq:recover_Foissy_Hinf}
  \mathcal H_{\rm inf}:=
   \mathcal K(\mathbf F_{\rm pl},\Delta_0)
\end{equation}
is the unital infinitesimal Hopf algebra of planar rooted forests
introduced by Foissy~\cite{FoissyInf2009,FoissyOperads2010}.
\end{prop}
\begin{proof}
For the coproduct~\eqref{eq:Foissy_left_cuts}, an iterated cut is
specified by two cut positions $0\leq i\leq j\leq n$ in the total
order on vertices. The two possible ways of grouping the resulting
three consecutive pieces give the same sum, proving coassociativity.
Every cut position in the concatenation $FG$ occurs either inside
$F$ or inside $G$. The position separating $F$ from $G$ is counted
in both cases, and is subtracted once. This proves
\eqref{eq:Foissy_infinitesimal_identity}. Finally, the new root in
$B_+(F)$ is the last vertex in the combined order. Every
nontrivial left-admissible cut of $B_+(F)$ therefore comes from a
left-admissible cut of $F$; this gives the last identity in
\eqref{eq:Foissy_recursive_delta}. These descriptions are invariant
under relabellings, so they define a $0$-bimonoid in species.
The bosonic Fock functor carries the $q=0$ compatibility to the
unital infinitesimal identity
\eqref{eq:Foissy_infinitesimal_identity}. Connectedness provides
the infinitesimal antipode, as in Foissy's construction.
\end{proof}

\medbreak
\begin{example}\label{ex:Foissy_compare_two_copros}
Let $\bullet$ be the one-vertex planar rooted tree. Already on the
two-component forest $\bullet\bullet$, the two planar coproducts
differ:
\begin{align*}
\Delta_1(\bullet\bullet)
 &=\bullet\bullet\otimes1
     +2\,\bullet\otimes\bullet
     +1\otimes\bullet\bullet,\\
\Delta_0(\bullet\bullet)
 &=\bullet\bullet\otimes1
     +\bullet\otimes\bullet
     +1\otimes\bullet\bullet.
\end{align*}
Moreover, for the planar corolla $t=B_+(\bullet\bullet)$,
\begin{align*}
\Delta_1(t)
 &=t\otimes1+1\otimes t
   +2\,\bullet\otimes B_+(\bullet)
   +\bullet\bullet\otimes\bullet,\\
\Delta_0(t)
 &=t\otimes1+1\otimes t
   +\bullet\otimes B_+(\bullet)
   +\bullet\bullet\otimes\bullet.
\end{align*}
Thus the difference is not caused by noncommutativity alone: it
comes from restricting ordinary admissible cuts to left-admissible
cuts. For the infinitesimal antipode one also has
$S_0(\bullet)=-\bullet$ and $S_0(\bullet\bullet)=0$.
\end{example}

\medbreak
\subsection{The three classical universal properties}
\label{subsec:three_CK_universal_properties}

We conclude by expressing the relation between the species
universal property and its three classical specializations.

\begin{theorem}\label{thm:three_CK_recoveries}
Let $\mathcal H$ be one of the following operated algebras:
\begin{enumerate}[label=(\roman*)]
\item $(\mathcal H_{\rm CK},B_+)$, the free commutative algebra on
      non-planar rooted trees;
\item $(\mathcal H_{\rm NCK},B_+)$, the free associative algebra on
      planar rooted trees, with the ordinary cut coproduct;
\item $(\mathcal H_{\rm inf},B_+)$, the same free associative algebra,
      with the infinitesimal coproduct~\eqref{eq:Foissy_left_cuts}.
\end{enumerate}
Let $A$ be a unital algebra and $L:A\to A$ a linear map, and assume
that $A$ is commutative in case~\textup{(i)}. There exists a unique
unital algebra morphism $\phi:\mathcal H\to A$ satisfying
\begin{equation}\label{eq:recover_univ_Connes_Foissy}
  \phi\circ B_+=L\circ\phi.
\end{equation}
If $A$ is a Hopf algebra in cases~\textup{(i)} and~\textup{(ii)},
or a unital infinitesimal Hopf algebra in
case~\textup{(iii)}, and $L$ satisfies
\begin{equation}\label{eq:Foissy_cocycle_target}
  \Delta_A(L(a))=L(a)\otimes1+
     (\mathrm{id}\otimes L)\Delta_A(a),
\end{equation}
then $\phi$ is also a coalgebra morphism, and hence a morphism of
the corresponding Hopf algebras.
\end{theorem}
\begin{proof}
The existence and uniqueness of the algebra map follow by induction
from the two recursion rules
\begin{equation*}
  \phi(1)=1_A,\qquad
  \phi(t_1\cdots t_r)=\phi(t_1)\cdots\phi(t_r),
  \qquad
  \phi(B_+(F))=L(\phi(F)).
\end{equation*}
Every tree is uniquely obtained by grafting its ordered forest of
branches (or its unordered forest in the commutative case), and
every forest is uniquely a product of its component trees, with the
appropriate commutativity convention. These are exactly the
bosonic Fock images of Theorem~\ref{Universality} and
Proposition~\ref{prop:planar_universal_up}.

For the coalgebra assertion, proceed by induction on the number of
vertices. The compatibility for products follows because
$\Delta_A$ is multiplicative for the ordinary tensor product in
cases~\textup{(i)}--\textup{(ii)}, or for
$\cdot_\varepsilon$ in case~\textup{(iii)}. The compatibility for a
single tree $B_+(F)$ follows by applying
\eqref{eq:Foissy_cocycle_target} and the common recursive identity
\eqref{eq:CK_Hochschild}, \eqref{eq:NCK_Hochschild}, or the last
identity in~\eqref{eq:Foissy_recursive_delta}, respectively.
The units are preserved by construction; the counits are
preserved by induction, using $\varepsilon_A\circ L=0$, which
follows from~\eqref{eq:Foissy_cocycle_target}. The antipodes are
then preserved by uniqueness of the convolution inverse. Thus
$\phi$ is a morphism of the corresponding Hopf algebras.
\end{proof}

The $q$-deformed interpolation of Theorem~\ref{thm:qCK_Fock_interpolation}
extends the planar cases of the preceding theorem. The construction can be summarized as follows:
\begin{equation*}
\begin{array}{ccl}
\mathcal K(\wE\circ\taa,\Delta_{\rm CK})
  &\cong& (\mathcal H_{\rm CK},\Delta_{\rm CK}),\\[3pt]
\mathcal K(\wL\circ\taa_{\rm pl},\Delta_1)
  &\cong& (\mathcal H_{\rm NCK},\Delta_1),\\[3pt]
\mathcal K_0(\wL\circ\taa_{\rm pl},\Delta_0)
  &\cong& (\mathcal H_{\rm inf},\Delta_0).
\end{array}
\end{equation*}
For each $q$, the planar part also admits the braided deformation
$(\mathcal H_q,\Delta_q)$ of Subsection~\ref{subsec:qCK_deformed_Fock}.

\medbreak

In this way, the Connes--Kreimer and Foissy universal properties
arise from a single species-level construction, namely the free (commutative or
ordered) forests with a branching up operator, followed by the
appropriate Fock functor and compatibility for the coproduct.


\bibliographystyle{plain}
\bibliography{Bibliography}

@book{AM2010,
  title={Monoidal functors, species and {H}opf algebras},
  author={Aguiar, Marcelo and Mahajan, Swapneel Arvind},
  volume={29},
  year={2010},
  publisher={American Mathematical Society Providence, RI}
}

@article{Joyal,
  title={Une th{\'e}orie combinatoire des s{\'e}ries formelles},
  author={Joyal, Andr{\'e}},
  journal={Advances in mathematics},
  volume={42},
  number={1},
  pages={1--82},
  year={1981},
  publisher={Elsevier}
}

@book{BLL,
  title={Combinatorial species and tree-like structures},
  author={Bergeron, Fran{\c{c}}ois and Labelle, Gilbert and Leroux, Pierre},
  series={67},
  year={1998},
  publisher={Cambridge University Press}
}

@article{FoissyBerlin,
  title={An introduction to {H}opf algebras of trees},
  author={Foissy, Loic},
  journal={preprint},
  year={2013}
}

@book{Mendez,
  title={Set operads in combinatorics and computer science},
  author={M{\'e}ndez, Miguel A.},
  year={2015},
  publisher={Springer}
}

@inproceedings{hairer2015geometric,
  title={Geometric versus non-geometric rough paths},
  author={Hairer, Martin and Kelly, David},
  booktitle={Annales de l'IHP Probabilit{\'e}s et statistiques},
  volume={51},
  pages={207--251},
  year={2015}
}

@article{AM2013,
  title={{H}opf monoids in the category of species},
  author={Aguiar, Marcelo and Mahajan, Swapneel},
  journal={{H}opf algebras and tensor categories},
  volume={585},
  pages={17--124},
  year={2013}
}

@incollection{connes1999hopf,
  title={{H}opf algebras, renormalization and noncommutative geometry},
  author={Connes, Alain and Kreimer, Dirk},
  booktitle={Quantum field theory: perspective and prospective},
  pages={59--109},
  year={1999},
  publisher={Springer}
}

@article{chapoton2007operads,
  title={Operads and algebraic combinatorics of trees.},
  author={Chapoton, Fr{\'e}d{\'e}ric},
  journal={S{\'e}minaire Lotharingien de Combinatoire [electronic only]},
  volume={58},
  pages={B58c--27},
  year={2007},
  publisher={Universit{\"a}t Wien, Fakult{\"a}t f{\"u}r Mathematik}
}

@article{livernet2006rigidity,
  title={A rigidity theorem for pre-{L}ie algebras},
  author={Livernet, Muriel},
  journal={Journal of Pure and Applied Algebra},
  volume={207},
  number={1},
  pages={1--18},
  year={2006},
  publisher={Elsevier}
}

@book{ML71,
	author = {Mac Lane, Saunders},
	title = {Categories for the {W}orking {M}athematician (Second ed.)},
	series = {Graduate Texts in Mathematics.},
	volume = {5},
	year = {1998},
	publisher = {Springer},
	ISBN = {0-387-98403-8},
}

@article{gucta2002symmetric,
  title={Symmetric {H}ilbert spaces arising from species of structures},
  author={Gu{\c{t}}a, M{\u{a}}d{\u{a}}lin and Maassen, Hans},
  journal={Mathematische Zeitschrift},
  volume={239},
  number={3},
  pages={477--513},
  year={2002},
  publisher={Springer}
}

@article{mendez2014antipode,
  title={An antipode formula for the natural {H}opf algebra of a set operad},
  author={M{\'e}ndez, Miguel A. and Liendo, Jean Carlos},
  journal={Advances in Applied Mathematics},
  volume={53},
  pages={112--140},
  year={2014},
  publisher={Elsevier}
}

@article{FoissyQ2003,
  author={Foissy, Lo{\"i}c},
  title={Quantifications des alg{\`e}bres de {H}opf d'arbres plans d{\'e}cor{\'e}s et lien avec les groupes quantiques},
  journal={Bulletin des Sciences Math{\'e}matiques},
  volume={127},
  number={6},
  pages={505--548},
  year={2003}
}

@article{FoissyInf2009,
  author={Foissy, Lo{\"i}c},
  title={The infinitesimal {H}opf algebra and the poset of planar forests},
  journal={Journal of Algebraic Combinatorics},
  volume={30},
  number={3},
  pages={277--309},
  year={2009}
}

@article{FoissyOperads2010,
  author={Foissy, Lo{\"i}c},
  title={The infinitesimal {H}opf algebra and the operads of planar forests},
  journal={International Mathematics Research Notices},
  volume={2010},
  number={3},
  pages={395--435},
  year={2010},
  doi={10.1093/imrn/rnp132}
}

@article{LRCofree,
  title={Algebres de Hopf colibres},
  author={Loday, Jean-Louis and Ronco, Mar{\i}a},
  journal={Comptes rendus. Math{\'e}matique},
  volume={337},
  number={3},
  pages={153--158},
  year={2003}
}

@inproceedings{Moerdijk,
title = "On the Connes-Kreimer construction of Hopf algebras",
author = "I. Moerdijk",
year = "2001",
language = "Undefined/Unknown",
pages = "311--321",
editor = "J. Greenlees",
booktitle = "Homotopy methods in algebraic topology",
publisher = "American Mathematical Society",
address = "United States",

}

@article{foissy2002algebres,
  title={Les algebres de Hopf des arbres enracin{\'e}s d{\'e}cor{\'e}s, I},
  author={Foissy, Lo{\i}c},
  journal={Bulletin des sciences mathematiques},
  volume={126},
  number={3},
  pages={193--239},
  year={2002},
  publisher={Elsevier}
}

 \end{document}